\documentclass[3p,serif]{elsarticle}

\usepackage{graphicx,float}
\usepackage[dvipsnames]{xcolor}

\usepackage{mathrsfs}
\usepackage{amsfonts}
\usepackage{amssymb}
\usepackage{amsmath,amsthm}
\usepackage{dsfont}
\usepackage{rotating}
\usepackage{multirow}
\usepackage{subfigure}

\DeclareMathAlphabet\mathbfcal{OMS}{cmsy}{b}{n}  % bold typing with \mathcal

\newcommand{\x}{\mathbf{x}}

\newtheorem{theorem}{Theorem}
\newtheorem{definition}[theorem]{Definition}
\newtheorem{proposition}[theorem]{Proposition}

\newtheorem{remark}[theorem]{Remark}
\newtheorem{corollary}[theorem]{Corollary}

\newcommand{\hh}{\textrm{h}}
\newcommand{\diag}{\text{diag}}	

\graphicspath{{Results/}}

\newcommand{\revp}[1]{\textcolor{black}{#1}}% REVIEWER 1
\newcommand{\revs}[1]{\textcolor{black}{#1}} % REVIEWER 2
\newcommand{\revt}[1]{\textcolor{black}{#1}} % REVIEWER 3

\newcommand{\revPO}[1]{\textcolor{black}{#1}}

\begin{document}
	
\begin{frontmatter}

	\journal{Journal of Computational Physics}

        \title{Fully well-balanced  entropy controlled \revp{discontinuous Galerkin spectral element method} for shallow water flows: \\global flux quadrature and cell entropy correction}
	
	\author[vit]{Yogiraj Mantri}
	\ead{yogiraj.mantri@vit.ac.in}
	\author[mn]{Philipp  $\ddot{\text{O}}$ffner}
	\ead{poeffner@uni-mainz.de}
	\author[Inria]{Mario Ricchiuto}%\corref{mycorrespondingauthor}}
	\ead{mario.ricchiuto@inria.fr}

	\address[Inria]{Inria, Univ. Bordeaux, CNRS, Bordeaux INP, IMB, UMR 5251, \\200 Avenue de la Vieille Tour, 33405 Talence cedex, France}
	\address[vit]{Vellore Institute of Technology, India}
	\address[mn]{Johannes Gutenberg-University, Mainz, Germany}

	%-------------------------------------------
	% ABSTRACT
	%
	\begin{abstract}

\revs{
We present a novel approach for solving the shallow water equations using a discontinuous Galerkin spectral element method. 
%Our method incorporates a discrete fully well-balanced property as well as  control of entropy production. 
The method we propose has three main features. First, it enjoys  a discrete well-balanced property, in a spirit similar to the one of e.g. \cite{Castro2020}. 
As in the reference, our scheme  does not require  any a-priori knowledge of the steady equilibrium, moreover it does not involve the explicit
solution  of any local auxiliary problem to approximate such equilibrium. The scheme is also 
arbitrarily high order, and verifies a continuous in time cell entropy equality. The latter becomes an inequality as soon as  additional dissipation is added to the method.
The method is constructed starting from a global flux approach  in which an additional flux term is constructed as the primitive of the source.
We show that, in the context of nodal spectral finite elements, this  can be translated into 
 a simple modification of the integral of the source term. We prove that, when using Gauss-Lobatto nodal finite elements
  this modified integration is equivalent  at steady state to a high order Gauss collocation method applied to an ODE
  for the flux.    \revp{This  method is superconvergent at the collocation points, thus providing    a discrete well-balanced  property
  very similar in  spirit  to the one proposed in \cite{Castro2020}, albeit not needing   the  explicit computation of a local approximation of    the steady state. }
To control the entropy production, we introduce artificial viscosity corrections at the cell level and incorporate them into the scheme. We provide theoretical and numerical characterizations of the accuracy and equilibrium preservation of these corrections. Through extensive numerical benchmarking, we validate our theoretical predictions, 
with considerable  improvements in accuracy for steady states, as well as enhanced robustness for more complex scenarios.
}

%
%\revo{where the well-balancing property is approximately rather than exactly fulfilled meaning it is satisfied in terms of the order of the scheme.}
	\end{abstract}	
	%------------------------------------------
	% KEY WORDS
	%
	%==========================================================================================
\end{frontmatter}

%%%%%%%%%%%%%%%%%%%%%%%%%%%%%%%%%%%%%%%%%%%%%%%%%%%%%%%%%%%%%%%%%%%%%%%%%%%%%%%%%%%%%%%%%%%%%%

%--------- SECTION ---------------------------------------------------------------------------
\section{Introduction}

Hyperbolic balance laws play a fundamental role in  various phenomena in natural science and engineering. 
A system of balance law is given by the following form 
\begin{equation}
	\partial_t U+\nabla\cdot F(U)=  S(U;\varphi(x)) \label{balance_law},
\end{equation}
where $U$ \revs{contains the} conserved variables,  $F$ is  the flux function, and $S$ denotes  the source term, which depends on the solution as well as on some external 
data (e.g. bathymetry, friction coefficient map, etc.) dependent on space. 
The numerical approximation of \eqref{balance_law} is  a very active research topic, with a lot of fundamental
contributions  developed to models such as 
the Euler equations with  gravity \cite{chertock2018well, GCD18, thomann2020speed,wilcox2022entropy}, 
or shallow water (SW) equations with various sources accounting for topography variations, friction and/or Coriolis forces
in Cartesian  or curvilinear coordinates \cite{berthon:hal-00956799, arpaia:hal-02422335, chertock2022well, ciallella2022arbitrary,  MICHELDANSAC2017115,arpaia:hal-03421078,carlino:hal-03850195}.
Even in one space dimension, the challenge of devising  well-balanced numerical approximations    
agnostic of the form of the steady state is still open. 
 There is already quite a large literature on the subject,
with several different approaches \revp{to manage this} issue. 
This paper focuses on the shallow water equations \revPO{in Cartesian coordinates},  including  all the effects mentioned \revPO{above}.\\

The  source term in   \eqref{balance_law}  leads to  a  rich  set of possible solutions, and in particular 
to a large number of different forms of steady equilibria between the source terms and the flux derivatives.
\revPO{Many of such equilibria have some interest in themselves. Many
 physical applications } involve   small perturbations of such equilibria. 
The ability of  a numerical method % to preserve  or at least
 to resolve with enhanced accuracy such steady states is a unanimously acclaimed design criterion usually referred to as \revPO{well-balanced or
 full well-balanced depending on whether the property  applies to a specific equilibrium, or to all  steady states. }
 %
% 
% what one knows of the equilibrium: does one need to know already a specific solution; does one need to know a family of solutions (from some relation, or invariant set) ; is the method agnostic of the exact solution and still well balanced for generic steady states;
%? is each method approximately or exactly well-balanced (according to [4] for instance)? if approxi- mately: in which sense ?
%? is there any linear or non-linear auxiliary problem to be solved ?
%? is the scheme high-order accurate
\revp{Giving a full review of the subject is way beyond our scope.   Questions relevant to the work of this paper   are the following:
\begin{itemize}
\item  what do we know of the equilibrium we wish to preserve: do we know it  already in  (some) explicit  form, 
do we   know if it belongs to a family of solutions verifying  some algebraic constraint, or  do we want  our method to be 
agnostic of  it, and still be well-balanced  in some sense;
\item do we seek to preserve the analytical steady state or some  approximation of it, and which one;
\item can we use some   auxiliary problem to improve  the well-balanced character of the scheme;
\item what is the accuracy aimed for, is the order  of the method arbitrary.
\end{itemize}
Note that we do not address the property of genuinely multidimensional  well-balanced, which is related to the notion of the preservation of solenoidal involutions.
We   also leave out   extensions including dry areas, for which the only relevant steady state is the lake at rest, at least in one space dimension.}

\revp{Concerning the other aspects, some have been discussed widely in literature.
For example, when steady (or even unsteady) equilibria are  explicitly known one can use a simple    and efficient idea
due to \cite{doi:10.1137/060674879}, which consists in evolving a discrete error with respect to the given equilibrium.
This boils down to removing from the discrete equations of any scheme  the discrete expression
 corresponding to the application of the scheme  itself to the given  solution. 
 This idea has been adapted to many discretization approaches and applied
 to several systems of equations ranging from the Euler equations with gravity
 \cite{GCD18,BERBERICH2021104858}, to the shallow water on manifolds \cite{carlino:hal-03850195},
 to the MHD equations \cite{BERBERICH2021104858,birke2023wellbalanced},  to hyperbolic reformulations of the Einstein equations of relativity \cite{dumbser2023wellbalanced},  to cite a few. This approach can be applied to study perturbations
 of a given equilibria with any scheme,  and to any order of accuracy.}\\
\revp{%A more complex case is the one in which one does no know explicitly the  relevant equilbrium,  but it is} characterized by a set of \revp{(generally nonlinear) algebraic 
%relations defining invariants which are constant.
\revp{A more intricate scenario arises when the relevant equilibrium is not explicitly known, but it can be characterized by a set of (generally nonlinear) algebraic relations that define constant invariants.}
 %In other words,  in this case the equilibrium can be defined by the set of relations
%$V=V_0$, where $V$ is a complete set of  variables which, ideally, can be used to obtain the conserved quantities $U$. 
\revp{To put it differently, in this scenario, the equilibrium state can be described through the collection of relationships expressed as $V=V_0$, where $V$ represents a comprehensive set of variables that, ideally, can be leveraged to derive the conserved quantities $U$.}
In this case, one tries to devise a      discretization} consistent  in some way with an approximation of the modified form of \eqref{balance_law} reading (in one space dimension),
\begin{equation}
	\partial_tU+ A_V(U,\varphi(x)) \partial_x V =  0  \label{balance_law-1}.
\end{equation}
\revp{This idea has been adapted within  a variety of methods  going from finite volume to finite elements,
and  embedded  either in the   polynomial approximation, and/or
in the definitions of the  discrete divergence and of the numerical flux and source  terms  \cite{10.1137/S1064827503431090,duenas11,R11,R15,XING2014536,cheng16,MICHELDANSAC2017115,berthon:hal-00956799,GOSSE2000135}.
While not requiring the full a-priori knowledge of the steady state, one still requires  the existence and the  knowledge of the invariant set $V(U)$. In this situation also,
discretizations with very high order of accuracy have been proposed in literature.}

%In the finite volume setting to construct well-balanced schemes is 
%obtained by enhancing the Riemann solvers providing the numerical fluxes with additional contact-like steady waves
%across which a cell integrated balance flux difference-source can be imposed.

\revp{When  nothing is known of the steady state,  two interesting approaches are considered here.
The first,  introduced in \cite{Castro2020},    is a discrete generalization of the idea of \cite{doi:10.1137/060674879}.
As the latter, it is quite general   and it     has  been  applied across various numerical frameworks such as finite elements, finite differences and finite volume 
\cite{Castro2020,math9151799,math10010015,GOMEZBUENO2021125820,GOMEZBUENO202318,macca22}. 
As in      \cite{doi:10.1137/060674879}   the  idea  is to evolve a perturbation with respect to the equilibrium solution.
However, when the latter is unknown  the authors of  \cite{Castro2020} propose to solve locally the auxiliary Cauchy problem 
$U'(x) = A^{-1}(U)S(U,x)$, where $A(U)$ is the flux Jacobian $\partial_UF$  to obtain an estimate of a local admissible steady state.
The admissibility is linked to the match with the mesh data.  }
%be implemented in conjunction w any numerical setting 
%(finite elements, finite difference, finite volume, etc.) \revPO{This makes no sense: 20 was after 18. I guess you meant something like this: \\
%The initial one is a discretized extension of the concept presented in \cite{doi:10.1137/060674879}, which was ultimately introduced in \cite{Castro2020}. Similar to the fundamental methodology in \cite{doi:10.1137/060674879}, this approach can be. \\
%Additionally, I do not understand the following sentences and also we should shorten them make them more easy to read.... here we have to discuss... and where does number two start....
%  }. The main idea is to define a general auxiliary problem 
%which is used to obtain locally an enhanced approximation of a  steady state compatible with the mesh data.
%%This approximation needs to be very accurate, anyways more accurate of the underlying discretization. 
%This approximation needs to be very accurate, anyways more accurate than  the underlying discretization. 
%For balance laws, the main idea   is to use as auxiliary problem  the Cauchy problem $U'(x) = A^{-1}S(U,x)$,
%where $A$ is the flux Jacobian $\partial_UF$. }
%This is quite a flexible setting and has been applied to a variety of schemes and models,
%and with different strategies to obtain the solution of the (local) Cauchy problem
%\cite{Castro2020,math9151799,math10010015,GOMEZBUENO2021125820,GOMEZBUENO202318}.
  \revp{This method is not an exactly well-balanced, however it has also a clear definition of the notion of  discrete equilibrium,
  associated to the solution of the local Cauchy problem.   Its drawback is that  it requires the explicit resolution of the latter as 
an  auxiliary problem.}\\
%
%\revPO{I would split the follwoing sentence in two, something like: 
\revp{A second approach   which is  also agnostic of the   steady equilibrium, is provided by the class of methods  known as global flux schemes. These schemes, initially introduced in the research by \cite{GASCON2001261}, were initially employed to develop nonlinear shock-capturing techniques for balance laws in \cite{CASELLES200916} and \cite{donat11}. 
More recently, they have been utilized as  means to achieve fully well-balanced methods, as seen in works like \cite{chertock2018well, CCHKT19, MantriNoelle2021}.
}
%Lastly, \revp{within the methods agnostic of the steady state,} we mention the family of so-called global flux schemes, which 
%\revp{date  back to the work of } \cite{GASCON2001261},
%and have been applied first to design non-linear shock capturing methods for balance laws in \cite{CASELLES200916,donat11}, and more recently
%as a path to obtain fully well-balanced methods (see e.g. \cite{chertock2018well,CCHKT19,MantriNoelle2021}).  
This approach relies on the observation that
(up to a constant) the non-local operator
\begin{equation}
	G(U,x) = F(U(x))   - \int\limits_{x_0}^xS(U(s),\varphi(s))ds  \label{balance_law-2}
\end{equation}
provides a natural invariant for the steady equations. So recasting \eqref{balance_law} as 
\begin{equation}
	\partial_tU+\partial_x G(U,x) =0 \label{balance_law-3},
\end{equation}
provides a reasonable path to obtain fully well-balanced schemes. As shown in \cite{Abgrall2022}, this approach also has relations
with flux difference splitting  and residual distribution methods consistently embedding the source integral in the splitting
and can thus be related to many classical versions of this idea \cite{Roe87,bv94,VAZQUEZCENDON1999497,pares_castro_2004}.
\revp{It  has clear relations with methods which define the discrete source term as  the derivative 
of a steady equilibrium flux, e.g.  \cite{PARES2021109880}.}
While providing very interesting results for many different models and solutions, one of the drawbacks of this approach is that, 
differently from the \revp{previous one,}
% correction method with explicit evaluation of the  solution of the steady Cauchy problem, there is so far
there is no
%no work on the 
precise characterization of the meaning of the steady solution. For this reason in \cite{macca22},  the authors
have chosen to combine the use of \eqref{balance_law-3} with the correction \revp{approach of \cite{Castro2020} 
which allows a more clear control to the notion of} %to have a more precise definition of 
well-balancing.\\

In addition to the above constraints, system  \eqref{balance_law} is usually endowed with 
an entropy pair  $(\eta(U),F_{\eta}(U))$ verifying the additional constraint 
\cite{dafermos20210hyperbolic, harten1983symmetric} 
\begin{equation}
	\partial_t\eta+\nabla\cdot F_\eta(U)   \stackrel{(\leq)}{=}  S_{\eta}(U;\varphi(x)) \label{entropy_balance_law},
\end{equation}
where $\eta=\eta(U)$ denotes the mathematical entropy (a convex function), $F_\eta(U)$ is the entropy flux, and  $S_{\eta}(U;\varphi(x))$ represents dissipation/production term. 
In \eqref{entropy_balance_law}   the equality holds for smooth solutions while weak admissible solutions  are characterized by   the inequality.  A degree of control on the production of entropy is thus a desirable feature in numerical schemes.  
A numerical method is called entropy conservative if it fulfills the equality in \eqref{entropy_balance_law} 
and  entropy dissipative if it ensures  the inequality  for one specific entropy   \cite{tadmor2003entropy}. 
There exist many different techniques which have been applied to obtain entropy conservation (dissipation), e.g. 
artificial viscosity and correction techniques \cite{abgrall2018general,abgrall2019analysis, abgrall2022reinterpretation, guermond2018second}, the summation-by-parts framework combined with entropy conservative (EC) two-points fluxes (flux differencing) \cite{chenreview, fisher2013discretely, renac2019entropy}, or multi-point approaches via  with  combination of EC fluxes \cite{fjordholm2012arbitrarily, lefloch2002fully}.  Some of this work has been generalized to embed some notion of  well-balanced.  Here, we refer to  \cite{fjordholm2011wellbalanced,ranocha2017shallow,  gassner2015shallow} for the literature on shallow water equations.\\ 

In this paper, we propose a high order \revp{discontinuous Galerkin spectral element method (DGSEM)} formulation for balance laws which embeds a fully discrete general well-balanced criterion agnostic of the exact steady state. The proposed construction exploits the idea of a global flux formulation to infer an ad-hoc quadrature strategy called here global flux quadrature. This is then used  to establish a one-to-one correspondence
between the discretization of the local steady Cauchy problem, and a discretization of the non-local integral operator underlying the
definition of the global flux. If the  discretizations exploit the same data on the same stencil, the steady solution can 
be obtained indifferently by means of one of the two methods.   This equivalence allows to construct balanced schemes without explicit knowledge of the steady state, and  without the need of solving explicitly the local Cauchy problem.  In this paper, we use a Gauss-Lobatto DGSEM setting which allows a natural connection to continuous collocation methods for integral equations. 
Thus we are able to fully characterize the discrete steady solution, and moreover provide  a superconvergence result \revp{at the collocation points}. The notion of entropy control is also included in the construction via appropriately designed artificial viscosity corrections at the cell level
following \cite{abgrall2018general,abgrall2022reinterpretation, abgrall2022relaxation,gaburro2023high, offner2020stability}.
 The accuracy and equilibrium preservation of these corrections are characterized theoretically and numerically. 
 
 The paper is organized as follows: The main notation for the shallow water equations are recalled in section \S2, where we also
 recall  several families of steady equilibria, depending on the terms included in the sources, and the definition of the mathematical entropy.
 In section \S3, we  recall the basics of  the DGSEM approach underlying the paper, introducing  some of its properties and the notation necessary for the remainder of the paper.
 The main idea of the paper is discussed in section \S4 devoted to the explicit derivation of the global flux quadrature approach, and to the characterization of its properties
for the particular case of the Gauss-Lobatto DGSEM method. The entropy correction is studied in section \S5, while a simple multidimensional extension is proposed in section \S6.
Section \S7 provides a thorough verification of the theoretical expectations, as well as some applications to challenging cases  and to multidimensional problems showing the
great potential of the approach proposed. Some conclusive remarks and future perspectives are drawn in section \S8.

\section{Shallow water equations}\label{se_SW}

This paper focuses on the shallow water (SW) equations, widely used in geophysical applications. 
Despite their relative simplicity \revp{with respect to} other systems, they  offer an excellent test bed for well-balanced methods due to the
variety of source terms they embed. The two dimensional form of the  system reads 
\begin{equation}\label{SW0}
	\partial_t\left(\begin{array}{c} h \\ hu \\ hv   \end{array}\right)+\partial_x \left(\begin{array}{c} hu \\ hu^2 + p(h) \\ huv   \end{array}\right) +\partial_y \left(\begin{array}{c} hv \\ huv \\ hv^2+p(h)   \end{array}\right) 
	=  -h \left(\begin{array}{c}  0\\ \partial_x\varphi +c_f u +  \omega v\\ \partial_y\varphi  + c_f v - \omega u\end{array}\right) ,
\end{equation}
where $h$  denotes the water depth, $\vec{\mathsf{v}}=(u,v)^T$ is the  horizontal velocity, $ p = gh^2/2$ is hydrostatic pressure with $g$ the gravity acceleration, 
$\varphi = g b$ denotes the gravitational potential with bottom topography $b(x,y)$, 
$c_f=c_f(h,\vec{\mathsf{v}})$ is friction coefficient, and  
$\omega$ denotes the  Coriolis coefficient.  It is customary and useful to introduce the additional variables  $\zeta= h+b$, 
representing  the free surface elevation, and the total energy density $E= g\zeta+k$, with \revp{$k=u^2/2 + v^2/2$ being the kinetic energy.}   
\revPO{In this work we will only consider continuous bathymetry, so the data $\varphi$ in the source is also assumed to be continuous.
Preliminary extensions of ideas similar to those discussed here to the more general case are considered in \cite{ctr23}.}\\
An interesting variant of  \eqref{SW0} is the pseudo-one dimensional rotating SW system proposed in \cite{MR2460785}
\begin{equation}
	\partial_t\left(\begin{array}{c} h \\ hu \\ hv   \end{array}\right)+\partial_x \left(\begin{array}{c} hu \\ hu^2 + p(h) \\ huv   \end{array}\right)  
	=  -h \left(\begin{array}{c}  0\\ \partial_x\varphi +c_f u +  \omega v\\    - \omega u\end{array}\right)   \label{SW1}.
	\vspace{0.5cm}
\end{equation}
This is a one dimensional system which has all the richness of the multidimensional one in terms of sources and steady states, as we will  shortly discuss. \\

The shallow water equations are endowed with a convex   entropy pair $(\eta, F_{\eta})$
given by
 \begin{equation}\label{SW-eta0}
\eta  = p(h) + hk\;,\quad
  F_{\eta}  = hu\, ( gh + k  ). 
 \vspace{0.25cm}
\end{equation}
 In the case of   \eqref{SW1} the associated balance law reads,  
  \begin{equation}\label{SW-eta1}
\partial_t \eta   + \partial_x  F_{\eta} \stackrel{(\leq)}{=} -  hu\partial_x\varphi - \underbrace{c_f hu^2}_{\mathcal{D}_f}
\end{equation}
with $ \mathcal{D}_{f} \ge 0$  implicitly defined above representing the dissipation due to friction. 
For time independent potentials, the form of the entropy production term allows to introduce a 
total entropy pair $(\eta_\varphi, F_{\eta_{\varphi}})$
which  is given for  \eqref{SW1} via 
 \begin{equation}\label{SW-eta2}
\eta_{\varphi} = p(h) + hk + h\varphi\;,\quad
  F_{\eta_{\varphi}} = hu\, ( gh + k  + \varphi ) =hu\, ( g\zeta + k   )
 \vspace{0.25cm}
\end{equation}
which satisfies  the simpler   balance 
 \begin{equation}\label{SW-eta3}
\partial_t \eta_{\varphi} + \partial_x  F_{\eta_{\varphi}} \stackrel{(\leq)}{=} - \mathcal{D}_f \le 0
\end{equation}
which reduces to a special conservation law for smooth solutions and in absence of friction. \revp{Note that no dissipation or production term 
are associated to the Coriolis terms, as the 
  $x$ and $y$   contributions to the energy balance   cancel identically when dotting the momentum equations by the velocity $(u,v)^T$.}

\subsection{%Incomplete zoology of steady equilibria
\revs{Partial taxonomy of steady equilibria}}
We recall here some of the classical equilibria of system \eqref{SW1} and \eqref{SW0}, which will be used in the process of numerical validation.
We limit the description  to five types of solutions, involving three combinations of sources, but many more can be imagined.

\paragraph{Frictionless one dimensional equilibria} This is the most classical case under consideration. It  involves two families
of steady equilibria. For smooth solutions, the most general form of these equilibria is
defined by the relations
%can characterize it by the two conditionsIn general, it  can be  obtained by considering smooth steady solutions, for which we can combine the first equation in \eqref{SW1} with the steady limit of \eqref{SW-eta3}
%to obtain the invariants
 \begin{equation}\label{steady0}
 h(x)u(x)=q_0\,,\;\; g\zeta(x) + k(x) = E_0,
\end{equation}
to obtain the invariants $q_0,E_0$,
with $q_0$ and $E_0$  being the given values of the volume flux and total energy density. Of course this is also a special
 solution of  \eqref{SW0} if $v=0$, otherwise a similar one can be defined in the frame of reference with the $x$ axis aligned to the velocity.
 The   moving equilibrium \eqref{steady0} can  be determined analytically by solving
 a cubic equation, given   $q_0$ and $E_0$  and the bathymetry $b(x)$ (see e.g. \cite{nxs07}).

A very well known particular case is the so-called \emph{lake at rest} state defined by
 \begin{equation}\label{steady1}
 u(x)=0\,,\;\; \zeta(x)   =\zeta_0.
\end{equation}
This  is the most widely used  analytical state  to  construct well-balanced methods, and first historically studied \cite{bv94}.

Note that at the analytical level, for smooth cases anyways, we can replace \eqref{steady0} by a global flux definition reading:
 \begin{equation}\label{steady2}
 hu=q_0\,,\;\;  Q:= hu^2 + p(h) +\int_{x_0}^{x}gh \partial_x b = Q_0.
\end{equation}
For smooth cases the two definitions coincide. 

\paragraph{Frictionless pseudo-one dimensional equilibria with Coriolis effects} This is a more general case of the previous one, including transverse 
velocity and Coriolis effects. \revs{ It also involves} two families of solutions, depending on  whether $u=0$ or not.
When $u=0$ we have   solutions which are defined by 
 \begin{equation}\label{steady4}
 hu=0\,,\;\; \zeta   = \zeta_0 - \dfrac{1}{g} \int_{x_0}^x \omega v(s)ds. 
\end{equation}
In this case the transverse velocity $v$ is essentially a free parameter. 
This is some sort of generalized analytical lake at rest state, accounting for movement in directions transverse  to those of the main bathymetric variations.

\revPO{
  For the moving case, combining the first and last equation in \eqref{SW1}, with the steady limit of \eqref{SW-eta3} one can show that  \eqref{steady3} holds, provided $u\ne 0$. 
 } 
 \begin{equation}\label{steady3}
 hu=q_0\,,\;\; g\zeta + k = E_0 \,,\;\;
 v - \omega x = v_0.
\end{equation}
Equation \eqref{steady3} needs to be solved analytically  given values of $q_0$, $E_0$ and $v_0$,  and  the bathymetry $b(x)$. 
Particular  smooth solutions  can also be obtained by defining $b(x)$ for given variations of $h$ and $u$. An example will be presented in the results section.
As before, for smooth cases we can equivalently provide a  global flux version of  \eqref{steady3}, which  reads
 \begin{equation}\label{steady5}
hu=q_0\,,\;\;
Q:= hu^2 + p(h) +\int_{x_0}^{x}(gh \partial_x b +\omega h v) = Q_0 \,,\;\;
 V := huv -\int_{x_0}^{x} \omega hu  = V_0.
\end{equation}
For smooth cases the solutions defined by \eqref{steady5} or \eqref{steady3} are equivalent. \revp{A detailed discussion for such problems and more can be found in recent papers \cite{Desveaux2022fully, chertock2018wellbalancing}.}

\paragraph{One dimensional equilibria with friction} This is a  moving equilibrium with friction and slope variations  but no Coriolis force. 
It cannot in general  be defined  analytically. It is characterized by the  relations ($v=0$ since it is one dimensional)
 \begin{equation}\label{steady6}
hu = q_0\,,\;\;   g\zeta  +  k  =  E_0 - \int_{x_0}^xc_f u \,ds .
\end{equation}
A trivial particular case for it is $b'(x)=\xi_0=$const, which admits as exact solution a special state $h=h_0$ and $u=u_0$
verifying   the algebraic relations 
 \begin{equation}\label{steady7}
h_0u_0 = q_0\,, \;\;  c_f(h_0,u_0)u_0 = -g\xi_0.
\end{equation}
These can be solved for given values of $q_0$ and of the slope $ \xi_0$, and the friction law (see e.g. \cite{R15}).
More general examples obtained integrating \eqref{steady6} are proposed e.g. in \cite{MICHELDANSAC2017115}.
 
 The global flux definition of these states differs a little bit from \eqref{steady6} and is given by
 \begin{equation}\label{steady8}
 hu=q_0\,,\;\;  Q:= hu^2 + p(h) +\int_{x_0}^{x}(gh \partial_x b + c_f(h,u)hu) = Q_0.
\end{equation}

\section{\revp{DGSEM} discretization: main notation and setting}\label{sec:dgsem}

\subsection{Main notation}\label{sec:notation}
\begin{figure}[h]
\begin{center}
\includegraphics[width=0.3\textwidth]{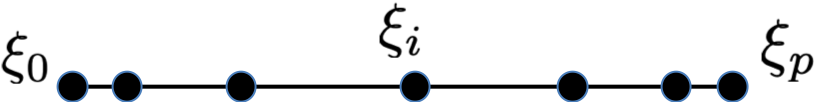}\includegraphics[width=0.3\textwidth]{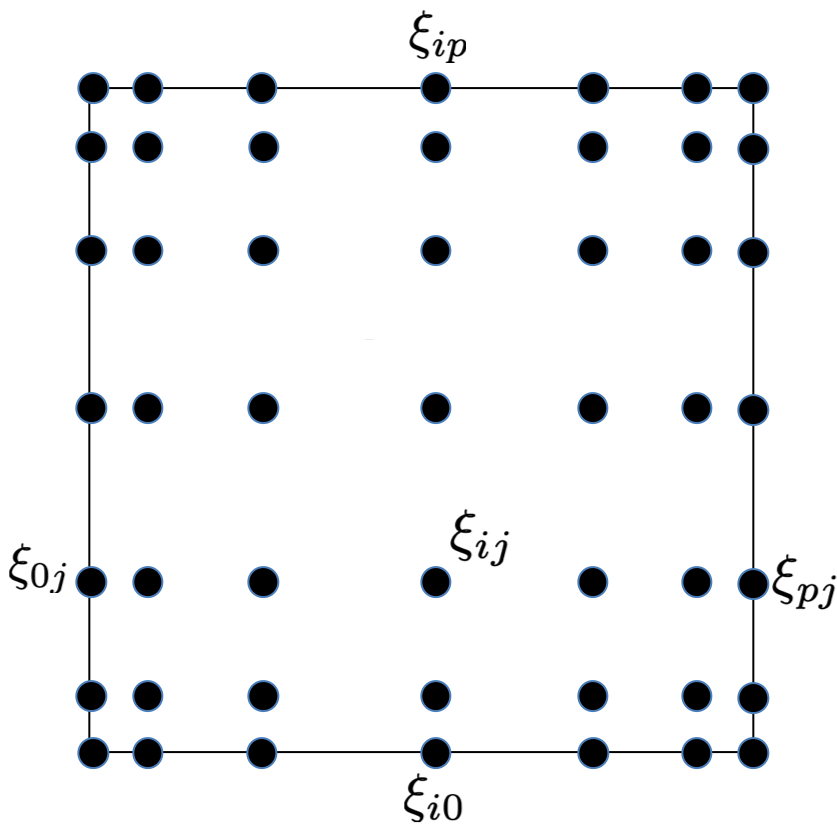}
\caption{\label{fig:elements}One and two dimensional Gauss-Lobatto \revp{nodes. The polynomial approximation is of degree $p$ in each direction separately for the $p+1$ GL points.}}
\end{center}
\end{figure}

We consider a tessellation of the $d$-dimensional spatial domain $\Omega$ in non overlapping elements $K$, obtained as tensor products of 1d elements. In other words \revPO{in two dimensions} $K=K_x\times K_y$
%\times K_z$,
with  $K_j$ a segment of \revt{width $\hh$} in the  direction $j$ (cf figure \ref{fig:elements}).
As it is classical, we consider in each direction a linear map \revPO{onto the classical unit reference element}  $x(\xi): K_x\mapsto[0,1]$,  and  similarly in the other directions. 
In the multidimensional case we will denote by $\{\xi_{j}\}_{j=1,d}$ the 
reference coordinates. 
On each 1d  reference element we consider a  standard Gauss-Lobatto (GL) collocated finite element approximation spanned by 
$\{ \phi_i(\xi) \}_{i=0,p}$  degree $p$ one dimensional Lagrange basis functions corresponding to the  $p+1$ GL points 
 $\{ \xi_i \}_{i=0,p}$.  For a given function $u$ we thus set in 1d
 \begin{equation}\label{eq:uh}
 u_{\hh} :=  \sum\limits_{j=0}^p\phi_j(x(\xi)) u_j
 \end{equation}
 while in 2d we have the usual tensor product representation
 \begin{equation}\label{eq:uh1}
 u_{\hh} := \sum\limits_{\revt{i,j=0}}^p\phi_i(y(\zeta))\phi_j(x(\xi)) u_{ij}
 \end{equation}
interpolating the solution on  a grid of Gauss-Lobatto points $\{\xi_{ij}\}_{i,j=0,p}$  as in figure \ref{fig:elements}, with $\xi_{ij}=\xi_i\zeta_j$.

%The temporal domain $[0,T]$ is classically discretized by temporal slabs $[t^n, t^{n+1}]$, with $\Delta t^n =t^{n+1} -t^n $, and  \revs{$\Delta t =\max\limits_n\Delta t^n$}.
%The superscript $^n$  will be dropped when no confusion is possible.

      \begin{remark}[Nodal approximation   for the data]   \label{remark_nodal_data}
{\it \revp{As already mentioned in section \S2,   we only  consider continuous data, and in particular for continuous bathymetry.  
Since   for   continuous   $\phi(x)$ the discrete spectral element
projection  $\phi_{\hh}$ will boil  down to standard nodal interpolation
on Gauss-Lobatto points.  This leads to  an approximation with \emph{no jumps at element boundaries} despite of the locality of the overall approximation,
and of the fact that  jumps may occur for  the solution $U_{\hh}$. Initial work  including discontinuous solutions as well as  data (e.g. bathymetry) with jumps is discussed in \cite{ctr23}. The extension to the 
current DGSEM framework is ongoing.}}
\end{remark}

\subsection{DGSEM for 1d conservation laws}
We consider  for the moment the approximation of solutions of
\begin{equation}\label{eq:dgsem0}
\partial_t U + \partial_xF(U) =0.
\end{equation}
On \revs{a generic element} $K$, we start from the standard DG variational form 
\begin{equation}\label{eq:dgsem1}
|K|\int\limits_0^1\phi_i(\xi) \partial_t U_{\hh} - \int\limits_0^1\partial_{\xi}\phi_i(\xi) F_{\hh}  +   
  (\phi_i \hat F_{\hh}(U_{\hh},U_{\hh}^+))_{\xi=1}-  (\phi_i \hat F_{\hh}(U_{\hh},U_{\hh}^+))_{\xi=0} =0.%\label{DGvarform}
\end{equation}
Following the classical spectral element approach, the quadrature used is the one associated to the GL nodes also used for the
nodal approximation. The resulting semi-discrete equations can be written in matrix form as \cite{Kopriva2010,heshavenWarburton},
\begin{equation}\label{eq:dgsem2}
    \dfrac{\textrm{d}  \mathbf{U}}{\textrm{d}t}  -    \widetilde{D}_x^T \mathbf{F}
% - \widetilde{D}_y^T \tilde{\mathbf{F}}^i_{y}
+    \mathcal{M}^{-1}
 \mathcal{B} \widehat{\mathbf{F}} =0,%+    \mathcal{M}^{-1} \mathcal{B} \widehat{\mathbf{F}}_y =0
\end{equation}
where  $\mathcal{M}=\diag(\{ w_i\}_{i=0,p})$  with   \revt{ $w_i $ the Gauss-Lobatto quadrature weights\footnote{\revp{Note that the GL quadrature is of degree $2p-1$ for $p+1$ points. Therefore, the last entry of the mass matrix is lumped, cf. \cite{Gassner2013, ranocha2016}.} }, 
$ \widetilde{D}_x = \mathcal{M} D_x \mathcal{M}^{-1}$   with $(D_x)_{ij} =\partial_{\xi}\phi_j(\xi_i)$, and with  }
 $ \mathcal{B} = \diag(-1, \dots,1)$ the matrix sampling boundary values. The arrays
 $ \mathbf{U},\,  \mathbf{F},\,$ and  $\widehat{\mathbf{F}}$ respectively contain nodal values of the solution, the flux and the numerical flux.

\subsection{\revp{Summation-by-parts}, strong form, and residual distribution} 
 The nodal DGSEM method is known to enjoy a summation-by-parts \revp{(SBP) property mimicking integration-by-parts on the discrete level. In terms of the operator notation \revt{focusing on a reference element}, this means that
 \begin{equation*}
 \mathcal{M} D_x +  D_x^T\mathcal{M} = \mathcal{B} \Longleftrightarrow   D_x^T\mathcal{M} = \mathcal{B} -\mathcal{M} D_x \Longleftrightarrow   D_x^T= \mathcal{B} \mathcal{M}^{-1}  -\mathcal{M} D_x\mathcal{M}^{-1}    \Longleftrightarrow   D_x^T= \mathcal{M}^{-1}  \mathcal{B}  -  \widetilde{D}_x \end{equation*}
 By applying summation-by-parts on the second term of \eqref{eq:dgsem2}, the 
  semi-discrete   equation \eqref{eq:dgsem2}  can be recast to 
 }  
 \begin{equation}\label{eq:dgsem3}   
     \dfrac{\textrm{d}  \mathbf{U}}{\textrm{d}t} +    \revp{D}_x \mathbf{F}  + \mathcal{M}^{-1}  \mathcal{B}( \widehat{\mathbf{F}} -\mathbf{F})=0.
      \end{equation} 
In other words, we can recast  the DGSEM  semi-discrete equations using a residual or a fluctuation distribution  form. \revp{Therefore, we multiply the mass matrix $\mathcal{M}$ again and consider the $i-th$ degree of freedom, i.e. the $i$-th GL point. The scheme is given by,} 
  \begin{equation}\label{eq:dgsem4}   
w_{i}\dfrac{dU_{i}}{dt} +   \Phi_{i}+  \Psi_{i}^{L} + \Psi_{i}^{R} =0,
\end{equation}
with cell and face residuals arising for the GL quadrature of 
\begin{equation}\label{eq:dgsem5}    
\begin{split}
\Phi_{i}: =&\int\limits_K\phi_i \partial_xF_{\hh}, \\[5pt]
\Psi_{i}^{L}:= [\phi_i(\hat F_{\hh}-F_{\hh})]_{\xi=0}&\;,\;\;
\Psi_{i}^{R}:= [\phi_i(\hat F_{\hh}-F_{\hh})]_{\xi=1}.
\end{split}
\end{equation}
Consistency and accuracy conditions can be characterized in terms of the properties of the split residuals, as discussed in thorough detail
in \cite{AR:17,Abgrall2022}. More importantly, from the point of view of preservation of  steady states, this form is very interesting. Indeed, due to the nodal Gauss-Lobatto
approximation, continuous solutions with  the constant flux $F_{\hh}=F_0$ are exact discrete steady states.

Note that for general numerical fluxes of the form
\begin{equation}\label{eq:dgsem6}
\hat F_{\hh} = \alpha F_{\hh}^+ + (1-\alpha) F_{\hh} +  \mathcal{D} ( U_{\hh}^+ -U_{\hh}),
      \end{equation} 
            where $\alpha$ can be a scalar or a matrix weight, and $\mathcal{D}$ a positive definite matrix,
    the DGSEM semi-discrete equation \eqref{eq:dgsem2} can also be written as 
  \begin{equation}\label{eq:dgsem7}   
     \dfrac{\textrm{d}  \mathbf{U}}{\textrm{d}t} +    \widetilde{D}_x \mathbf{F}  +  \mathcal{M}^{-1}  \mathcal{B}(\alpha[\![\mathbf{F} ]\!])+    \mathcal{M}^{-1}  \mathcal{B}(\mathcal{D}[\![\mathbf{U} ]\!]) =0
      \end{equation}    
      having introduced the interface jumps  $[\![\cdot]\!] =(\cdot)^+-(\cdot)^-$.

% 
%\begin{remark}{(2D extension on tensor product elements)}
%%
%All of the  above formulas can be trivially extended to the multidimensional case using the tensor product representation
%sketched in section for details on the extension of the DGSEM approximation. .  This is relatively classical, and 
%the final prototype has the exact same properties as those discussed above. To save space we focus the discussion in most of the paper on the one-dimensional case,
%and we refer the reader to  
%\end{remark}
%  
%
% 

\section{Discrete full well-balancing and global flux quadrature}\label{sec:full}

 We now consider the numerical approximation of   a  (hyperbolic) system of balance laws
\begin{equation}\label{eq:dgsem-gf0}
	\partial_tU+\partial_x F(U)=  S(U;\varphi(x))  .
\end{equation}
A steady state $U^*(x)$ of \eqref{eq:dgsem-gf0}  can be   described  by the  non-linear ODE 
\begin{equation}\label{eq:dgsem-gf1}
 F'(U^*)(x)=  S(U^*;\varphi(x))  .
\end{equation}
At the continuous  level,   \eqref{eq:dgsem-gf1} can be equivalently  expressed as a non-local  integral equation
\begin{equation}\label{eq:dgsem-gf2}
 F(U^*)(x) = F_0 + \int_{x_0}^{x} S(U^*;\varphi)(s) 
\end{equation}
for a given initial integration point $x_0$, for example the left-end of the domain, and with appropriate definition of the initial state $F_0$.
Integral relations  as \eqref{eq:dgsem-gf2}  are  the most general way of defining reference solutions, as shown
by some of the definitions of section \S2.1.   At the discrete level, however,
given a high order strategy to discretize the local equation  \eqref{eq:dgsem-gf1}, it is in general not true that  the discrete equations
for all point (or average) values   correspond to a   discretization of \eqref{eq:dgsem-gf2}.
We thus give the following definition in the spirit of \cite{Castro2020} and subsequent works.
\begin{definition}[Discrete fully well-balanced] \revt{Consider a scheme  discretized  
using some mesh data  (point values, cell averages, etc.). 
Assume that  using the same data we can construct a discrete approximation
of any  continuous exact steady state with enhanced accuracy compared to the scheme, either in terms of error convergence rate, or of error magnitude, or both.
If this approximation also satisfies exactly  the  stationary  algebraic equations  obtained from the  discretization of 
\eqref{eq:dgsem-gf0}, then we say that the scheme is discretely fully well balanced.
In other words,  a discretely fully well-balanced scheme is superconvergent for any given  continuous exact steady state.}
%
%
%at steady state there is a higher order app
% data evolved by the scheme 
%
%the nodal values provided define an approximation 
% solution can be
%fully well-balanced if (a) the discrete stationary solutions exist (i.e. the ODE system corresponding
%to the semi-discrete method has equilibria) and (b) every
%approximated by a discrete one with enhanced accuracy. Moreover, the method is said to be exactly
%well-balanced for a particular stationary solution if its point values or cell averages (depending on
%the method) constitute a discrete stationary solution.
%the steady discrete solution for all the local point (or averaged) values  can be equivalently obtained   
%as from the steady limit of the  scheme  itself, or from a   direct discretization of   the non-local  equation \eqref{eq:dgsem-gf2}.
%A discrete well-balanced scheme is  exactly well-balanced with respect to the discrete solution  defined by the  discretized integral equation.
\end{definition}
The definition above assumes that given a set of point/averaged values, we are able to construct with the same data both a local
approximation of the balance law, reducing to a discretization of \eqref{eq:dgsem-gf1} at steady state, 
as well as a consistent discretization of \eqref{eq:dgsem-gf2}, and that   the two  are equivalent at steady state.

Fully well-balanced schemes verifying this definition have been proposed e.g.  in \cite{math9151799,math10010015}.
In these works, the authors exploit the   corrected approach recalled in the introduction,  initially proposed in \cite{Castro2020},
combined with a collocation method to formulate the discrete integral equations. This approach thus requires the explicit computation
of $U^*$. In this work, we proceed differently, and show that for   a specific approximation of the integral of the source arising in the finite element 
statement  the above  equivalence holds. This provides an approach which does not necessarily require the evaluation of $U^*$.\\

We start from the   global   formulation in which we replace the source by the derivative of an unknown flux $R$ in the DG variational form
\begin{equation}\label{eq:dgsem-gf3}
\int\limits_{K}\phi_i S = \int\limits_{K}\phi_i \partial_x R.
\end{equation}
We then evaluate $R$ from some approximation of the integral relation
\begin{equation}\label{eq:dgsem-gf4}
\begin{split}
 R(U;\varphi(x)) =& \;r_0 +\int_{x_0}^{x} S(U;\varphi)(s) ,
 \end{split}
\end{equation}
where $x_0$ denotes the left end of the spatial domain, and $r_0$ the left value of the source flux.
Note that this value is essentially an integration constant not affecting the solution. It  is  set to  $r_0=0$ in the following. 
Introducing the  global flux $G=F-R$, we then  simply apply the standard DGSEM approach to the global flux form of    \eqref{eq:dgsem-gf0}  
\begin{equation}\label{eq:dgsem-gf5}
	\partial_tU+ \partial_x G(U;\varphi(x))= 0.
\end{equation}
This is obviously a good setting to obtain    consistency with \eqref{eq:dgsem-gf2}, since at steady state $F=R$. The details  of this consistency  depend  on   how  $R(U;\varphi(x))$ is evaluated   at mesh nodes,
or otherwise said on how \eqref{eq:dgsem-gf4} is approximated.
In the  DGSEM setting, the most natural  approach is to introduce   discrete elemental approximations in the local polynomial space for all the quantities involved, and in particular for
$S_{\hh}$,  $R_{\hh}$, and $\varphi_{\hh}$, as described in section \S\ref{sec:dgsem}.  Then, in each  element $K$ we compute the nodal values $\{R_i\}_{i=0,p}$ as 

 \begin{enumerate}
 \item Set $R_0   = r^-$ 
 \item  For $i=1,p$ set  $R_i = R_{i-1}-\int\limits_{x_{i-1}}^{x_{i}}S_{\hh}(x)$
 \end{enumerate}
with  $r^-$ a local  integration constant  to be defined.

Using the \revp{spectral element} expansion for   $S_{\hh}$ we can recast
the previous expressions   in a  compact form by introducing
   the following  $(p+1)\times (p+1)$  integration matrix $\mathcal{I}$ 
\begin{equation}\label{eq:dgsem-gf6}
 \mathcal{I}_{jk}  := \int\limits_0^{\xi_j}\phi_k(\xi)\,d\xi. 
 \end{equation}
With the previous definition we now have
\begin{equation}\label{eq:dgsem-gf7}
 \mathbf{R} =  \mathbf{R}^-   - \hh\, \mathcal{I}  \mathbf{S},
 \end{equation}
 where now $ \mathbf{R} $ is the array of nodal values of $R$, $\mathbf{R}^-$ has  entries  all equal to  $r^-$,
 and $ \mathbf{S}$ contains the nodal values of the source. 
 The semi-discrete global flux DGSEM equations can now be readily written in matrix form (cf.  \eqref{eq:dgsem3}) as
\begin{equation}\label{eq:dgsem-gf8}
     \dfrac{\textrm{d}  \mathbf{U}}{\textrm{d}t} +    \widetilde{D}_x \mathbf{G}  + \mathcal{M}^{-1}  \mathcal{B}( \widehat{\mathbf{G}} -\mathbf{G}) =0
      \end{equation} 
with $\mathbf{G}  =\mathbf{F} - \mathbf{R} $, or equivalently in the fluctuation form  \eqref{eq:dgsem4} by setting
\begin{equation}\label{eq:dgsem-gf9}
\begin{split}
\Phi_{i}: =&\int\limits_K\phi_i \partial_xG_{\hh}\;,\quad\left\{
\begin{array}{lll}
\Psi_{i}^{L}:= [\phi_i(\hat G_{\hh}-G_{\hh})]_{\xi=0},\\[5pt]
\Psi_{i}^{R}:= [\phi_i(\hat G_{\hh}-G_{\hh})]_{\xi=1}.
\end{array}
\right.
\end{split}
      \end{equation}

To proceed further we  use an explicit form of the numerical  fluxes  similar to \eqref{eq:dgsem6}, namely
\begin{equation}\label{eq:dgsem-gf10}
\hat G_{\hh} = \alpha G_{\hh}^+ + (1-\alpha) G_{\hh} + \mathcal{D}( U_{\hh}^+ -U_{\hh})
      \end{equation} 
Making explicit use of \eqref{eq:dgsem-gf7} and of \eqref{eq:dgsem6} we can write global flux DGSEM as 
\begin{equation}\label{eq:dgsem-gf11}
 \dfrac{\textrm{d}  \mathbf{U}}{\textrm{d}t} +    \widetilde{D}_x \mathbf{F}  + \mathcal{M}^{-1}  \mathcal{B}( \alpha [\![\mathbf{F}]\!]) 
  + \mathcal{M}^{-1}  \mathcal{B}( \mathcal{D} [\![\mathbf{U}]\!]) 
 =
    \hh \, \widetilde{D}_x \mathcal{I}  \mathbf{S} +  \mathcal{M}^{-1}  \mathcal{B}(\alpha[\![\mathbf{R}]\!]) .
      \end{equation} 
The last term only involves jumps in the source flux, depending on the definition of the local integration constant $r^-$. 
  Consistently with \revp{Remark \ref{remark_nodal_data}},  in this work we have not accounted for this jump. 
\revPO{We refer to \cite{ctr23} for a way to include this contribution. Neglecting the jump}
   is equivalent to  the choice $r^- =(R_p)^{K^-}$, the last value of the left neighbouring element.
This leads to the following  \emph{DGSEM formulation with global flux quadrature of the  source}:
  \begin{equation}\label{eq:dgsem-gf12}
 \dfrac{\textrm{d}  \mathbf{U}}{\textrm{d}t} +    \widetilde{D}_x \mathbf{F}  + \mathcal{M}^{-1}  \mathcal{B}( \alpha [\![\mathbf{F}]\!]) 
  + \mathcal{M}^{-1}  \mathcal{B}( \mathcal{D} [\![\mathbf{U}]\!]) 
 =
    \hh \, \widetilde{D}_x \mathcal{I}  \mathbf{S}  .
      \end{equation} 
 The only remaining ingredient is the definition of the nodal value of the source. \\
      
      %\begin{remark}[Global flux quadrature]
  The above formula shows that  with our choices the global flux  DGSEM approach boils down to a very specific quadrature of the source with weights
       provided for each nodal degree of freedom by the  integration tableau $\mathcal{I}$. In particular, the new scheme
    only requires modifying the mass matrix in front of the source term. 
       \revPO{For this reason we find it more appropriate, instead of a global flux method, 
       to speak of  a DGSEM method with  \emph{global flux quadrature}. For a given function $f$, global flux quadrature corresponds to  the  approximation
  \begin{equation}\label{eq:dgsem-gf13}
  \int\limits_K\phi_i f_{\hh}\,dx = \hh \,(\widetilde{D}_x \mathcal{I}  \mathbf{f}  )_i
        \end{equation} 
%  rather than global flux method.  
%  Note that this view is  note only philosophical. 
  The formula above shows that the global flux quadrature is actually a fully local method. 
  Note that this  locality is also true for   \eqref{eq:dgsem-gf11}, as
long as one can express the jumps in the source flux $[\![\mathbf{R}]\!]$ as a function of the source term itself.
Preliminary work in this sense is discussed in \cite{ctr23}.} 
%\end{remark}

\subsection{Discrete equilibria and connection with collocation methods} For smooth solutions, the DGSEM formulation with global flux quadrature of the  source
has a neat connection with continuous collocation methods for ODEs.  Summarised by the following proposition. 

%
%We are now ready to provide a more precise
% characterization of the   the discrete steady states of remark \ref{rem:steadysol}.  We start by noting that the integration matrix $\mathcal{I}$ is
%nothing else that the Butcher tableau of the collocation  ODE integration method associated to the interpolation nodes.
%This tableau naturally emerges from integrating from one collocation point to the next the ODE
%\begin{equation}\label{eq:dgsem-gf14}
%U'(x) + S(U,x) =0     
%\end{equation}
%Given a discretization of the independent variable domain in slabs  $[x_k,x_k+\hh]$, 
%integration between the collocation points  within each slab leads to
%an implicit full tableau method which can be recast as   \cite{Prothero74}
%\begin{equation}\label{eq:dgsem-gf15}
%\mathbf{U}_{k+1} = \mathbf{U}^- - \hh\, \mathcal{I}\mathbf{S}  
%\end{equation} 
%with $\mathbf{U}_{k+1}$ containing all the values of $U$ at intermediate stages in the slab $[x_k,x_k+\hh]$, and  $\mathbf{U}^-$ the last available value from the previous slab.
%In particular, in the case of the Gauss-Lobatto nodes used here, the matrix $ \mathcal{I}$ is the integration tableau of the 
%$p+1$ stages RK-LobattoIIIA  ODE solver \cite{Prothero74}.   This allows to give a precise meaning to  the solutions of the scheme, 
%which is summarized in the following proposition.
\begin{proposition}[Global flux Gauss-Lobatto DGSEM: discrete steady states]\label{proposition_4} 
\revt{Provided that the mapping $F(U)$ is invertible and that the inverse $U(F)$ is bounded and  uniquely defined,
then the global flux Gauss-Lobatto DGSEM equations \eqref{eq:dgsem-gf12} admit a nodally continuous discrete steady state 
$U^*=U^*(F)$ obtained upon integration with the fully implicit    continuous collocation  RK-LobattoIIIA  method
 of the nonlinear ODE \eqref{eq:dgsem-gf1} with  $S(U^*;\varphi) =S(U^*(F);\varphi)$, and with initial condition $U_0$ satisfying $F(U_0) = F_0$ with $F_0$ given. }
      \begin{proof}
     \revt{Note that the entries of the matrix $ \mathcal{I}$ associated to the Gauss-Lobatto points/basis 
        is by definition   the integration tableau of the  $p+1$ stages fully implicit RK-LobattoIIIA  ODE solver \cite{Prothero74}.
         So the steady state  defined in the proposition can be written element by element as 
         $$
         \mathbf{F}=\hh \mathcal{I} \mathbf{S},
         $$
         which by using the definition of the source flux and the initial condition  implies  $F_i-R_i=F_0$ $\forall \, i$.
         The condition $[\![\mathbf{R}]\!]=0$   boils down (only at steady state) to  $[\![\mathbf{F}]\!]=0$, from which the continuity of $U^*$ follows from 
          the  properties of the map $U(F)$. Putting all this together shows that $U^*$ so defined is a solution of the steady discrete equations 
           \eqref{eq:dgsem-gf12}. }
      \end{proof}
\end{proposition}

The proposition states simply that  at steady state the DGSEM method with global flux quadrature provides  a direct approximation of \eqref{eq:dgsem-gf2} in which $S$ is replaced by the piecewise Lagrange polynomial
through the Gauss-Lobatto interpolation points. \revt{There are three very important aspects which need to be addressed: 
\begin{itemize}
\item  the implication of the above proposition on the accuracy of the discrete steady state $U^*$; 
\item the validity of the  hypotheses (e.g. on the map $U(F)$) required for the  analysis to be true;
 \item and the implication of this analysis on the implementation of the method, and   some clarifications
 with respect to the relations with schemes having similar properties proposed in \cite{Castro2020,math9151799,math10010015}.
\end{itemize}
}

 %the implication of this analysis on the implementation of the method, and   some clarifications
 %with respect to the relations with schemes having similar properties proposed 
 %e.g.   proposed by \cite{Castro2020,math9151799,math10010015}.}

\revt{Concerning the first aspect,} as a corollary of the last proposition we have the following property.

\begin{corollary}[Global flux Gauss-Lobatto DGSEM: superconvergence at steady state]\label{superconvergence} 
\revt{Consider a system of conservation laws for which one can exhibit a flux linearization  reading
\begin{equation}\label{eq:Flin}
F(U)-F(V) = A(U,V)(U-V)
\end{equation}
with $A$  be a matrix with entries uniformly bounded with respect to its arguments, and diagonalizable with eigenvalues $\{\lambda_j^A\}_{j\ge 1}$.
Provided that there exist a bounded strictly positive constant $C_A$ such that
\begin{equation}\label{eq:lambda}
\lambda_{\min}^A := \min_{j \ge 1}|\lambda_j| \ge C_{A} > 0,
\end{equation}
then under the hypotheses of  Proposition \ref{proposition_4}, the discrete steady state   $U^*$ 
is nodally superconvergent, and in particular the element endpoint values of $U^*$ are approximations of order $2p$ of any smooth  exact steady  solution, 
while the internal nodal values have accuracy  $\hh^{\min(p+2,2p)}$.}
       \begin{proof}
       \revt{There  are two parts to the proof. The first is to  invoke the properties  of the LobattoIIIA method to 
        argue that $F(U^*)$ has an accuracy of order  $\hh^{2p}$ at the elements endpoints and $\hh^{\min(p+2,2p)}$ at the
        internal nodes. One can refer to e.g. Theorem 7.10 in \cite{hairer} (see also \cite{Prothero74}) for the proof.
        The second part is to bound the error on the solution. To this end, we use the following 
        % \revPO{I have my doubt about this estimattion, I am missing the $\lambda^2$ and also we should lead $1/C_A$ then above or not???}
        $$
        (F(U_i^*) - F(U_i^{\text{exact}}))^2 = \left( A(U_i^*,U_i^{\text{exact}}) (U_i^* - U_i^{\text{exact}})\right)^2  \ge (\lambda_{\min}^A)^2 (U_i^* - U_i^{\text{exact}})^2
        $$
        which using \eqref{eq:lambda} readily leads to
        $$
        \|U_i^* - U_i^{\text{exact}} \| \le \dfrac{1}{C_{A}}  \|F(U_i^*) - F(U_i^{\text{exact}})\| 
        $$
        and thus the result.}
      \end{proof}
\end{corollary}

\revt{The corollary above shows the main added value of the global flux quadrature:  the consistency with 
direct integration of the ODE by means of a Gauss-Collocation method, by which we inherit its accuracy, and in particular, superconvergent properties.
This can be interpreted now in two ways. The first is to say, following the approach of \cite{Castro2020,math9151799,math10010015},
that the global flux Gauss-Lobatto DGSEM scheme is \emph{exactly  and fully well-balanced} with respect to  all  solutions $U^*$ characterized by Proposition \ref{proposition_4}
and Corollary \ref{superconvergence}. We can otherwise  understand the corollary as the scheme having enhanced
accuracy, with orders $\hh^{2p}$ for element boundaries and $\hh^{\min(2p,p+2)}$ for internal nodes, 
\emph{for all smooth steady states}.  }\\

\revt{Concerning now the validity of the hypotheses used, some of them have some simple and clear justifications. For example,  linearizations of the type \eqref{eq:Flin}
are known and are exhibited by many systems of balance laws. For the shallow water equations, 
the matrix $A$ can be obtained by means of a conservative  Roe-like linearization of the Jacobian of the flux $\partial_U F$,
see e.g.  \cite{bv94} (and \cite{SARMANY201386} \S3.1.1 for a multidimensional generalization). Concerning condition \eqref{eq:lambda},
for the classical shallow water equations it  implies that the analysis only holds far from critical points.
We will verify   that  in practice the method still performs well across such points. }

\revt{While there may be a strategy to avoid \eqref{eq:lambda} and prove the above result in a more general setting, 
we are going to see shortly that critical points also play a role for 
 the strongest hypothesis  made, which is the one on $F(U)$ being invertible,  and $U(F)$ being bounded and unique. 
 For several hyperbolic systems in 1D  precise  conditions can be provided to characterize this issue.
For  the shallow water equations, the study of such conditions is classical.
One can easily show  that by setting    $F=(q, M)$  in   the non-trivial case $q^2 >0$  given admissible data
verifying $ M >  3\sqrt[3]{gq^4}/2$, one can always compute two non-negative values of the depth
corresponding to these data,  a   unique sub-critical one and a unique super-critical one. So the inversion
can be performed as soon as one knows the nature of the flow. The case $M=3\sqrt[3]{gq^4}/2$
defines a critical solution (we refer e.g. to  \cite{CCHKT19} for details).
This shows that the inversion of the mapping requires a-priori knowledge of the  super- or sub-critical  of the flow. So if one was
able to apply the ODE solver   using the flux as main unknown,  flows with transition from sub- to super-critical may be problematic.}\\

\revt{ This brings us to the last aspect of the discussion.
%
% 
% given for example in
%We provide here a   different version of the study. Setting   $F=(q, M)$   it can be easily shown (see e.g. \cite{CCHKT19}) that the computing $U(F)$ is equivalent to 
%solving the cubic equation
%\begin{equation}\label{eq:cube}
%h^3 - \dfrac{2M}{g}h +\dfrac{2q^2}{g} =0\,.
%\end{equation}
%%For
% }
%
% 
%  \revt{\begin{proposition}[1D shallow water system: invertibility of $F(U)$] \label{prop:finv} Let $\sigma=\text{sign}(Fr-1)$ define the sub- or super-critical
%  nature of the solution, with $\sigma=0$ by definition defining the critical point, and where
%  $$
%  Fr := \dfrac{|u|}{\sqrt{gh}}.
%  $$
%   Given a  tuple $T=(q,M,\sigma)$,   with $M >0$ and  $\sigma\in\{-1,0,+1\}$, there exist a unique solution to  problem $\eqref{eq:cube}$ matching $T$ under the  
%    admissibility condition
%    \begin{equation}\label{eq:cube-adm}
%%    \sigma =0 \rightarrow M =\dfrac{3}{2}\sqrt[3]{gq^4}\,,\;\;
%    \sigma \ne 0 \rightarrow M >  \dfrac{3}{2}\sqrt[3]{gq^4}
%    \end{equation}
%  \end{proposition}
%  \begin{proof}
%  The proof is an extended and modified version of the arguments given  in \cite{nxs07} for a similar problem. It is reported in appendix for completeness.
%  \end{proof}
%}
%
%\revt{The proposition above shows that 
%The above analysis shows is that, under some conditions we just discussed, at steady state t
The use of the global flux quadrature allows to make a direct link between the 
DGSEM method and continuous collocation  RK methods 
applied to the  steady ODE.
% or equivalently Gauss-Lobatto integration directly applied  to \eqref{eq:dgsem-gf2} by sub-intervals of size $\hh$.
As already mentioned, other methods exploit collocation ODE integrators e.g. in \cite{math9151799,math10010015}. 
These methods use locally the ODE integration to construct a reference discrete  steady solution and correct the scheme so that it is exact in correspondence of such solution.
In both cases,  exact well-balancing is  only obtained in correspondence of certain steady  states as e.g. the lake at rest (cf. below).
Otherwise, the schemes preserve only approximate solutions corresponding to those obtained by the ODE integrators.
The strong point of both, the method proposed here and those in the references, is that nothing of these steady  solutions   needs to be known a-priori.
For  our proposed method however note that   the solution of the ODE integrator is \textbf{not} needed. 
The analogy is in fact   just a tool to characterize the steady states. As we will see in the numerical experiments, for example  there is absolutely no issue in applying the method
to trans-critical problems with excellent results, although the current proofs do not apply to such cases.
The fact that the solution of the ODE is not required to implement the method is a net 
 advantage of the scheme studied in this paper, compared to those proposed in  \cite{math9151799,math10010015}.
Conversely,  its drawback is that there is no flexibility in the choice of the collocation method.
In other words, once the polynomial degree is fixed, so is the best possible approximation of the steady state.
For the scheme in  \cite{math9151799,math10010015} the two ingredients are somehow independent, which is an advantage.
}

\subsection{\revPO{Modified evaluation of the nodal bathymetric source to get exact lake at rest states}}

The DGSEM  approach with global flux quadrature  described in the previous section requires the definition of the
nodal values of the source. For the quasi-1D shallow water system \eqref{SW1}, if $\mathbf{u}=(u,v)^T$ denotes the velocity vector, and $\mathbf{u}^{\perp}=(v,-u)^T$ its orthogonal,    both the friction term $c_{f} u$, and 
the Coriolis terms $\omega\mathbf{u}^{\perp}$ can be easily evaluated nodally. The  tricky part is how to define the value of the term $ h\partial_x \varphi=gh\partial_x b$
related to the potential effects due to  bathymetric variations.

In this work, we have used two approaches. The first is a straightforward   analytical evaluation of  $\partial_x b(x)$: 
  \begin{equation}\label{eq:S_node}
  S_l =-h_l\left( \begin{array}{c}
  0\\
  g\partial_xb(x_l) + c_f u_l +\omega v_l\\
  -\omega u_l
  \end{array}\right).
  \end{equation}
  To be  exactly well-balanced for the  hydrostatic equilibrium, we adapt an idea proposed in \cite{shu2006highorder}.
Denoting by $\zeta$ the free surface level  $\zeta =h+b$, we set
  \begin{equation}\label{eq:S_node-mod}
  S_l =-\left( \begin{array}{c}
  0\\
  g\zeta_l \partial_xb_{\hh}(x_l) - \partial_xp_{\hh}(b)(x_l) + c_f h_lu_l +\omega h_lv_l\\
  -\omega h_lu_l
  \end{array}\right),
  \end{equation}
  where now $b_{\hh}$ is the finite element expansion built using the nodal values of the bathymetry, and similarly for $p_{\hh}(b)=g(b^2/2)_{\hh}$.
With this definition we can prove the following simple result.

\begin{proposition}[Exact well-balanced for lake at rest] If the interpolation polynomial $b_{\hh}$   is continuous, then 
the DGSEM formulation with global flux quadrature \eqref{eq:dgsem-gf12},
and with nodal source given by \eqref{eq:S_node-mod} is exactly well-balanced for lake at rest states $\zeta=\zeta^*=\text{const}$,
and $\mathbf{u}=0$.
\begin{proof}
\revt{%To obtain the result we can analyze term by term scheme \eqref{eq:dgsem-gf12} to show that under the hypotheses made $\dfrac{d\mathbf{U}}{dt}=0$.
First of all note that since   $b_{\hh}$ is continuous (cf. Remark  1.)   and so is $\zeta^*$ which is constant,  $h$ will also be continuous. 
Since also $h\mathbf{u} =0$ is constant, then we have $[\![\mathbf{U}]\!] =0$ and   $[\![\mathbf{F}]\!] =0$ at all element boundaries.
The remaining term can be written as 
$$
\dfrac{d\mathbf{U}}{dt}=-\widetilde D_x( \mathbf{F} -   \mathbf{F}_0  -\hh\mathcal{I}\mathbf{S}   )
$$
with $ \mathbf{F}_0$ be  the array containing the values  $F(U_0)$ with $U_0$ containing the value of the unknowns in the first Gauss-Lobatto 
point within the element. Because of the hypotheses made, only the hydrostatic terms remain in $\mathbf{F}$ and $\mathbf{S}$.
For the source,  using the identity  $ \partial_x f_{\hh} = \sum_{l=0,p}\partial_x \phi_l(x)f_ l = \sum_{l=0,p}\phi_l(x) \partial_x f_{\hh}(x_l)$,
 we can easily show that 
\begin{equation*}
\begin{split}
(\hh\mathcal{I}\mathbf{S})_i =& \int_{x_0}^{x_i}g\zeta^*\sum\limits_{l=0,p}\phi_l(x_l)\partial_xb_{\hh}(x_l)  - \int_{x_0}^{x_i}g\sum\limits_{l=0,p}\phi_l(x_l)\partial_x  p_{\hh}(b)(x_l)\\
=& g\zeta^* \int_{x_0}^{x_i}\partial_x b_{\hh} - \int_{x_0}^{x_i}\partial_x p_{\hh}(b)=
g\zeta^*(b_i - b_0) - g\left(\dfrac{b_i^2}{2}-\dfrac{b_0^2}{2}\right).
%\int_{x_0}^{x_i}g\zeta^*\sum\limits_{l=0,p}\partial_x \phi_l(x_l) b_l  - \int_{x_0}^{x_i}g\sum\limits_{l=0,p}\partial_x\phi_l(x_l) p_{\hh}(b_l) 
\end{split}
\end{equation*}
%having used the approximation  identity $ \partial_x f_{\hh} = \sum_{l=0,p}\partial_x \phi_l(x)f_ l = \sum_{l=0,p}\phi_l(x) \partial_x f_{\hh}(x_l)$.
%
%is related to the fact that the  derivative of the finite element  polynomial is in a lower dimensional space than
%the  polynomial itself, and can be written trivially in the larger space as the   interpolation  of its nodal values. 
%As a consequence we get
%\begin{equation*}
%\begin{split}
%(\hh\mathcal{I}\mathbf{S})_i =& .
%\end{split}
%\end{equation*}
Using  the fact that $b_i+h_i=\zeta^*=b_0+h_0$ we deduce  $(\hh\mathcal{I}\mathbf{S})_i =p_{\hh}(h_i) -p_{\hh}(h_0)$, and thus the result.}
\end{proof}
\end{proposition}
This second approach will be  referred to in the results section as \emph{modified global flux  quadrature}, being obtained with a modification of the nodal source terms
allowing to guarantee the exact preservation of lake at rest states. 

\section{Entropy control via cell corrections:  flux vs energy conservation}\label{sec:EC_global}

We consider now the issue of the compatibility between the consistency with constant global flux,
underpinning notion of the previous section, and the consistency with a given pair entropy/entropy flux.
In absence of friction, smooth solutions of the  quasi one dimensional shallow water equations \eqref{SW1}  embed an additional conservation
for the total entropy $\eta_{\varphi}$ with flux $F_{\eta_{\varphi}}$  defined in \eqref{SW-eta2}.
Defining the total energy density $E=g\zeta + k$, we can readily show that analytical steady states of   \eqref{SW1}  (without friction)
are characterized by the invariants  
\begin{equation}\label{inv-eta0}
\begin{split}
hu = &\, q_0,\\
E =& \,E_0 ,\\
v -\omega x = & \,v_0.
\end{split}
\end{equation} 
These three invariants are compatible with both the steady equations and a constant distribution of the entropy flux $F_{\eta_{\varphi}}$.
The same equations, in global flux  form  provide a steady state described by the invariance of the global flux components:
\begin{equation}\label{inv-eta1}
\begin{split}
hu = &\, q_0,\\
hu^2 +p + r_u  =& \,Q_0, \\
huv  + r_{v}  = & \,V_0,
\end{split}
\end{equation} 
having denoted by $r_u$ and $r_{v}$ the components of the source flux arising from the second and third equations.
These three relations are by construction compatible with the discrete full well-balanced property, and the related superconvergence property of \revp{Corollary \ref{superconvergence}}.     At the continuous level all is fine. 

At the discrete level, however,   the global flux quadrature is exactly consistent only with  \eqref{inv-eta1},
and only within its error with the first. 
In other words,  when looking at the preservation of discrete steady states, \revp{we may have two different situations depending on the selected data initialization strategy:}
\begin{enumerate}
\item We use exact analytical expressions, satisfying  \eqref{inv-eta0}, to initialize the data.  In this case,  we may use the analytical expressions for the entropy fluxes  to construct entropy conservative/dissipative schemes. 
Scheme \eqref{eq:dgsem-gf12} will still  be discretely fully well-balanced within
the conditions of  \revp{Corollary \ref{superconvergence}.} 
\item We choose to  initialize with the global flux solution $U^*$ of \revp{ Proposition \ref{proposition_4}}. In this case, the resulting accuracy will still
be the one   of  \revp{Corollary \ref{superconvergence}}, but scheme  \eqref{eq:dgsem-gf12} would be exactly well-balanced with respect to this initial steady state.
For this initialization, however using the standard analytical expressions for the entropy fluxes,  would break this property. 
\end{enumerate}
The objective of the following subsections is to propose a modification of the schemes which allows to solve this dichotomy, and provide a degree of control 
on the entropy production of the method independently on the initialization procedure.

\subsection{Cell entropy correction method}
There exist several well established techniques to embed \revp{DGSEM} methods with degree of control of the entropy evolution.
We refer to the introduction for a short overview.   The approach used here is based on the idea of \cite{abgrall2018general,abgrall2022reinterpretation} to
introduce   local  corrections  in the form of a \revp{symmetric positive definite (SPD)} bilinear term with a free coefficient allowing to impose a desired constraint in terms of entropy balance.

To briefly describe the method, we start from the frictionless case. We consider  the fluctuation form \eqref{eq:dgsem4} (with \eqref{eq:dgsem-gf9}) to which we add the correction term:
\begin{equation}\label{eta-corr0}
w_i\dfrac{dU_i}{dt} + \Phi_i +\Psi_i^{\text{L}}+\Psi_i^{\text{R}} + \alpha_K \mathcal{D}_i^K=0\;,\quad
 \mathcal{D}_i^K:=\int\limits_K\partial_x\phi_i \, A_0 \partial_x W_{\hh} ,
\end{equation} 
where $W$ is the appropriate set of entropy variables such that 
\begin{equation}\label{eta-corr1}
\partial  \eta  = W^t\partial U,
\end{equation} 
and $A_0$ is the inverse of the Hessian    $A_0^{-1}=\partial_{UU}\eta = \partial W/\partial U$.  To obtain the entropy balance we need only to dot \eqref{eta-corr0}
 by $W_i^t$ and sum over $i=\{0,p\}$. \revp{The aim is to specify the correction term $ \alpha_K \mathcal{D}_i^K$, in particular  $\alpha_K $ that conservation is not violated, but entropy conservation is additionally ensured. }

 This readily leads to the cell entropy evolution equation
\begin{equation}\label{eta-corr2}
|K| \dfrac{d\bar \eta_K}{dt} +  \Phi_{\eta}^K +   \alpha_K  \|\partial_x W \|^2_{L^2_{A_0}\!(K)} =0 
\end{equation} 
having introduced the average cell entropy $ \bar\eta_K $, and cell entropy production of the scheme $ \Phi_{\eta}^K$:
\begin{equation}\label{eta-corr3}
\begin{split}
 \bar\eta_K := & \;\dfrac{1}{|K|}\int\limits_K\eta_{\hh}, \\[5pt]
 \Phi_{\eta}^K := & \sum\limits_{i=0,p}W_i^t\left(\Phi_i +\Psi_i^{\text{L}}+\Psi_i^{\text{R}},
 \right)
\end{split}
\end{equation} 
and norm \revp{$ \|g \|_{L^2_{A_0}\!(K)}:= \left( \int\limits_K g \, A_0 g\right)^{1/2}$.} 
From  \eqref{eta-corr2}, we can see that, unless the solution is locally constant, the entropy  production associated to the correction term multiplied by $\alpha_K$
is strictly positive. This allows to use the free parameter to gain direct control on entropy production. Let  $\hat F_{\eta}$ denote a consistent numerical entropy flux. \revp{Its calculation and specification will be part of the next subsection \ref{subsec_flux}. For now, we set }
  \begin{equation}\label{eta-corr4}
  \Psi_{\eta}^K:= \oint_{\partial K} \hat F_{\eta} ,
  \end{equation} 
  reducing in 1D to 
    \begin{equation}\label{eta-corr5}
  \Psi_{\eta}^K=  (\hat F_{\eta} )_K^R -(\hat F_{\eta} )_K^L .
  \end{equation}
  Also, following \cite{AR:17,Abgrall2022},   for a given smooth exact solution $U^e(x,t)$, and given smooth compactly supported test function $v(x,t)$ define the consistency error 
        \begin{equation}\label{eta-corr6}
        \mathcal{E}:= \sum\limits_{K}\sum\limits_{i\in K}v_i \left\{w_i\dfrac{dU^e_i}{dt} +\Phi_i(U^e_{\hh}) + \Psi_i^R(U^e_{\hh})+ \Psi_i^L(U^e_{\hh}) + \alpha_K(U^e_{\hh})\mathcal{D}_i^K(U^e_{\hh})\right\}.
        \end{equation}
\revPO{with  $U^e_{\hh}$ the projection of $U^e(x,t)$ on the local finite element space}.
 The  proposition below characterizes the
  construction used here, following  \cite{abgrall2018general,abgrall2022reinterpretation}.
  \begin{proposition}[Entropy correction: local/global entropy balance, consistency] \label{prop_entropy}
 The corrected DGSEM scheme obtained from \eqref{eta-corr0} setting 
      \begin{equation}\label{eta-corr7}
  \alpha_K =     \dfrac{\Psi_{\eta}^K - \Phi_{\eta}^K }{ \|\partial_x W_{\hh} \|^2_{L^2_{A_0}\!(K)}}  
  \end{equation}  
verifies the  elemental entropy balance 
   \begin{equation}\label{eta-corr8}
 |K| \dfrac{d\bar \eta_K}{dt}  + \Psi_{\eta}^K = 0
  \end{equation}  
  and, for homogeneous or periodic boundary conditions, the global entropy conservation  equation
     \begin{equation}\label{eta-corr9}
\sum\limits_K |K| \dfrac{d\bar \eta_K}{dt}   = 0.
  \end{equation}  
Moreover, the scheme   verifies a consistency estimate of the type $|\mathcal{E}|=\mathcal{O}(h^{p+1})$.
  \begin{proof}
  Properties \eqref{eta-corr7} and \eqref{eta-corr8} are a straightforward consequence of  \eqref{eta-corr6}. Concerning the consistency estimate,
   the proof  follows  e.g.  Section 3.2 in \cite{AR:17}. \revp{It is reported in appendix with  the necessary definitions from the reference}. 
  \end{proof}
  \end{proposition}

The above proposition shows that the correction approach allows to readily recover  a local entropy balance law with a simple choice of the scalar
coefficient $\alpha_K$.  Moreover, the correction does not spoil the formal consistency of the scheme which remains of order $h^{p+1}$. 
However, as discussed in the beginning of Section \S5,  
the well-balanced properties of the corrected scheme depend on the properties of the numerical entropy flux $ \hat F_{\eta} $. 

\revt{\begin{remark}[Division by zero]
For uniform flows, \eqref{eta-corr7}  gives a zero divided by zero singularity.
In practice, the denominator  is modified to avoid this singularity  as  $\max(\|\partial_x W_{\hh} \|^2_{L^2_{A_0}\!(K)},\epsilon_K)$, with
$\epsilon_K$  constant, and set to $10^{-8}$ in   the numerical experiments. The interested reader can refer to 
\cite{gaburro2023high,abgrall2022reinterpretation,abgrall2022relaxation} for similar modifications.
\end{remark}}

\begin{remark}[Friction and dissipation]\label{remark_3}
In presence of friction (or other dissipative relaxation terms), the above construction needs to be modified by including the effects of these terms 
in the definition of the entropy balance $\Psi_{\eta}^K$ which becomes
     \begin{equation}\label{eta-corr10}
     \Psi_{\eta}^K  := \oint_{\partial K} \hat F_{\eta}  + \mathcal{D}_f 
       \end{equation}  
       with $\mathcal{D}_f $ a consistent approximation of the dissipation
            \begin{equation}\label{eta-corr11}
       \mathcal{D}_f \approx 2 \int\limits_K hc_f k\ge 0
       \end{equation}  
       with  $k$ the kinetic energy. In this case, \revp{Proposition \ref{prop_entropy}} still holds, but \eqref{eta-corr9} becomes
                   \begin{equation}\label{eta-corr12}
      \sum\limits_K |K|\dfrac{d\bar\eta_K}{dt} = -       \sum\limits_K \mathcal{D}_f \le 0
       \end{equation}  
\end{remark}

\subsection{Entropy fluxes compatible with  global flux quadrature}\label{subsec_flux}

We propose here two possible definitions of the numerical entropy flux required in \eqref{eta-corr4} and \eqref{eta-corr5}.
The first is given by 
     \begin{equation}\label{eta-corr13}
 \hat F_{\eta}(U^+_{\hh}, U^-_{\hh}) = \lambda  F_{\eta}(U^+_{\hh}) +  ( 1- \lambda)  F_{\eta}(U^-_{\hh}) 
  \end{equation}  
with $\lambda\ge 0$ some scalar coefficient, and $ F_{\eta}(U)$ the \emph{analytical entropy flux}. \revp{For numerical tests in Section \S7, we choose $\lambda$=$\frac{1}{2}$.}
We then define  the local entropy balance as (we consider the one dimensional case for simplicity)
     \begin{equation}\label{eta-corr14}
     \begin{split}
\Psi_{\eta}^K = &  \hat F_{\eta}(U^+_{\hh}, U^-_{\hh})^R -  \hat F_{\eta}(U^+_{\hh}, U^-_{\hh})^L +  2 \sum\limits_{i=0,p}w_i c_{fi}h_ik_i \\=&
\lambda [\![  F_{\eta}]\!]^R  +  \int\limits_K \partial_x F_{\eta} 
   + (1-\lambda) [\![  F_{\eta}]\!]^L +  2 \sum\limits_{i=0,p}w_i c_{fi}h_ik_i,
  \end{split}
  \end{equation} 
  where the second identity is trivially obtained by adding and removing the internal values of $F_{\eta}$  on the element's boundaries.
Following the discussion
from the beginning of Section \S5, this definition is compatible with analytical steady states verifying  \eqref{inv-eta0}. It is however not
exactly compatible with the global flux solutions of the DGSEM scheme with global flux quadrature, which verify instead  \eqref{inv-eta1}.\\

To obtain a definition which is exactly compatible the global flux quadrature approach proposed in this work, we modify \eqref{eta-corr14}
as 
     \begin{equation}\label{eta-corr15}
     \Psi_{\eta}^K = \lambda [\![  F_{\eta}]\!]^R  +  \int\limits_K W^t_{\hh}  \partial_x G_{\hh}
   + (1-\lambda) [\![  F_{\eta}]\!]^L    
   \end{equation} 
  This definition embeds all the non-differential effects in the global flux term in the middle,  and  it is compatible with global flux solutions verifying  \eqref{inv-eta1}.
  In particular,   the resulting corrected scheme  is  exactly consistent with  the discrete  solutions $U^*$ of \revp{Proposition \ref{proposition_4}}.

\section{A 2D extension}

We briefly discuss here a possible extension to the two dimensional shallow water model \eqref{SW0}. 
This extension is based on a tensor product implementation of the one dimensional  DGSEM scheme,
as well as on a  dimensionally split generalization of the notion of global flux. The study of the genuinely multidimensional
case is object of current work and left out of this paper.  We refer to section \S\ref{sec:notation}
for the notation, and to  \cite{Kopriva2010,Gassner2013,GASSNER201639}  for more details on the implementation.\\

We start by rewriting \eqref{SW0}  as 
\begin{equation}\label{gf-SW0}
\partial_t U + \partial_x F_x  + \partial_y F_y = S\,,
\end{equation}
having set
\begin{equation}\label{gf-SW1}
F_x= \left[ \begin{array}{c} hu\\ hu^2 +p(h)\\ huv \end{array} \right]\,,\;\;
F_y= \left[ \begin{array}{c} hv\\  huv\\hv^2 +p(h) \end{array} \right]\,.
\end{equation}

The   dimension by dimension generalization used here is similar to the one proposed  initially
in \cite{GASCON2001261}, and  boils down to adding an
appropriately defined diagonal flux  tensor to the conservative flux.
In other words, we recast \eqref{SW0} as 
System \eqref{gf-SW0} can be written succinctly as 
\begin{equation}\label{gf-SW2}
	\partial_tU  + 	\partial_xG_x + \partial_y G_y =0
\end{equation}
where 
\begin{equation}\label{gf-SW3}
G_x=F_x +  [0\; R_x\;0]^T\,,\;\;  G_y=F_y+ [0\; 0\;R_y]^T\,,% \left[ \begin{array}{c} 0\\  0\\ R_y \end{array} \right]\,,
\end{equation}
having defined 
\begin{equation}\label{gf-SW4}\begin{split}
	R_x := \int_x h (\partial_x \varphi +c_f u+\omega v )d\xi\,,\;\;
	R_y:=\int_y h (\partial_y \varphi +c_f v-\omega u )d\eta.
\end{split}\end{equation}
The above definition allows to  construct direction by direction fully well-balanced  schemes 
which preserve  global fluxes/equilibrium variables $G_x$ in the $x$ direction, and $G_y$  in the $y$ direction.\\

We discretize  \eqref{gf-SW2} with a classical tensor DGSEM strategy using the Gauss-Lobatto points of the right picture
on figure \ref{fig:elements}. All computations done, and using a notation
similar to that of the previous sections, within each element $K$ we obtain $(p+1)^2$ equations  for the degrees of freedom $U_{ij}$
which can be compactly written as 
\begin{equation}\label{gf-SW5}
\begin{split}
 \dfrac{d\mathbf{U}}{dt}  +  &\widetilde{\mathbf{D}}_x \mathbf{F}_x + \mathcal{M}^{-1}_x  \mathcal{B}_x( \alpha [\![\mathbf{F}_x]\!]_x) 
+ \mathcal{M}^{-1}_x  \mathcal{B}_x( \mathcal{D} [\![\mathbf{U}]\!]_x)\\
  +&  \widetilde{\mathbf{D}}_y \mathbf{F}_y  + \mathcal{M}^{-1}_y  \mathcal{B}_y( \alpha [\![\mathbf{F}_y]\!]_y)  
  + \mathcal{M}^{-1}_y  \mathcal{B}_y( \mathcal{D} [\![\mathbf{U}]\!]_y)
  = \hh\widetilde{\mathbf{D}}_x\mathbf{I}_x \mathbf{S}_x +\hh\widetilde{\mathbf{D}}_y\mathbf{I}_y \mathbf{S}_y
\end{split}\end{equation}
where now the array $\mathbf{U}$ contains line ordered values of the solution at the Gauss-Lobatto collocation points,
$\mathbf{F}_x$ and $\mathbf{F}_y$ the  flux values at collocation points, while $\widetilde{\mathbf{D}}_x$ and $\widetilde{\mathbf{D}}_y$
denote  mono-dimensional derivative matrices with a structure depending on the  numbering used for $\mathbf{F}_x$ and $\mathbf{F}_y$.
If  $\mathbf{F}_x$ is line ordered, and $\mathbf{F}_y$ is ordered by columns, both derivative matrices  have a block diagonal structure with block entries given by the corresponding   one dimensional matrices. 
As before with $ \mathcal{M}$  we denote the one dimensional mass matrices in each direction, and $\mathcal{B}_{x/y}$  are the matrices selecting  values
on the left/right and top/bottom element boundaries, while  $[\![\cdot]\!]_{x/y}$  denote the horizontal/vertical jumps at the left/right and top/bottom element boundaries.
Finally, the right hand side contains the source terms, integrated using one dimensional global flux quadrature involving the derivative matrices, and block matrices
$\mathbf{I}_{x/y}$ whose  entries are the same as the one dimensional ODE collocation methods, and whose structure depends on the ordering of the entries of the 
$\mathbf{S}_{x/y}$ arrays which are values at collocation points of the sources
\begin{equation}\label{gf-SW6}
S_x:=-[0\; h(\partial_x\varphi+c_fu+\omega v) \;0]^T\,,\;\;
S_y:=-[0\; 0\; h(\partial_y\varphi+c_fv-\omega u )]^T
\end{equation}
The entropy correction in 2D is evaluated with natural extensions of the one dimensional formulas, see e.g.  \cite{gaburro2023high} and references therein for details.

\section{Numerical results} \label{sec:Numerics}

%  \begin{itemize}
%\item  Bathymetry  
%\item   Friction 
%\item Coriolis 
%\item 2D examples and applications
%\end{itemize} 

We have thoroughly tested the scheme proposed to verify all the theoretical properties presented in the previous sections, as well as to
investigate its application to  resolution of complex solutions  and to test its robustness. 
We also provide results on 2D flows which show the advantage of  the  global flux quadrature, despite of the simplicity of the extension proposed.
%
%In the following section, we test the scheme discussed in the previous section's for shallow water equations for flows near equilibrium in one and two dimensions. At first, we test the order of convergence for an analytical equilibrium as given in Corollary 4.5 for the well-balanced scheme by calculating error using all nodes and end nodes within the cell and also compare it to the non-well-balanced scheme. Next we compare the three schemes-- non-well-balanced scheme, well-balanced scheme with and without entropy correction for perturbation to the steady state solutions. Following this, we  show examples for  entropy evolution using these schemes to see the effect of entropy correction term. We also give examples to see the 2D extension of the schemes for flows near steady state and also for few standard test cases for 2D flows.

\subsection{Verification of Corollary \ref{superconvergence}}

We start by verifying the theoretical predictions on the superconvergence of the scheme. We performed two sets of tests.
The first consists in computing explicitly the discrete  steady solution (or global flux solution) $U^*_{\hh}$ by solving,
element by element, the nonlinear  algebraic equations  associated to the condition $G_{\hh}=G_0$. 
The second set of tests involves initializing the solution with the exact analytical steady state, and letting the scheme 
evolve it until a given finite time. In both cases we report the grid convergence for finite element approximations of degree $p=1, 2, 3, 4$.
For the second test, we compare the global flux quadrature to the classical approach in which the source term is integrated following a classical DGSEM implementation as 
$$
\int_K\phi_i S_{\hh} \rightarrow \mathcal{M} \mathbf{S}\,.
$$  
Only moving equilibria are considered, since the lake at rest is exactly preserved using \eqref{eq:S_node-mod}. Small perturbations of the latter
are studied later.

\paragraph{Frictionless one dimensional equilibria} We start from the  first family of equilibria  from section \S2.1,  analytically defined by \eqref{steady0}.
We compute solutions for two classical smooth states: a sub-critical one with $(q_0,E_0)= (4.42\text{m}^2/\text{s},22.05535\text{m}^2/\text{s}^2)$;
 a super-critical one  for which  $(q_0,E_0)= (4.42\text{m}^2/\text{s},28.8971\text{m}^2/\text{s}^2)$  (the gravity constant is set to $g=9.80665\text{m}/\text{s}^2$ here). 
The computational domain is $25$m long, \revp{$x\in [0,25]$}, and the bottom topography  is taken as (we omit the units to avoid  clutter)
\begin{align}
	b(x)=\begin{cases}
		0.2-0.05(x-10)^2 & x>8 \text{ and }x<12\\
		0&  \text{else}
	\end{cases}\label{bathymetry_1}\,.
\end{align}

\begin{figure}[H]
%	%\centering	
	\subfigure[All nodes]{\includegraphics[width=0.5\textwidth]{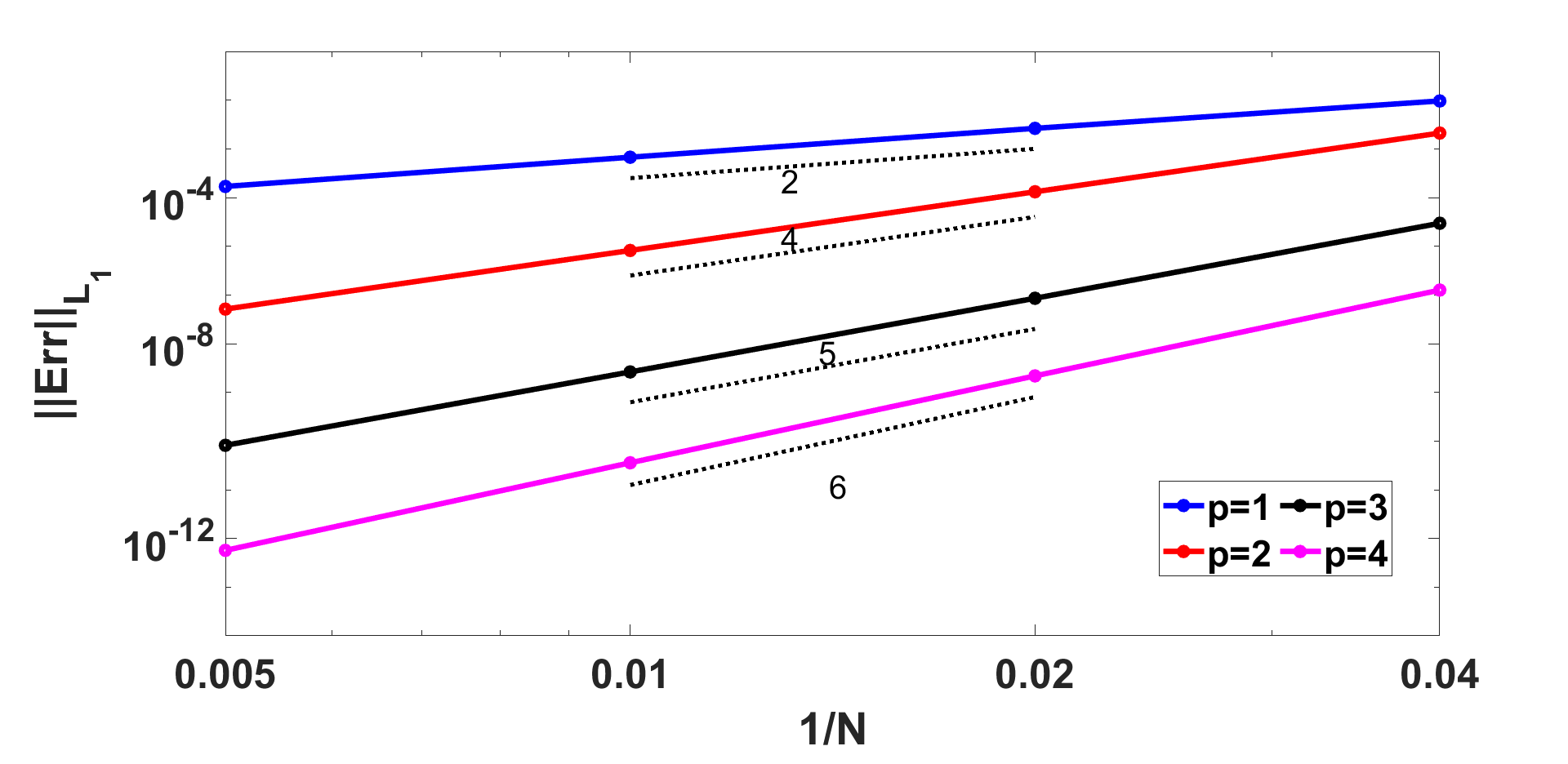}}	
	\subfigure[End nodes]{\includegraphics[width=0.5\textwidth]{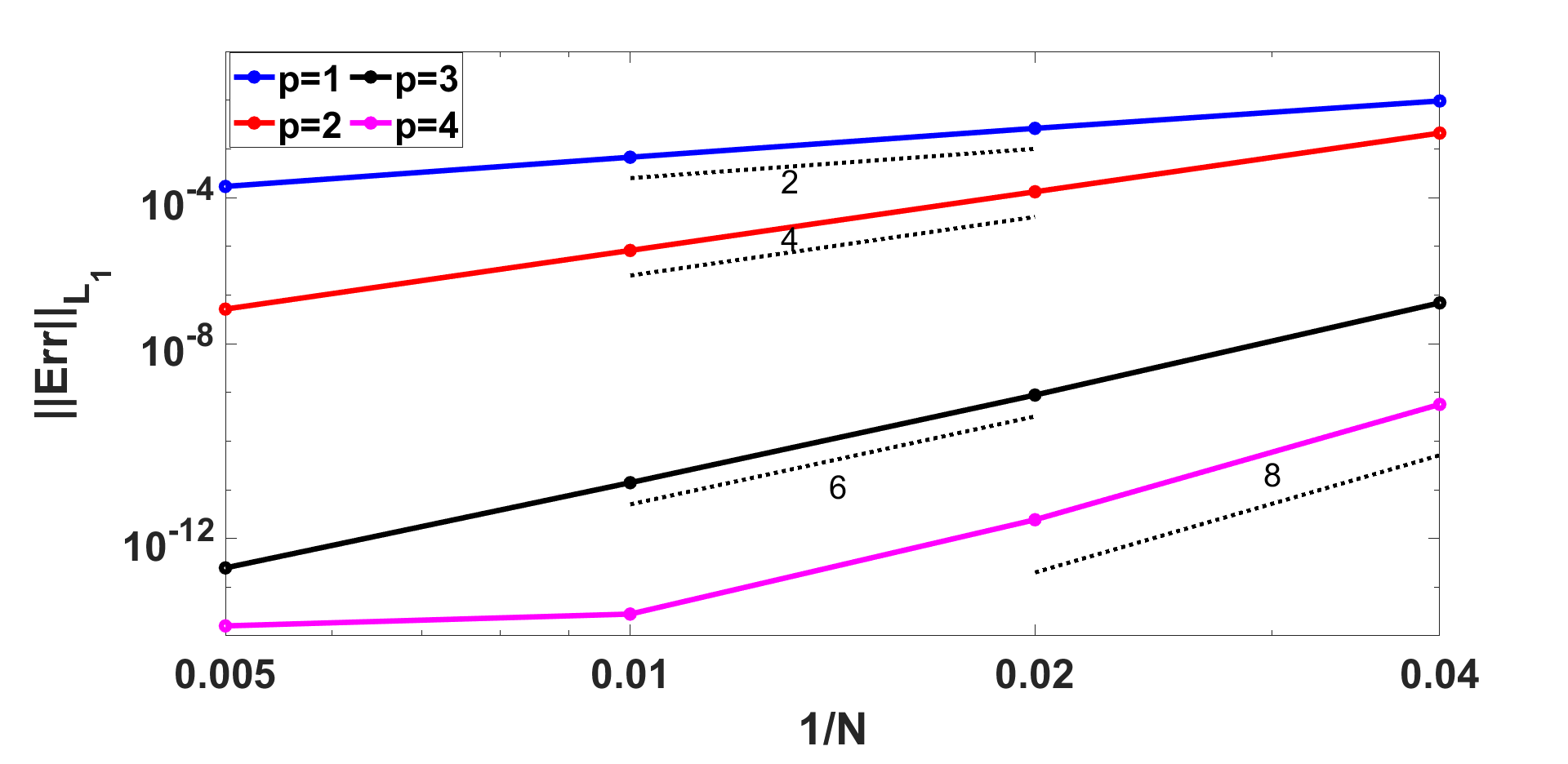}}
	\caption{Frictionless one dimensional  flow: sub-critical case. Error of the  global flux solution for polynomial spaces of degree $p=1,2,3,4$, measured with all nodes (left) and only  end nodes (right).}
	\label{Fig:conv_allP_subcrit_ic}
\end{figure}

\begin{figure}[H]
%	%\centering	
	\subfigure[All nodes]{\includegraphics[width=0.5\textwidth]{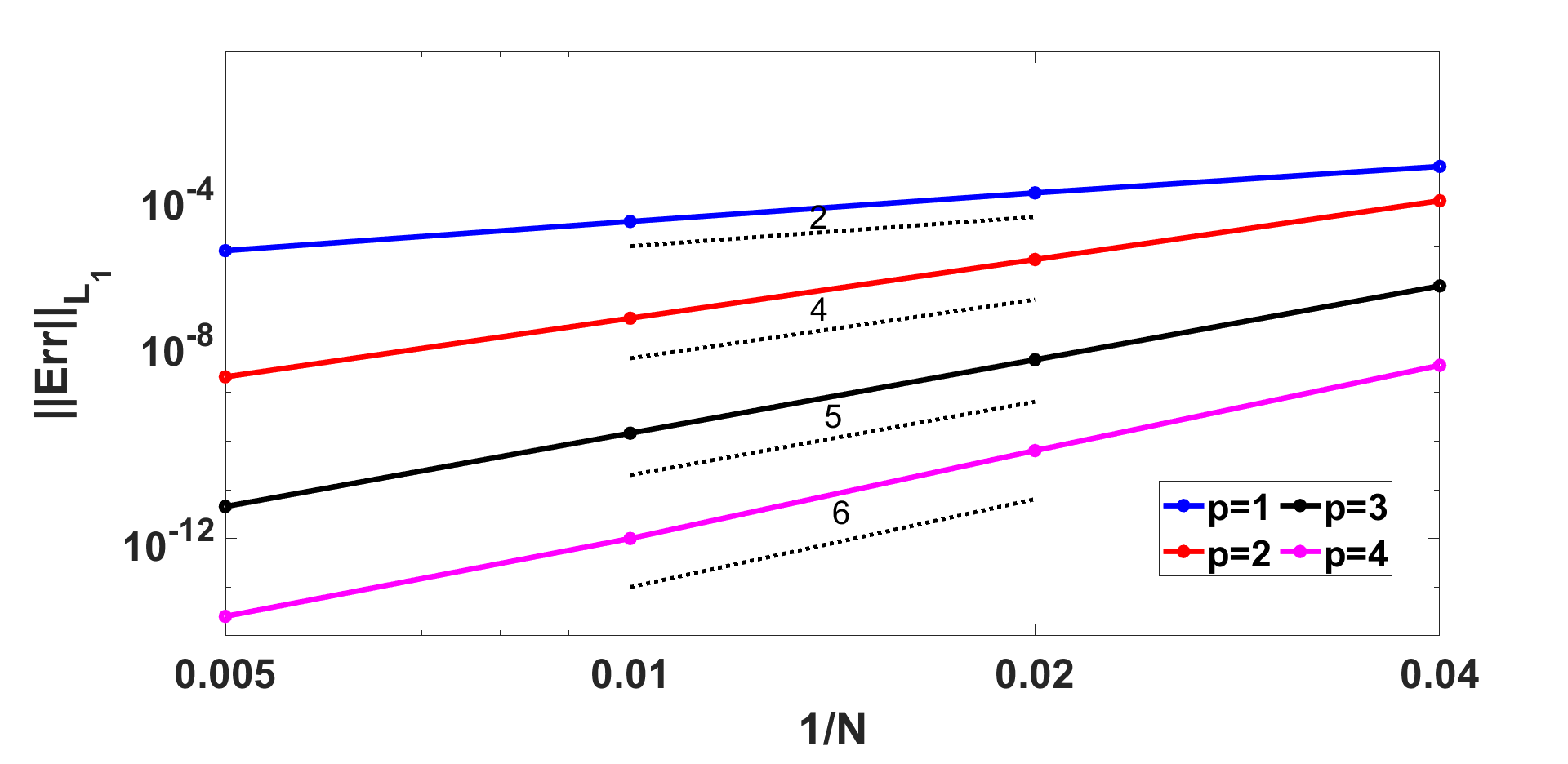}}	
	\subfigure[End nodes]{\includegraphics[width=0.5\textwidth]{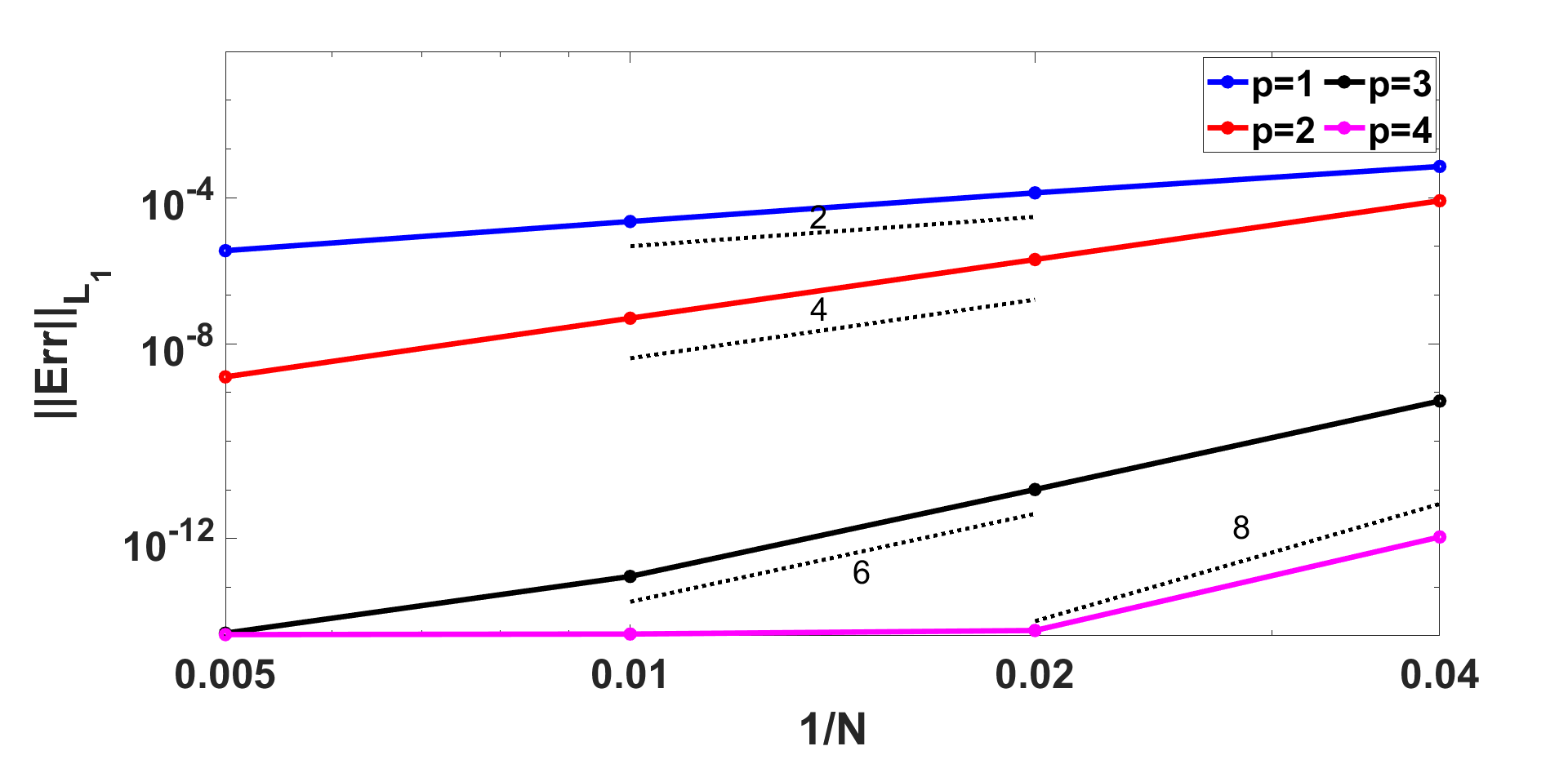}}
		\caption{Frictionless one dimensional  flow: super-critical case. Error of the  global flux solution for polynomial spaces of degree $p=1,2,3,4$, measured with all nodes (left) and only  end nodes (right).}
	\label{Fig:conv_allP_supercrit_ic}
\end{figure}

% 
%Consider the 1D shallow water equations with bottom topography, $b$ and without friction $c_f=0$ and Coriolis $\omega=0$. The analytical equilibrium for such a flow is given by,
%\begin{align}
%	hu=\text{constant},\quad \frac{u^2}{2}+g(h+b)=\text{constant}.\label{an_equilibrium_bathymetry}
%\end{align}
As already mentioned,  first we study the difference between the discrete steady state corresponding to constant global fluxes (global flux solution) and the analytical one. \revp{The constant global flux solution $V_0$ is obtained such that we get the state $(q_0,E_0)$ on the left boundary of the domain, i.e. $V_0=(4.42\text{m}^2/\text{s}, 29.3815\text{m}^3/\text{s}^2)$ for sub-critical flow and  $V_0=(4.42\text{m}^2/\text{s},31.7365\text{m}^3/\text{s}^2)$ for super-critical flow. The boundary conditions are also such that the solution satisfies the respective steady state on the boundary. }
Following the predictions of \revp{Corollary \ref{superconvergence}}, we compute two versions of the discrete  $L_1$-error 
$$
\|\text{Err}\|_{L_1} :=\dfrac{1}{M}\sum\limits_{j=1,M}|U^*_j -U^e(x_j)|\,
$$
the first with $j$ running over all the degrees of freedom of all the elements, the second running over only the element end-points.
 The logarithmic plots for $1/N$-$L_1$ error for $h$ with $p=1,2,3,4$ and $N=25,50,100,200$ are as shown in figure \ref{Fig:conv_allP_subcrit_ic} and figure \ref{Fig:conv_allP_supercrit_ic} for sub- and super-critical case respectively.

%<<<<<<< HEAD
%\begin{figure}[H]
%	%\centering	
%	\subfigure[All nodes]{\includegraphics[width=0.5\textwidth]{oneD/Coriolis/error_coriolis_h_all_an_gf}}	
%	\subfigure[End nodes]{\includegraphics[width=0.5\textwidth]{oneD/Coriolis/error_coriolis_h_end_an_gf}}
%	
%	
%	\subfigure[All nodes]{\includegraphics[width=0.5\textwidth]{oneD/Coriolis/error_coriolis_v_all_an_gf}}	
%	\subfigure[End nodes]{\includegraphics[width=0.5\textwidth]{oneD/Coriolis/error_coriolis_v_end_an_gf}}
%	
%	\caption{Comparison of difference between analytical and global flux equilibrium in DG space with $p=1,2,3,4$ for shallow water equation with Coriolis measured with all nodes and end nodes in every cell}
%	\label{Fig:conv_allP_coriolis_ic}
%\end{figure}
%
%=======
We can see that the convergence rates from figure \ref{Fig:conv_allP_subcrit_ic} and figure \ref{Fig:conv_allP_supercrit_ic} are as predicted in \revp{Corollary \ref{superconvergence}}, i.e. by using all nodes to measure the $L_1$ error the convergence is of order $p+2$, whereas by using only the end nodes, the order of convergence is $2p$. It is to be noted that with $p=1$, nodes for Gauss Lobatto are at the ends of the cells, and hence the order of convergence is $2$, even in figure \ref{Fig:conv_allP_subcrit_ic}(a) and figure \ref{Fig:conv_allP_supercrit_ic}(a). It can also be noted that with $p=4$ the error at endpoints  reaches  machine accuracy on coarse meshes, which of course affects the  slopes which appear somewhat reduced. \\
%>>>>>>> 8633838fbb44d758f3df865c7cef7dc43c507d06

\begin{figure}[H]
\centering\subfigure[WB, end nodes]{\includegraphics[width=0.5\textwidth]{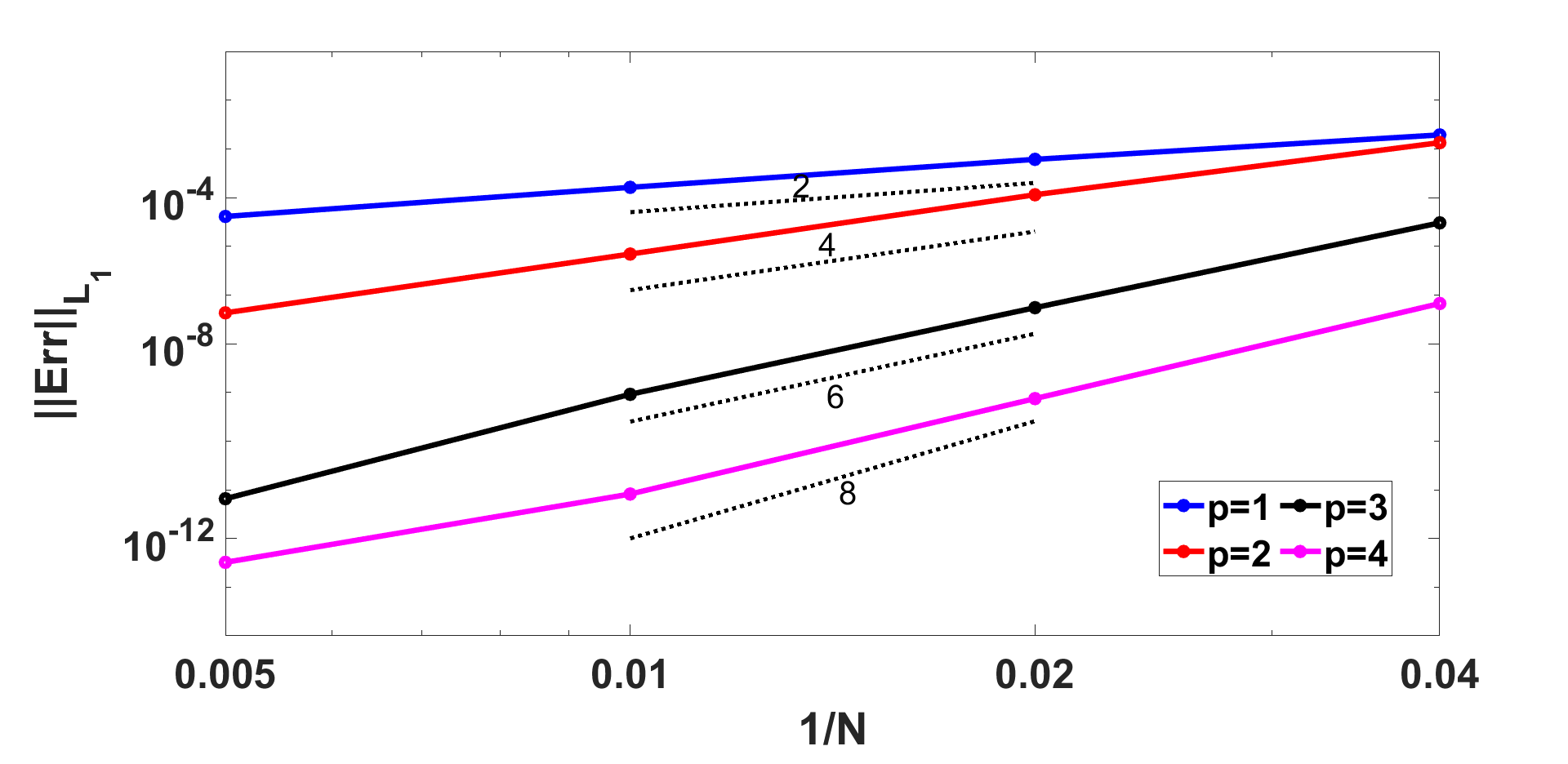}}	
	
	\subfigure[WB, all nodes]{\includegraphics[width=0.5\textwidth]{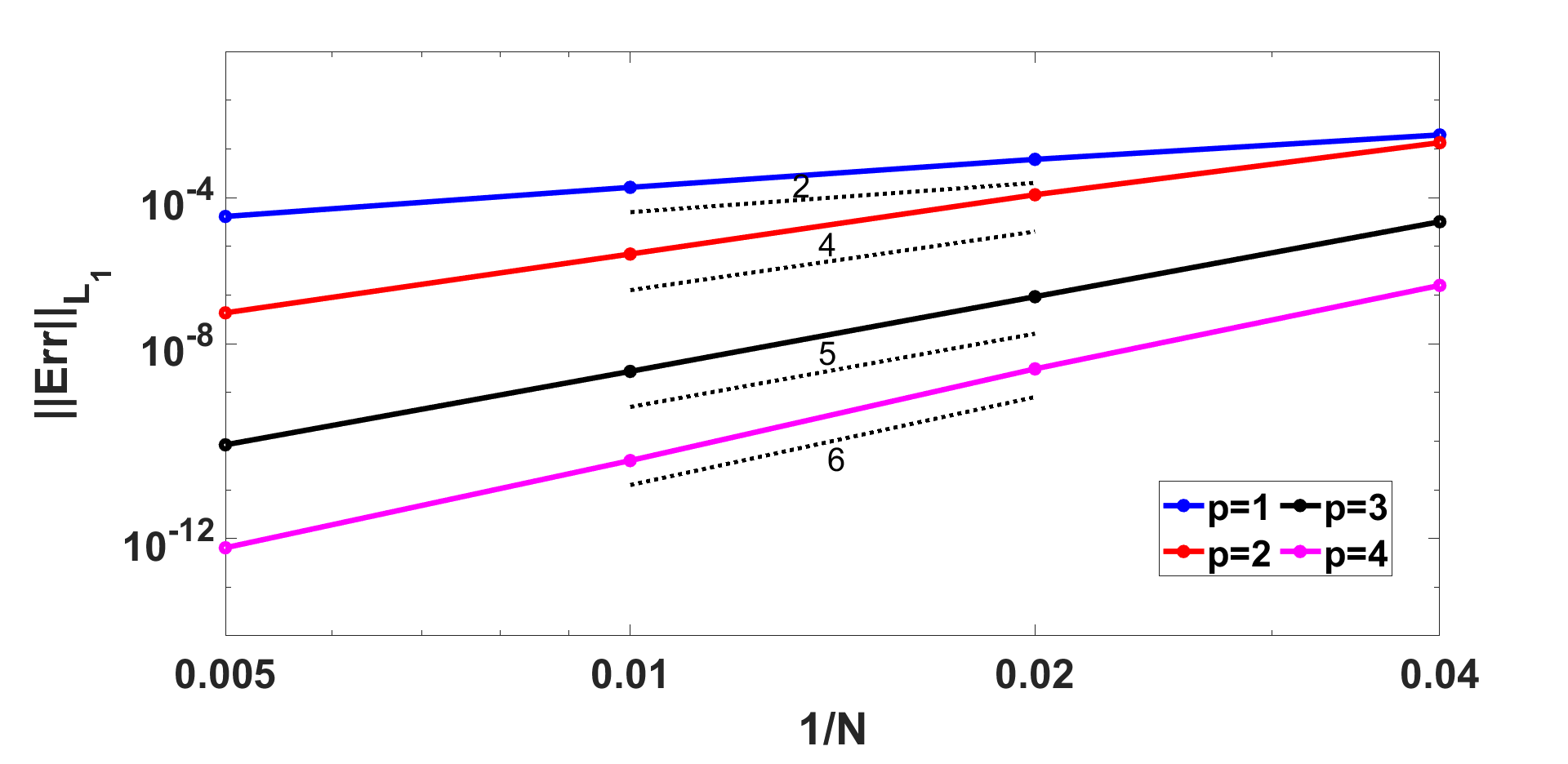}}\subfigure[NWB]{\includegraphics[width=0.5\textwidth]{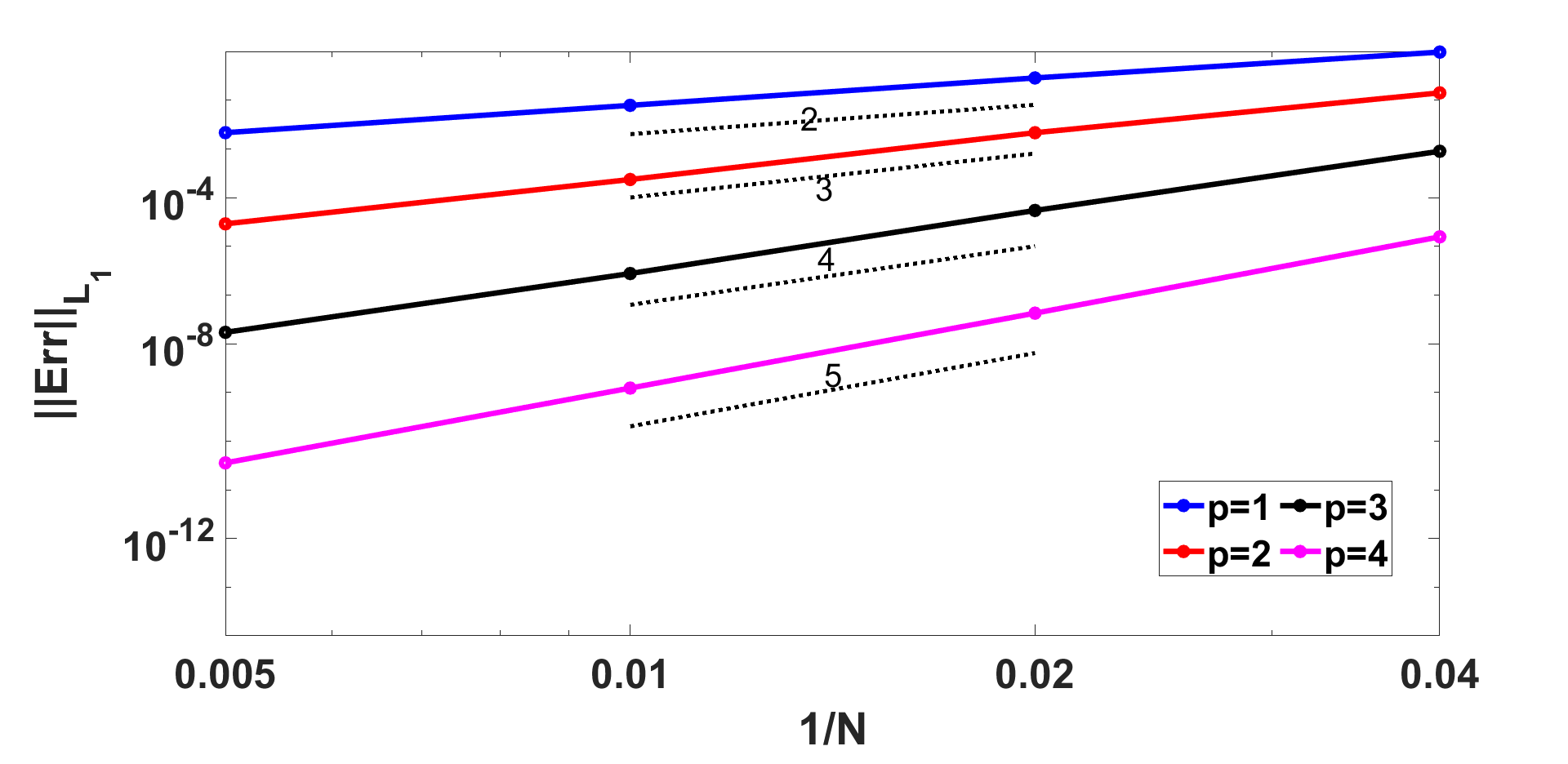}}
	\caption{Frictionless one dimensional  flow: sub-critical case. Error at finite time  for polynomial spaces of degree $p=1,2,3,4$.
	Top: end-nodes error. Bottom: global flux quadrature  (left) and non-well-balanced method (right).}
	\label{Fig:conv_allP_subcrit}
\end{figure}

\begin{figure}[H]
	\centering\subfigure[WB, end nodes]{\includegraphics[width=0.5\textwidth]{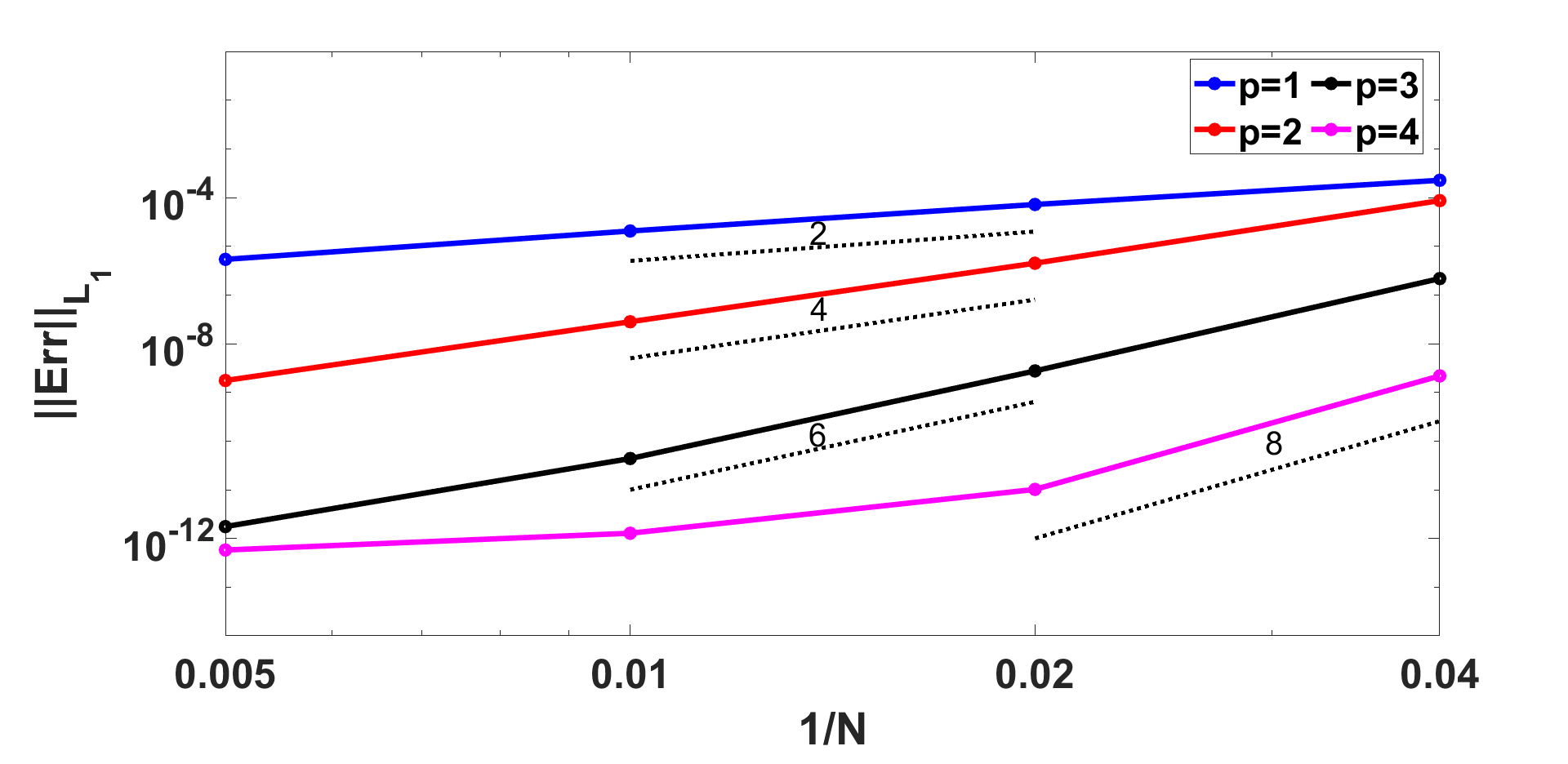}}	
		
	\subfigure[WB, all nodes]{\includegraphics[width=0.5\textwidth]{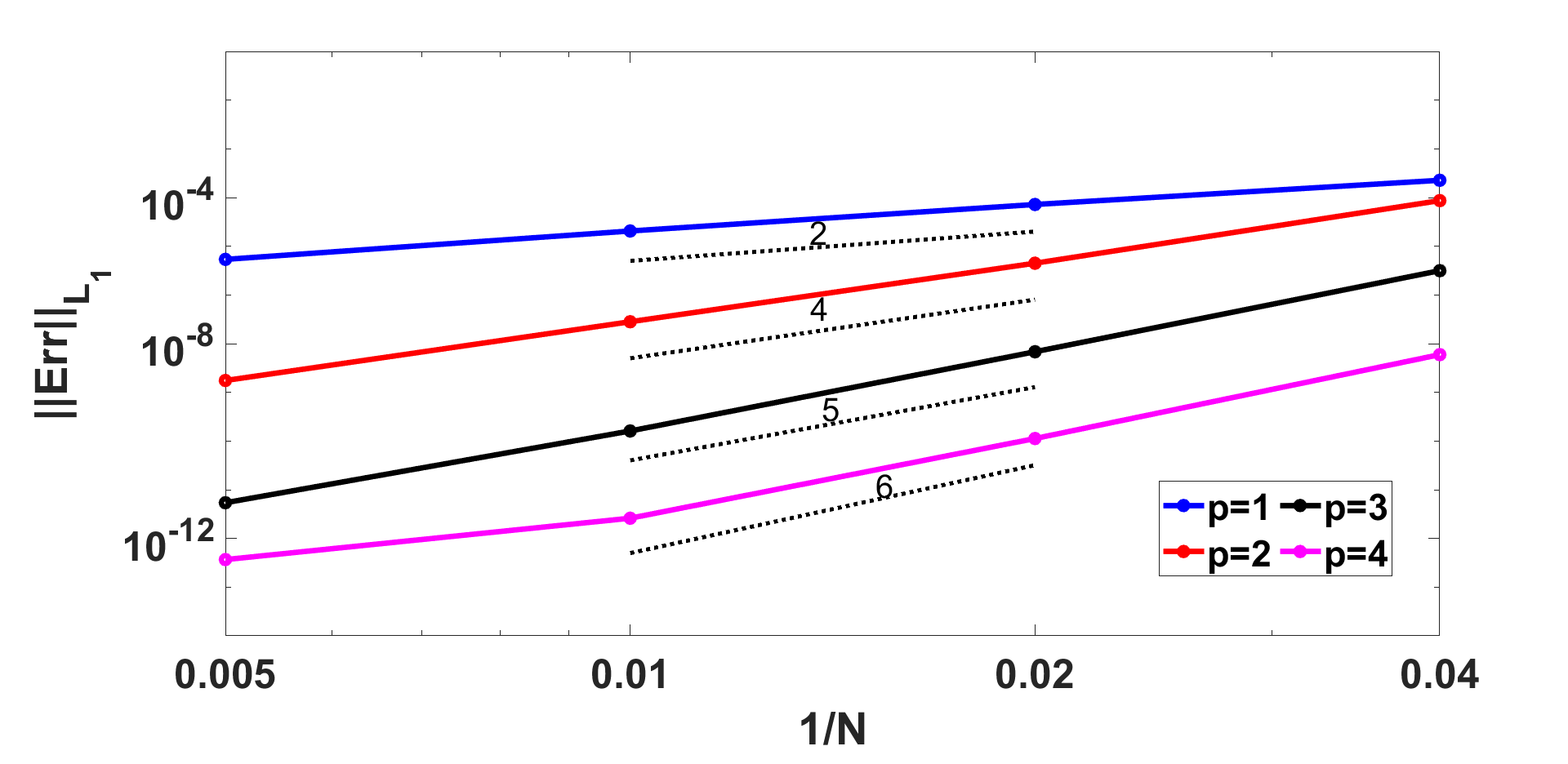}}\subfigure[NWB]{\includegraphics[width=0.5\textwidth]{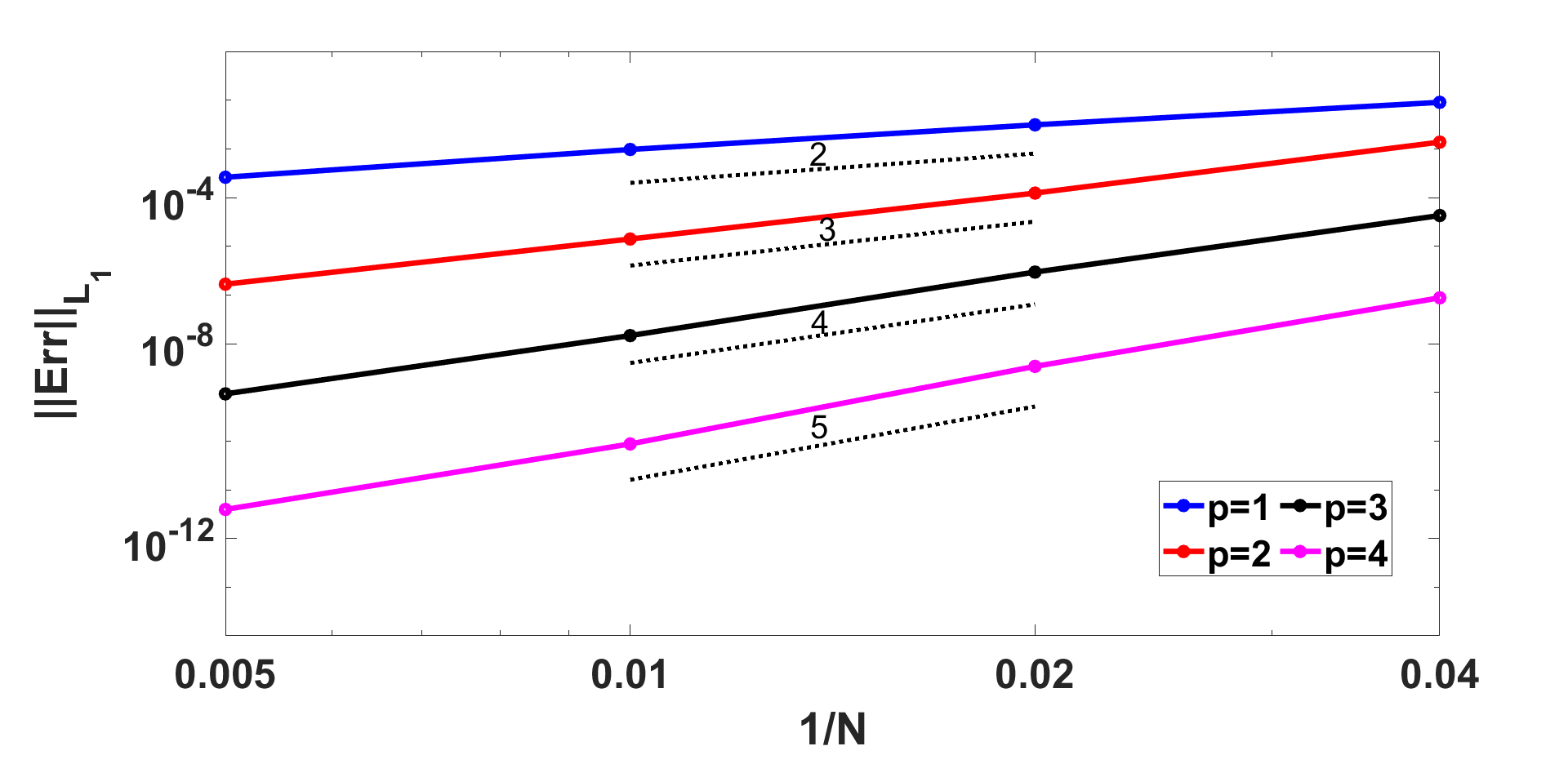}}
	\caption{Frictionless one dimensional  flow: super-critical case. Error at finite time  for polynomial spaces of degree $p=1,2,3,4$.
	Top: end-nodes error for the global flux quadrature method. Bottom: total error for the global flux   (left) and non-well-balanced method (right).}
	\label{Fig:conv_allP_supercrit}
\end{figure}

Next, we initialize the solution at degrees of freedom using the analytical value, and we    run both the global flux quadrature based   and non-well-balanced schemes 
up to $T=2$s.  \revPO{Note that this is an entirely different exercise, as we now let the scheme perturb the exact equilibrium,
and we measure that the magnitude of this error remains within the accuracy  foreseen by Corollary \ref{superconvergence}.}
 The logarithmic plots for $1/N$- $L_1$ error with $p=1,2,3,4$  and $N=25,50,100,200$ with error for well-balanced schemes measured using both all and end nodes and non-well-balanced scheme with all-nodes is as shown in figures \ref{Fig:conv_allP_subcrit} and  \ref{Fig:conv_allP_supercrit}.

The results show that, without computing the global flux solution, the proposed method still delivers  superconvergent results following the  predictions of \revp{Corollary \ref{superconvergence}}.
Moreover, we can see that the gain in accuracy  with respect to  the analytical solution on a given mesh compared to the classical DGSEM implementation is of more than 
two orders of magnitude. This brings the error of higher order approximation very rapidly to very low values, allowing to resolve very small perturbations as we 
will see shortly.

\begin{figure}[H]
	%\centering	
	\subfigure[$||E_h||_{L_1}$ All nodes]{\includegraphics[width=0.5\textwidth]{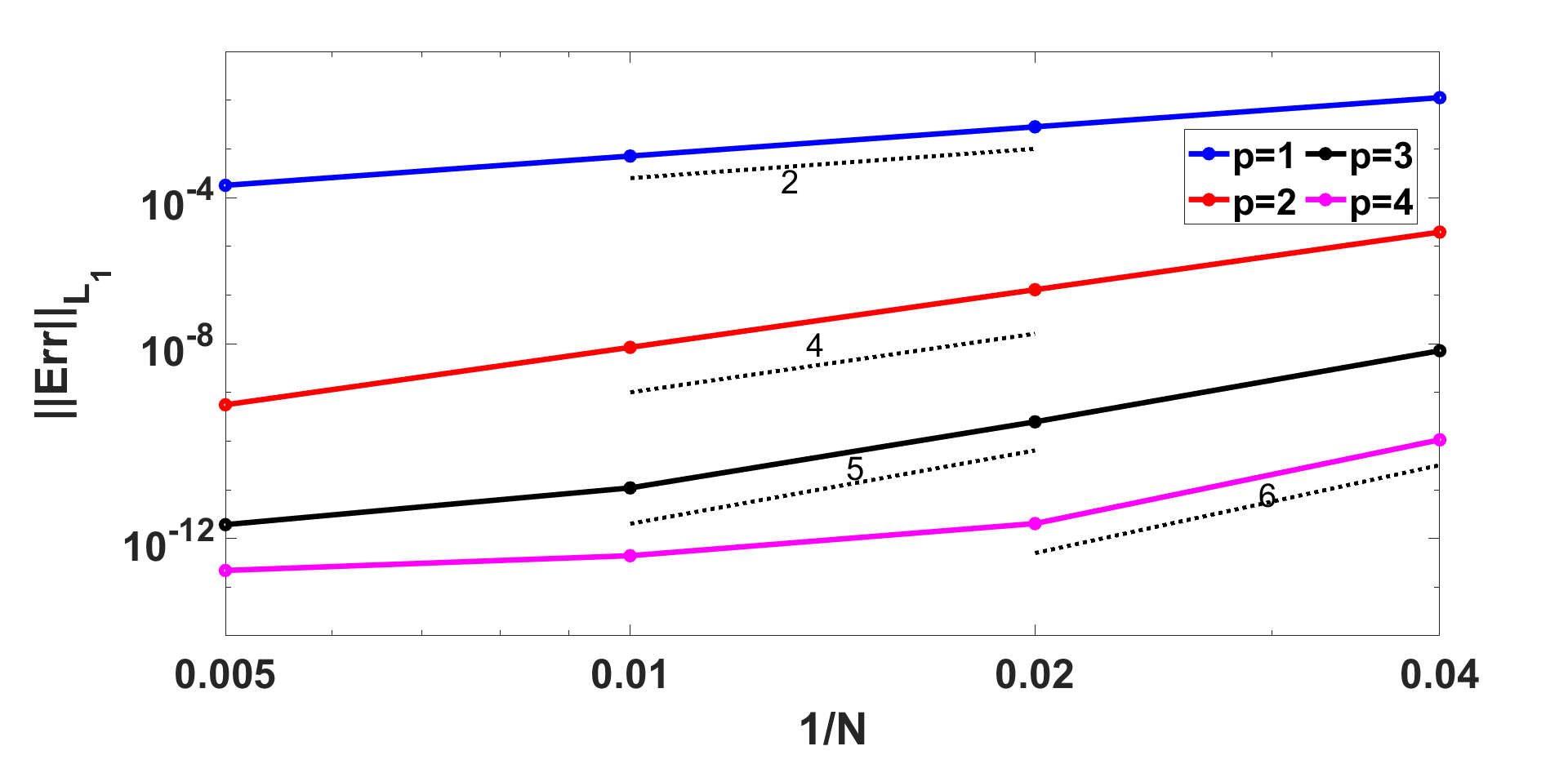}}	
	\subfigure[$||E_h||_{L_1}$ End nodes]{\includegraphics[width=0.5\textwidth]{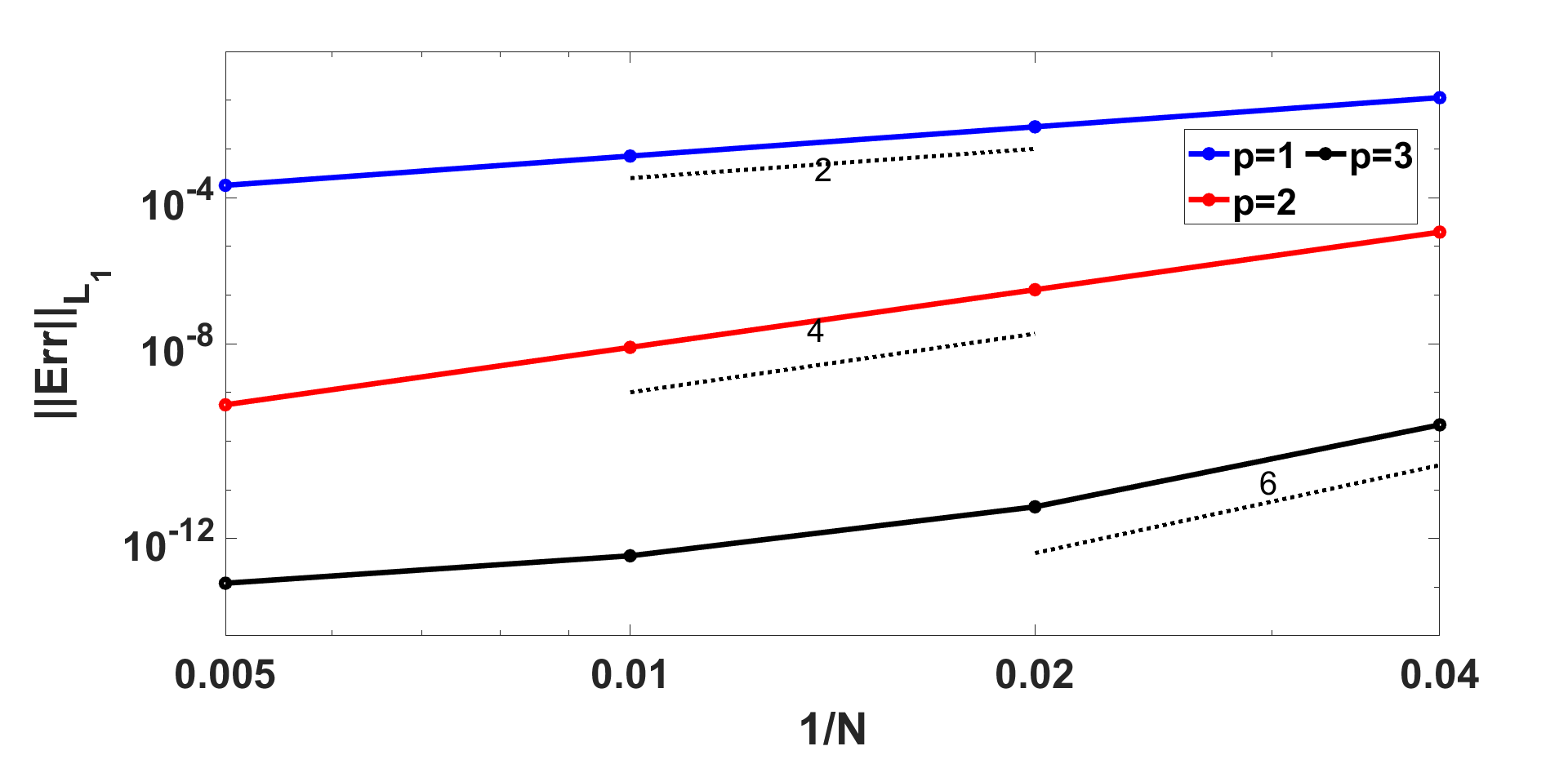}}
	\subfigure[$||E_{hv}||_{L_1}$ All nodes]{\includegraphics[width=0.5\textwidth]{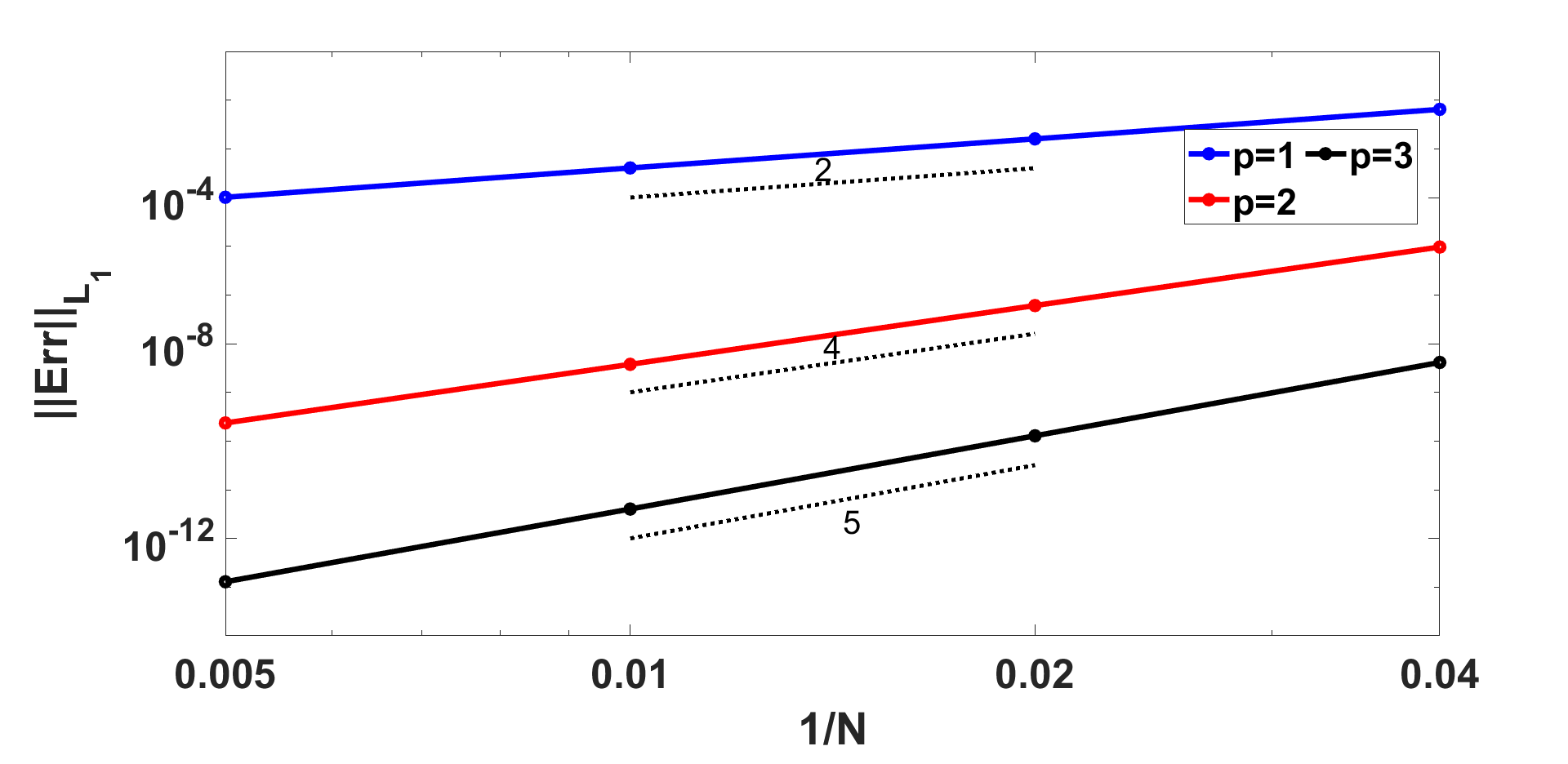}}	
	\subfigure[$||E_{hv}||_{L_1}$ End nodes]{\includegraphics[width=0.5\textwidth]{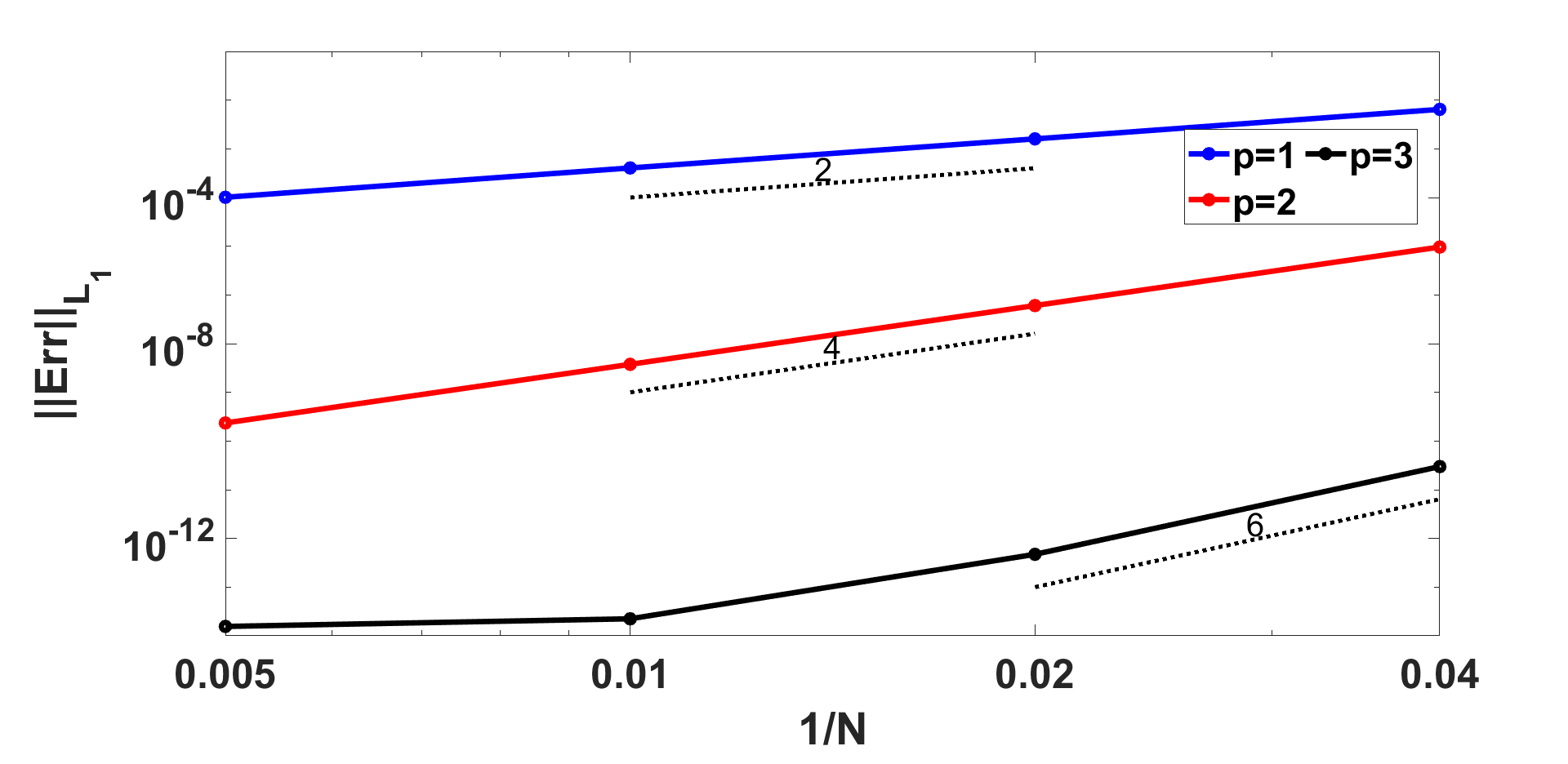}}
	
	\caption{Pseudo-1D  equilibria with Coriolis effects. Error convergence of the global flux solution with $p=1,2,3,4$ for depth (top) and transverse  momentum (bottom). Error at all the degrees of freedom (left), and only elements endpoints (right).}
	\label{Fig:conv_allP_coriolis_ic}
\end{figure}

\paragraph{Frictionless pseudo-one dimensional equilibria with Coriolis effects} We consider next the solutions defined by \eqref{steady3}.
We choose the manufactured state also used in \cite{desveaux21} given within  $x\in [0,1]$m by  (here the gravity is set to $g=1\text{m}/\text{s}^2$ \revp{and Coriolis coefficient $\omega=1\text{s}^{-1}$}) (we omit units for simplicity):
\begin{align}
	h=e^{2x},\quad hu=1,\quad hv=-\omega x e^{2x}\label{eq:SWE_Coriolis_equilibrium_mov_an}
\end{align}
with bathymetry
$$
b(x)=-\frac{1}{2}\omega^2x^2-e^{2x}-\frac{1}{2}e^{-4x}
$$ 
%The global flux equilibrium is calculated using the analytical solution for $hu$ and $v$ and calculating $h$ such that the global flux is constant,
%\begin{align}
%	{hu^2}+\frac{h^2}{2}+\int(g h b_x-\omega hv) dx=\text{constant}.\label{eq:SWE_Coriolis_equilibrium_mov_gf}
%\end{align} 
 We perform  the same comparisons done for the previous case. \revp{The discrete  global flux steady state is calculated using the values of $h,hu,hv$ on the left boundary.} % at $x=0$.} 
 We report  in figure \ref{Fig:conv_allP_coriolis_ic}
the error convergence plots for the global flux solutions, using the discrete norms involving all the degrees of freedom, as well as only the endpoints.
For completeness we report both the solutions for $h$ and for the transverse  momentum which is not constant.
As before, the convergence rates  confirm \revp{Corollary \ref{superconvergence}}. For this case,  with $p=4$ the errors reach  machine accuracy extremely fast,
already on the coarsest mesh for most cases. These cases are not shown in the plots for clarity. The $p=3$ errors show a similar trend at end-nodes. 
This makes these schemes in practice almost exactly well-balanced with the analytical state.\\

\begin{figure}[H]
          \centering\subfigure[$||E_h||_{L_1}$ WB, all nodes]{\includegraphics[width=0.5\textwidth]{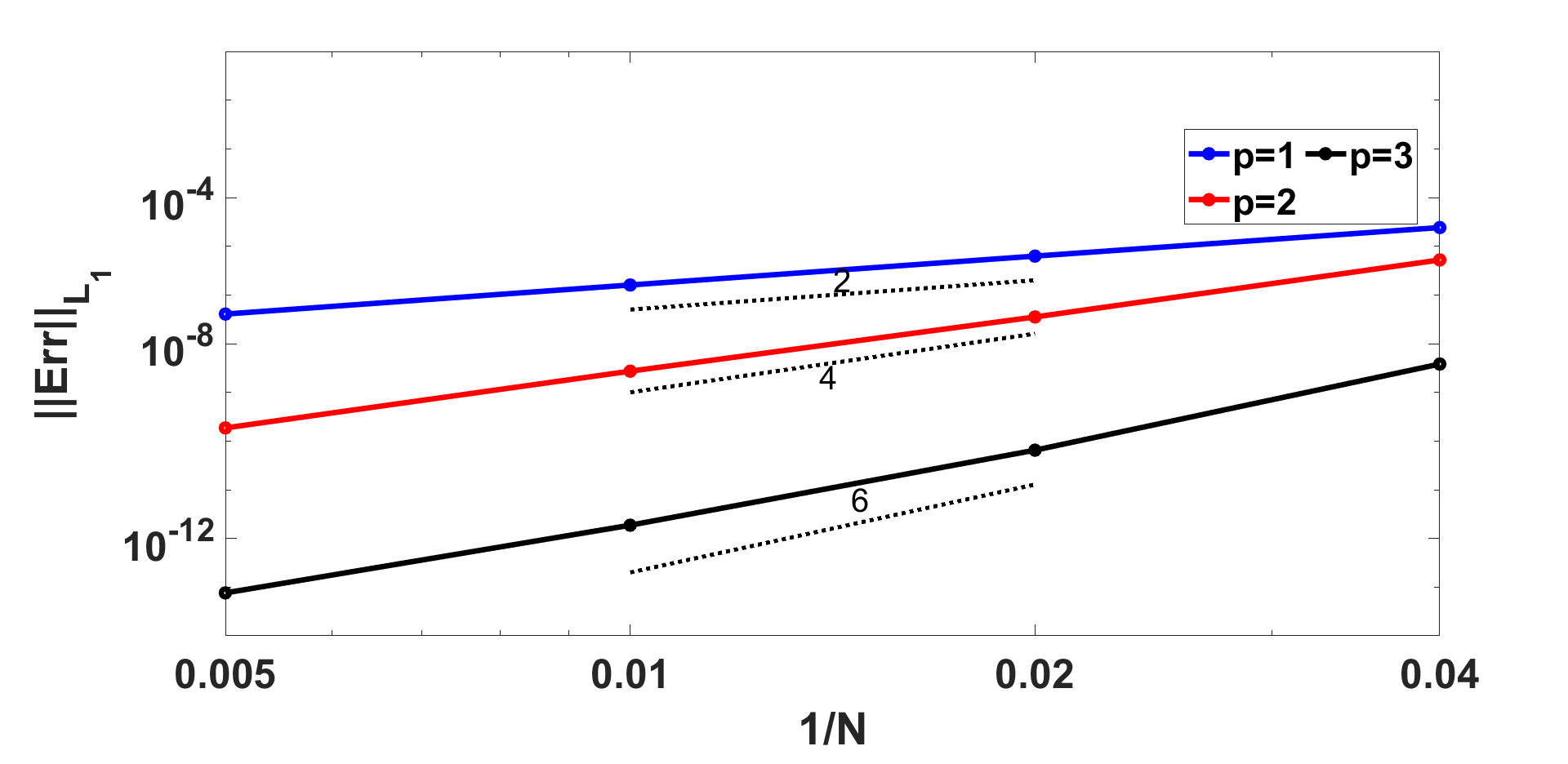}}
%	\subfigure[$||E_{hv}||_{L_1}$ WB, all nodes]{\includegraphics[width=0.5\textwidth]{oneD/Coriolis/error_coriolis_v_all_wb}}	
\subfigure[$||E_h||_{L_1}$ WB ]{\includegraphics[width=0.5\textwidth]{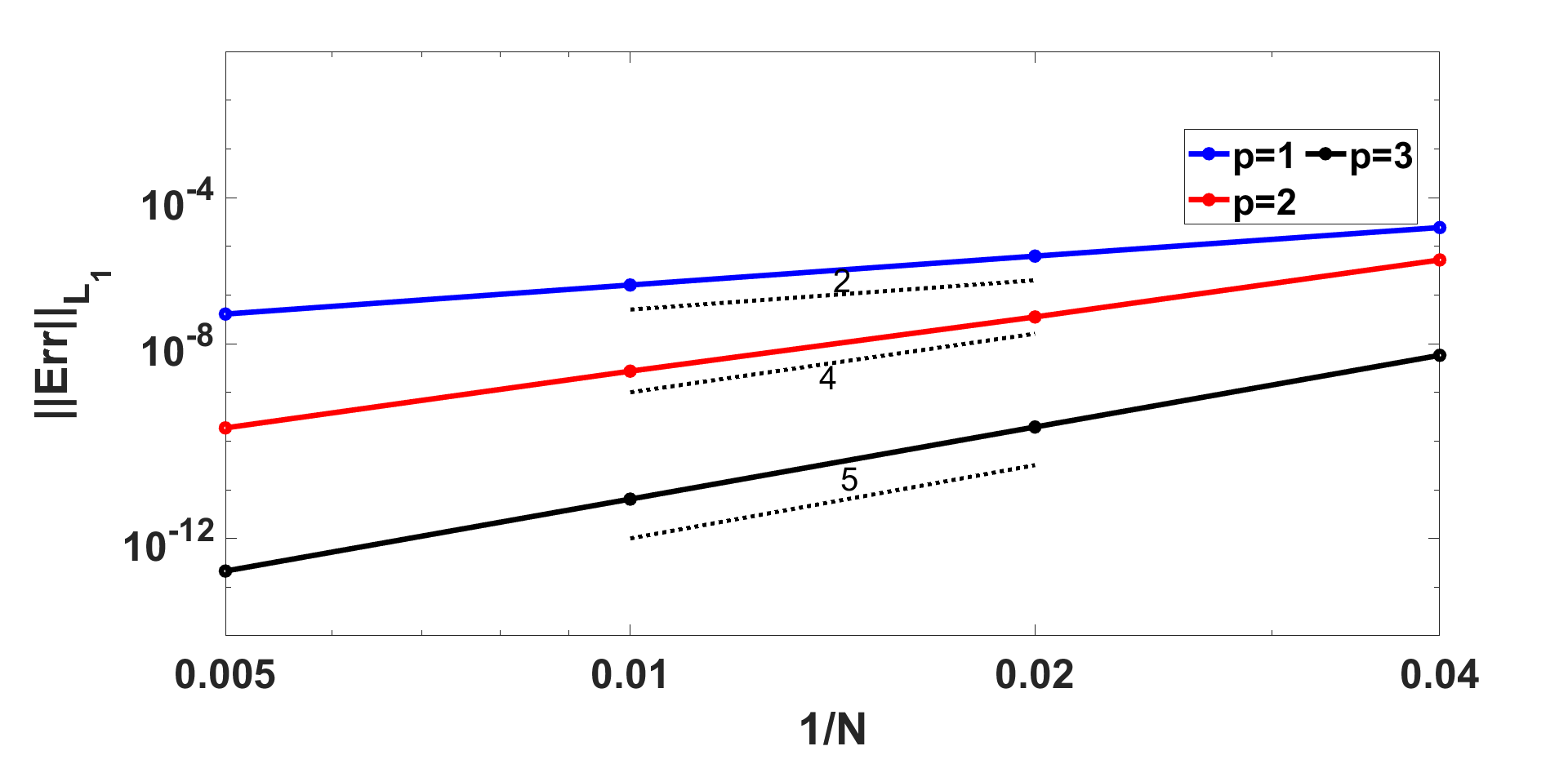}}\subfigure[$||E_h||_{L_1}$  NWB]{\includegraphics[width=0.5\textwidth]{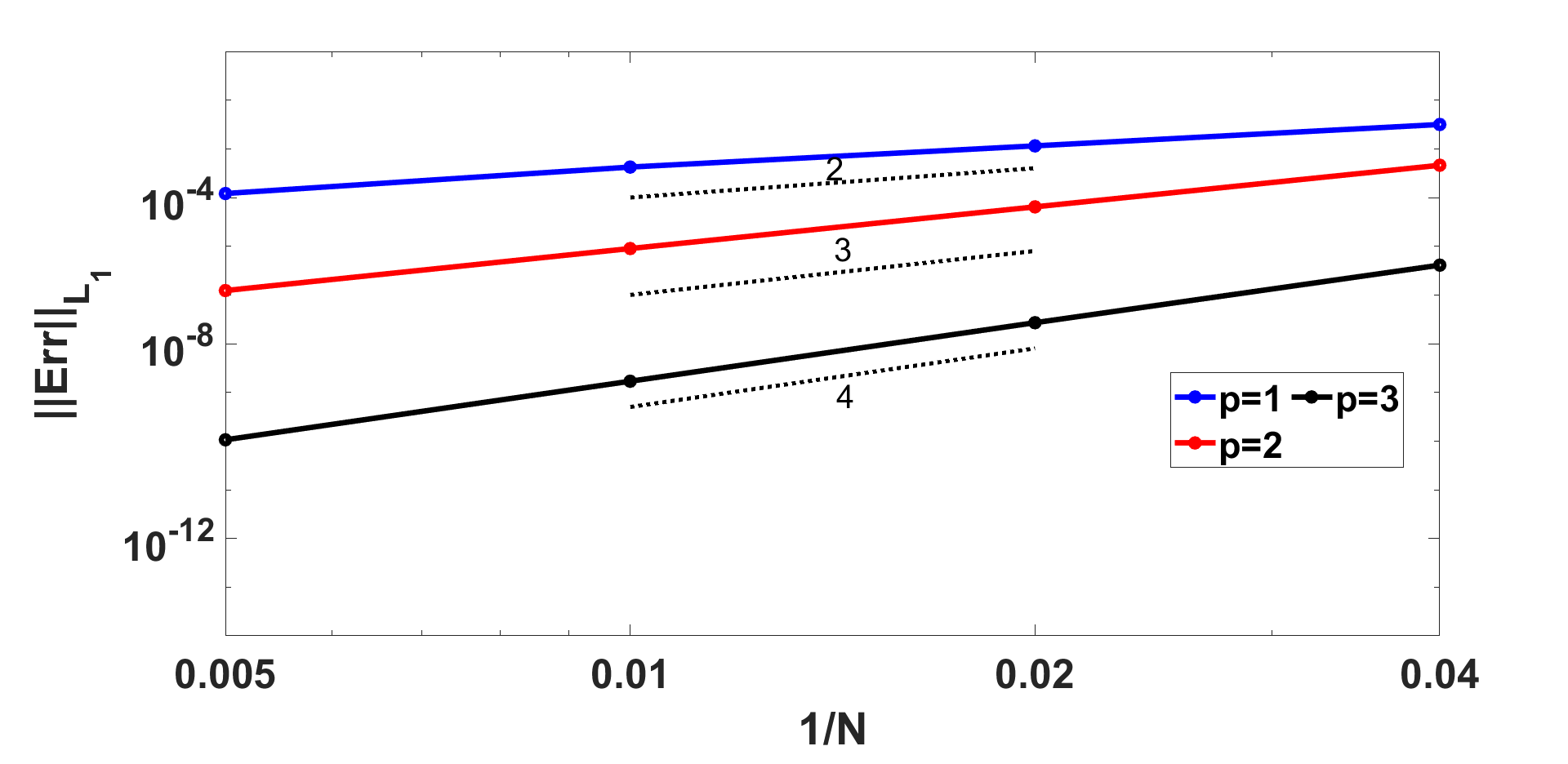}}	
%<<<<<<< HEAD	
%	\subfigure[$||E_{hv}||_{L_1}$, NWB]{\includegraphics[width=0.5\textwidth]{oneD/Coriolis/error_coriolis_v_all_nwb}}
	\caption{Pseudo-1D  equilibria with Coriolis effects. Error convergence of the global flux solution with $p=1,2,3,4$ for the depth.  
	  Top: End-nodes error for the global flux quadrature method. Bottom: Total error for the global flux   (left) and non-well-balanced method (right).}
%=======
%	\subfloat[WB, all nodes]{\includegraphics[width=0.33\textwidth]{oneD/Coriolis/error_coriolis_v_all_wb}}	
%	\subfloat[WB, end nodes]{\includegraphics[width=0.33\textwidth]{oneD/Coriolis/error_coriolis_v_end_wb}}	
%	\subfloat[$NWB$]{\includegraphics[width=0.33\textwidth]{oneD/Coriolis/error_coriolis_v_all_nwb}}
%	\caption{Comparison of error for $h$ and $hv$ using all nodes and end nodes at each cell and all nodes for non-WB with $p=1,2,3,4$ for super-critical flow}
%>>>>>>> 8633838fbb44d758f3df865c7cef7dc43c507d06
	\label{Fig:conv_allP_coriolis}
\end{figure}

\begin{figure}[H]
	\centering\subfigure[$||E_{hv}||_{L_1}$ WB  all nodes]{\includegraphics[width=0.5\textwidth]{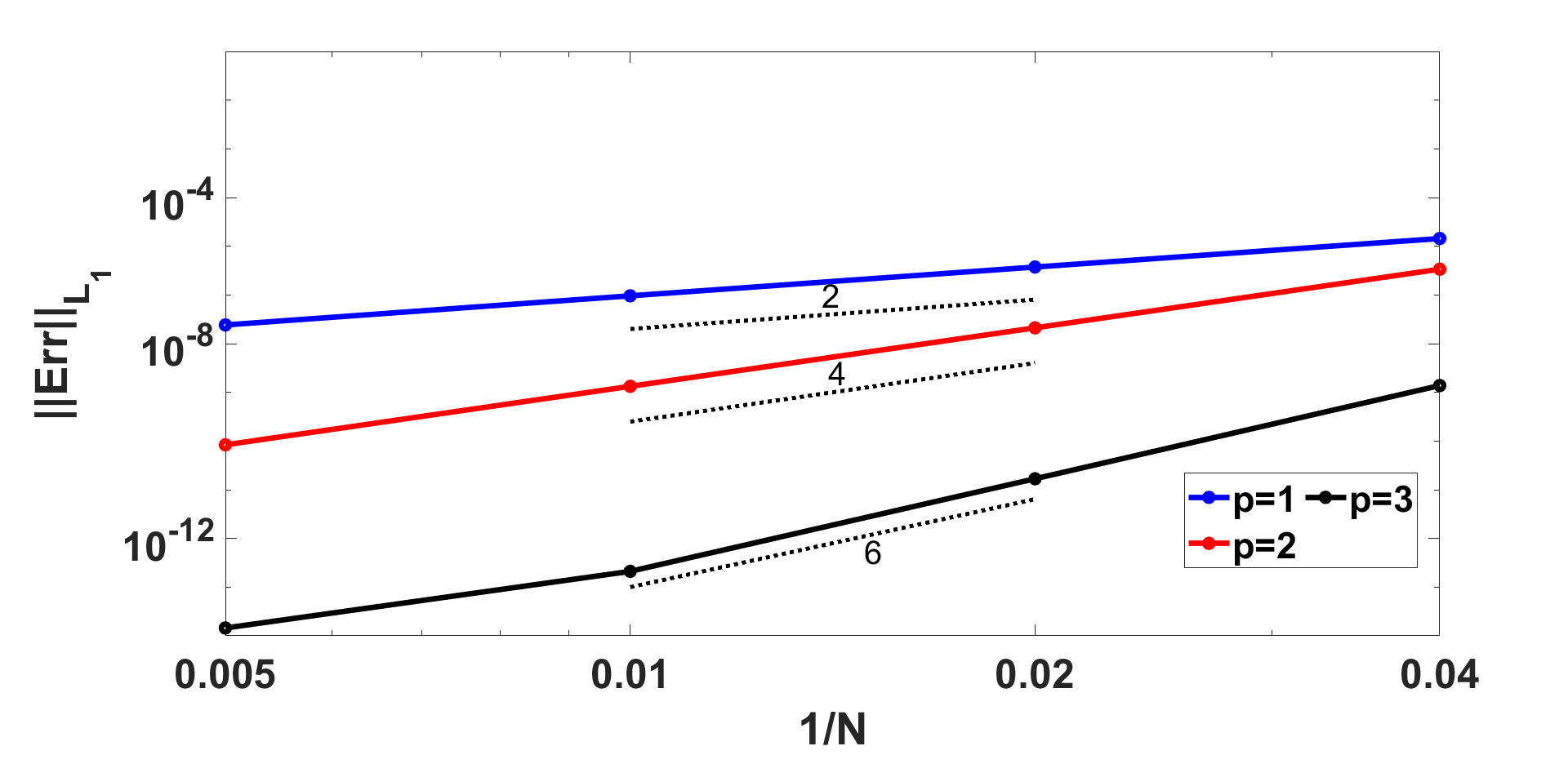}}	
%	\subfigure[$||E_h||_{L_1}$, WB, end nodes]{\includegraphics[width=0.5\textwidth]{oneD/Coriolis/error_coriolis_h_end_wb}}
	\subfigure[$||E_{hv}||_{L_1}$ WB all nodes]{\includegraphics[width=0.5\textwidth]{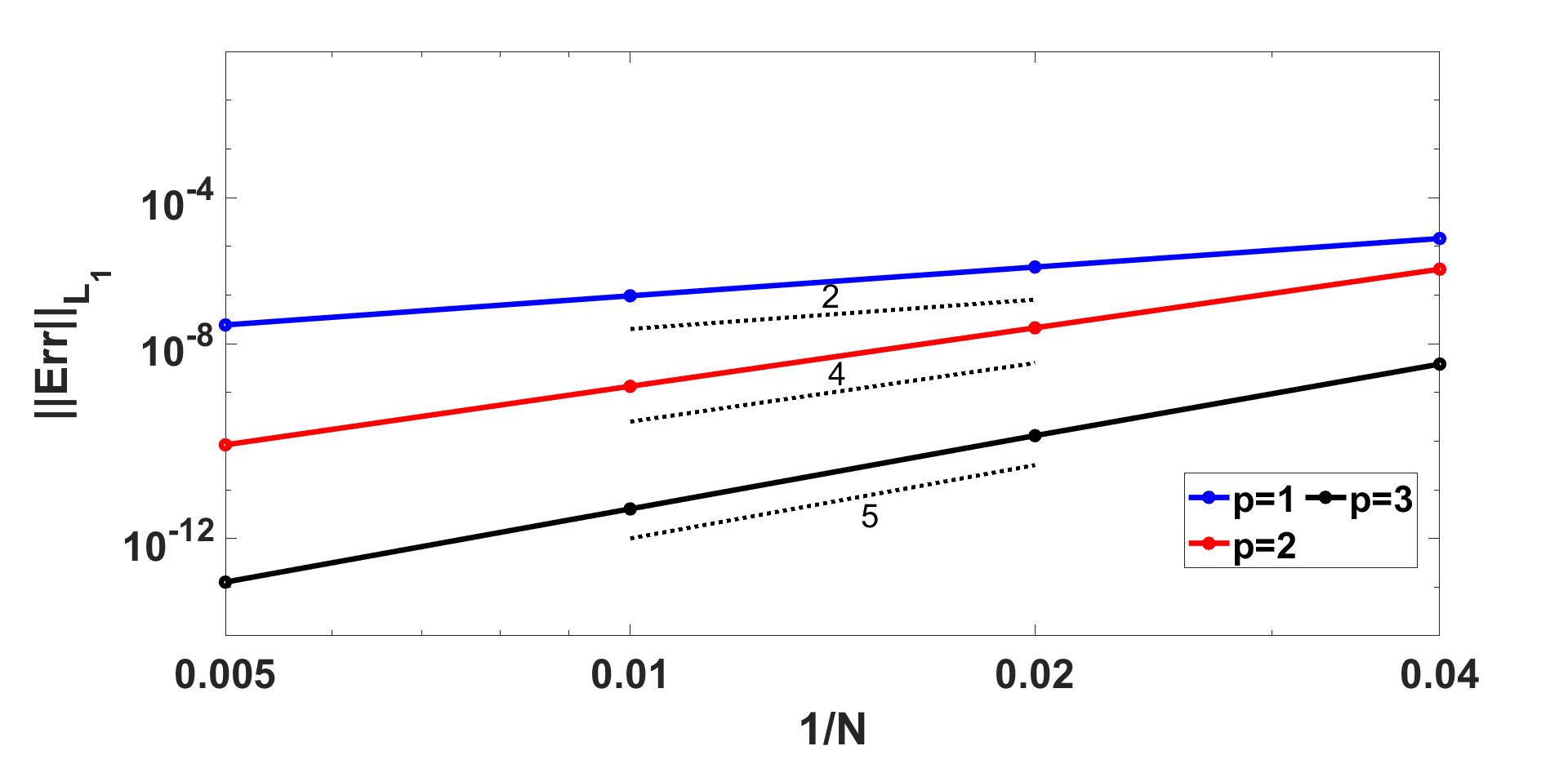}}\subfigure[$||E_{hv}||_{L_1}$  NWB]{\includegraphics[width=0.5\textwidth]{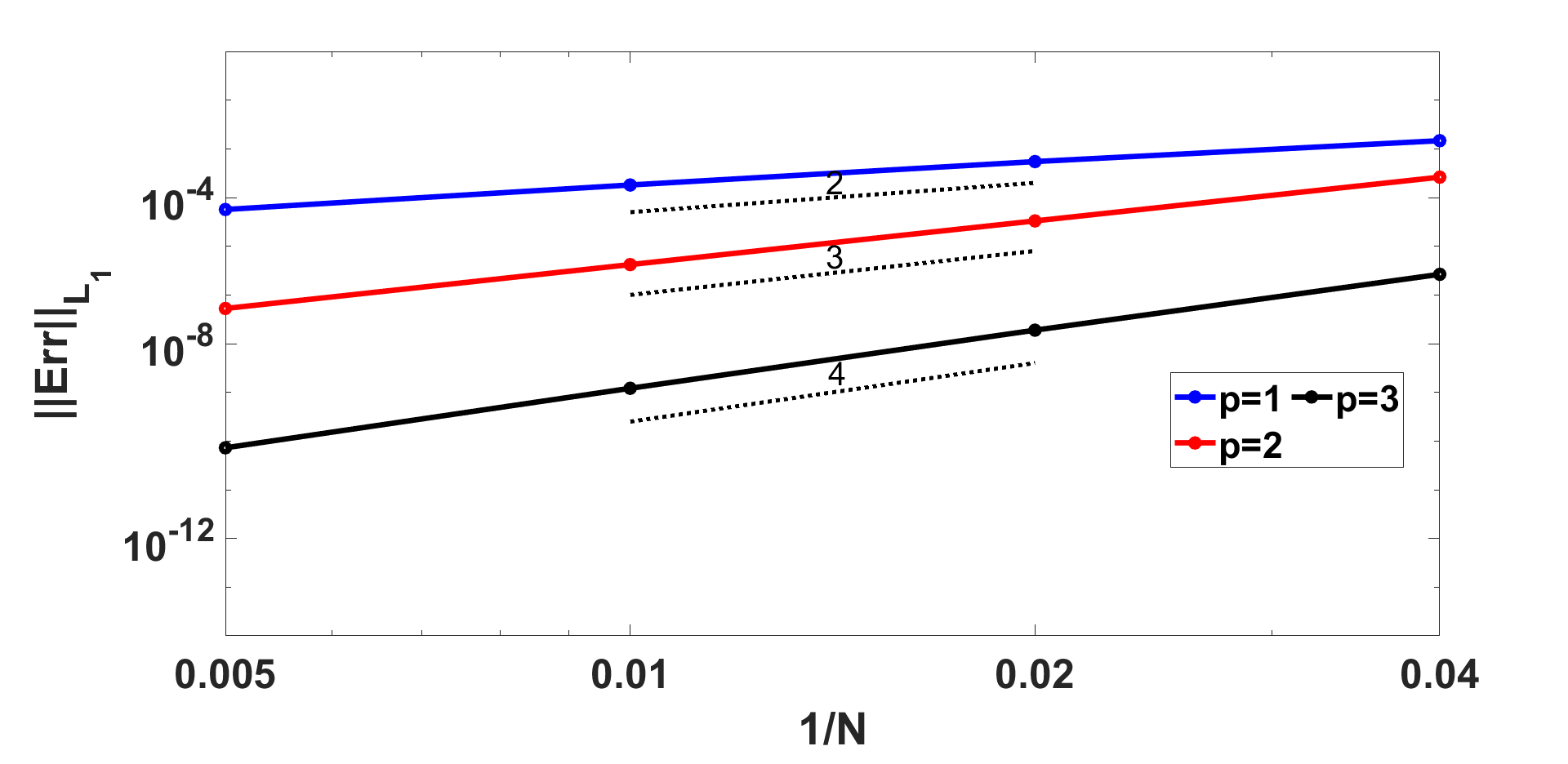}}
	\caption{Pseudo-1D  equilibria with Coriolis effects. Error convergence of the global flux solution with $p=1,2,3,4$ for the tranverse momentum.  
	  Top: End-nodes error for the global flux quadrature method. Bottom: Total error for the global flux   (left) and non-well-balanced method (right).}
%=======
%	\subfloat[WB, all nodes]{\includegraphics[width=0.33\textwidth]{oneD/Coriolis/error_coriolis_v_all_wb}}	
%	\subfloat[WB, end nodes]{\includegraphics[width=0.33\textwidth]{oneD/Coriolis/error_coriolis_v_end_wb}}	
%	\subfloat[$NWB$]{\includegraphics[width=0.33\textwidth]{oneD/Coriolis/error_coriolis_v_all_nwb}}
%	\caption{Comparison of error for $h$ and $hv$ using all nodes and end nodes at each cell and all nodes for non-WB with $p=1,2,3,4$ for super-critical flow}
%>>>>>>> 8633838fbb44d758f3df865c7cef7dc43c507d06
	\label{Fig:conv_allP_coriolis1}
\end{figure}

Next we run the well-balanced and non-well-balanced schemes with initial condition given by analytical equilibrium \eqref{eq:SWE_Coriolis_equilibrium_mov_an}. The $L_1$ error at $T=1$ with all nodes for both  schemes  is as shown in figure \ref{Fig:conv_allP_coriolis} for the depth, and in figure  \ref{Fig:conv_allP_coriolis1} for the transverse  momentum.
We only report the results for   $p=1,2,3$. We omit the $p=4$ results as for which the global flux quadrature method is 
within machine accuracy of the exact solution. We can see that the results once again verify the rates predicted by  \revp{Corollary \ref{superconvergence}}, with convergence order of    $p+2$ for the total error, and  $2p$  for the end points error when using the global flux quadrature.  As for the previous case, 
the comparison with the non-well-balanced DGSEM implementation shows a gain in error of several orders of magnitude.

\subsection{Perturbations of steady states}
\subsubsection{Tests setup and solution visualization}
\revPO{In the following sections we consider a classical exercise consisting in studying the evolution of  small perturbations of steady state equilibria. The initial condition for  all the tests is given by
\begin{equation}\label{perturbation}
\begin{split}
	h(0,x)=&h^*(x)+\xi e^{-\frac{(x-x_0)^2}{100}}\\
		u(0,x)=&u^*(x) 
		\end{split}
\end{equation}
where $h^*(x) $ and $u^*(x)$ are the steady distributions of water depth and speed.}\\ 

\revt{The definition of the last two plays a critical role. Very often  in literature one finds  reference to \emph{well prepared} initial conditions without
clearly specifying what exactly this preparation of the initial data consists of. Here we consider two cases. The first is the use of the analytical solution, if known. 
This choice is independent of the scheme, and moreover allows to test the schemes wrt the preservation of the exact equilibrium which is the ideal situation.}

\revt{Another possible choice is to use a well defined  and well behaved discrete approximation of the exact equilibrium. For example,
in the works \cite{Castro2020,math9151799,math10010015} the authors use as a natural choice the enhanced  discrete approximation
also employed to modify the polynomial reconstruction of the scheme. This initialization is done by construction of the discrete steady state of the well balanced
scheme proposed in the references. A similar procedure can be used in our case, and would consist in using the solution provided by the ODE integrator LobattoIIIA
of Proposition \ref{proposition_4} and Corollary \ref{superconvergence}. These are reasonable choices, however slightly lacking objectivity when comparing  different schemes, as they
involve the discrete steady solution of one scheme in particular, which is thus favoured by this choice. An example is reported in the appendix
to confirm this fact.  In practice, in the following tests $h^*(x)$  and $u^*(x)$ are  given by the exact analytical steady values.} \\

\revPO{Also note that in all  plots we visualize cell by cell the actual high order finite element polynomial  instead of some other interpolation of the
nodal values.   We feel that this representation is the most faithful to the nature of the high order discontinuous finite element approach of the paper. }

\subsubsection{Lake at rest}  
\begin{figure}[H]
%%<<<<<<< HEAD
%	\hspace{-0.6cm}\subfigure[$\eta=10^{-1}$]{\includegraphics[width=0.33\textwidth]{oneD/Bathymetry/height_lakerest_pert1zoom}}\hspace{-0.4cm}	
%	\subfigure[$\eta=10^{-3}$]{\includegraphics[width=0.33\textwidth]{oneD/Bathymetry/height_lakerest_pert3zoom}}	
%	\subfigure[$\eta=10^{-6}$]{\includegraphics[width=0.33\textwidth]{oneD/Bathymetry/height_lakerest_pert6zoom}}
%%=======
	\centering\subfigure[$\xi=10^{-1}m$]{\includegraphics[width=0.5\textwidth]{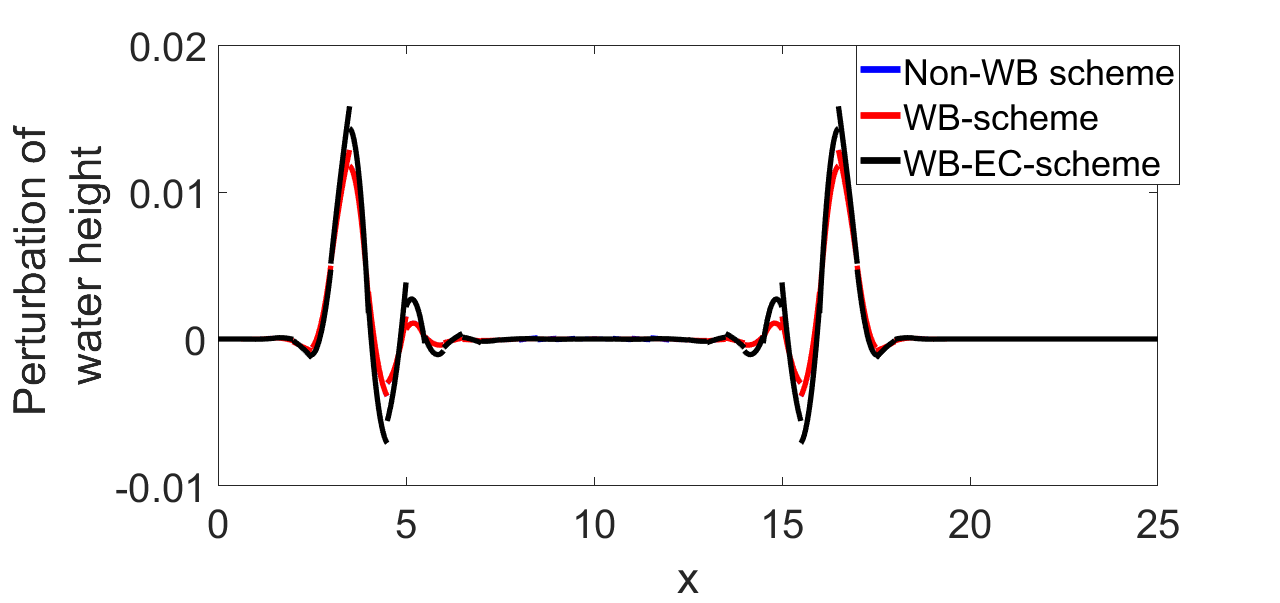}}\hspace{-0.4cm}	
	\subfigure[$\xi=10^{-3}m$]{\includegraphics[width=0.5\textwidth]{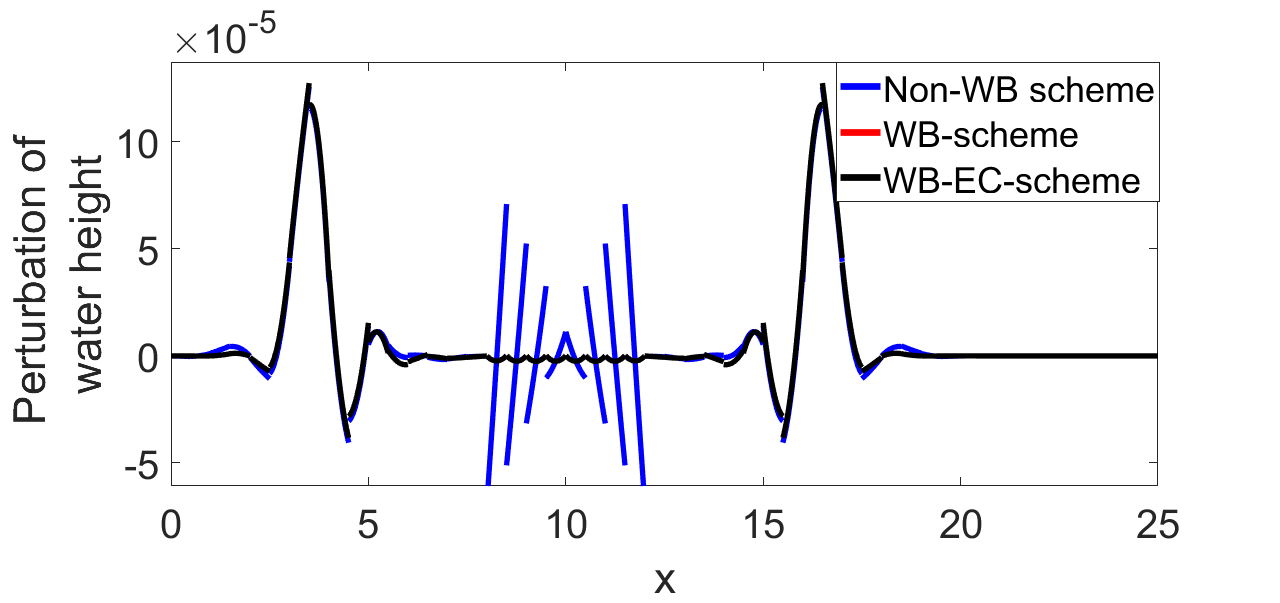}}	
	
	        \centering\centering\subfigure[$\xi=10^{-5}m$]{\includegraphics[width=0.5\textwidth]{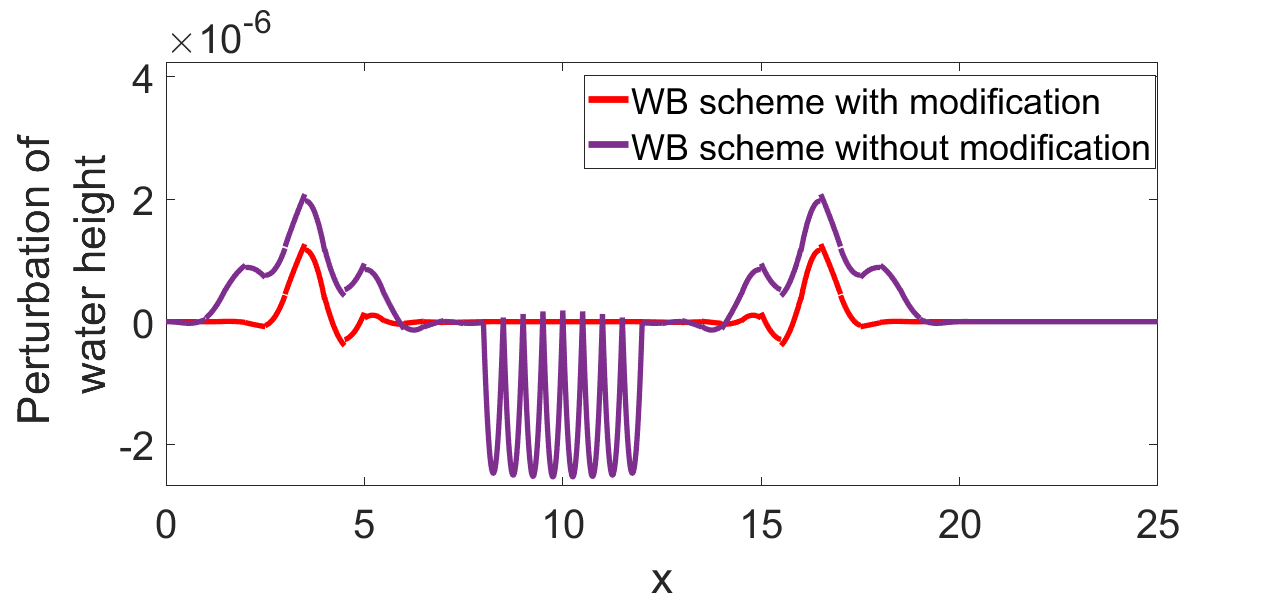}}
%	\subfigure[$\xi=10^{-6}$]{\includegraphics[width=0.33\textwidth]{oneD/Bathymetry/height_lakerest_pert6zoom}}
%>>>>>>> 8633838fbb44d758f3df865c7cef7dc43c507d06
	\caption{Perturbation of the lake at rest. Top:  NWB DGSEM (blue) and global flux quadrature approach  with (black) and without (red) cell entropy correction for $p=2$.
	 Bottom: Global flux quadrature without entropy correction using  the basic nodal source definition  \eqref{eq:S_node} (magenta) and  the modified  one  \eqref{eq:S_node-mod} (red)}
	\label{fig:lake_at_rest_pert}.
\end{figure}

For the first case we consider $u(0,x)=0$ and $h^*(x)+b(x)=2m$ where $b$ is as given in \eqref{bathymetry_1},
and the spatial domain is $x\in [0,25]m$. \revp{The boundary conditions are given by $h(t,0)=h(t,25)=2\text{m}$ and $q(t,0)=q(t,25)=0\text{m}^2/\text{s}$.}
We consider three values for the amplitude of  the perturbation  $\xi =10^{-1}m$, $ \xi=10^{-3}m$, and $\xi=10^{-5}m$. 
We set $x_0=10$ in \eqref{perturbation}, and plot the solution   at $T=1.5$. \revp{The tests in this section are performed on a mesh with  $N=50$ cells.}
We compare the  non-well-balanced scheme,  and the well-balanced scheme  using \eqref{eq:S_node-mod} 
with and without entropy correction. Note that for this case the choice of the analytical entropy flux  or of the modified one to define the cell correction has no effect.
This aspect will be studied later in the results section.
The results  obtained  here with $p=2$ is as shown in figure \ref{fig:lake_at_rest_pert}.

From the results in figure \ref{fig:lake_at_rest_pert}, one can see that  with $p=2$ the largest amplitude perturbation is well resolved by all methods. 
However, as we go lower with the amplitude of the perturbation we see the spurious oscillations of the non-well-balanced scheme.
For the last case, the results of the latter cannot even be plotted on the same scale since the amplitude of the error with respect to the steady solution is several
orders of magnitude higher than the perturbation. For completeness, we also report in the same figure a comparison of the global flux quadrature method with 
the straightforward definition of the nodal source  \eqref{eq:S_node}, and with the modified one  \eqref{eq:S_node-mod}. 
To compare the two \revPO{and highlight the advantage of the modified formulation} we have to use a much   smaller amplitude perturbation of $\xi=10^{-5}$m. The results for the \revPO{larger perturbations are}  almost superimposed. 
%This shows   the advantage of having  an exact definition of preservation  when extremely small perturbations are considered on coarse meshes.
%

\subsubsection{Frictionless one dimensional moving equilibria} 
We consider again the equilibria defined by 
\eqref{steady0}  setting respectively  $(q_0,E_0)= (4.42\text{m}^2/\text{s},22.05535\text{m}^2/\text{s}^2)$,
and   $(q_0,E_0)= (4.42\text{m}^2/\text{s},28.8971\text{m}^2/\text{s}^2)$.  The bathymetry  is  given in \eqref{bathymetry_1},
and the spatial domain is $x\in [0,25]m$. \revp{The boundary conditions are given by  $h(t,0)=h(t,25)=2\text{m},q(t,0)=q(t,25)=4.42\text{m}^2/\text{s}$ for sub-critical flow and $h(t,0)=h(t,25)=0.66\text{m}, q(t,0)=q(t,25)=4.42\text{m}^2/\text{s}$ for super-critical flow.}  As in the previous case, we consider the evolution
of perturbations of different amplitudes of  the analytical steady depth for $p=2$. We compare the non well-balanced  DGSEM,
with the one using global flux quadrature with and without cell entropy correction. 

\begin{figure}[H]
%<<<<<<< HEAD
	\centering\subfigure[$\xi=10^{-1}$]{\includegraphics[width=0.5\textwidth]{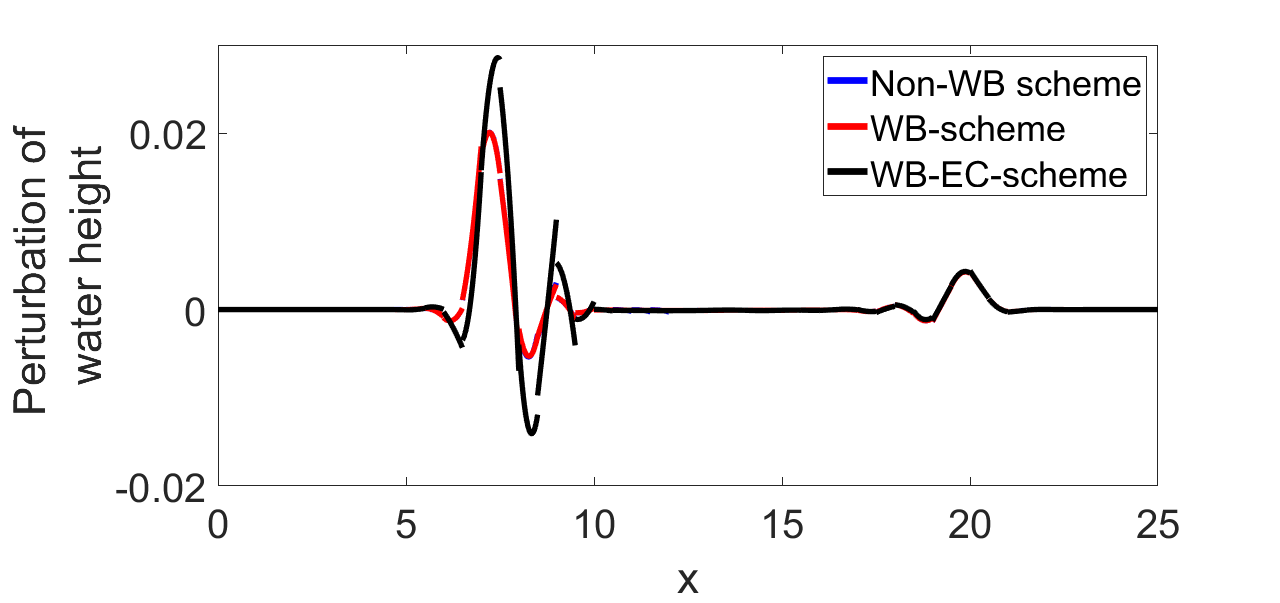}}\hspace{-0.4cm}	
	\subfigure[$\xi=10^{-3}$]{\includegraphics[width=0.5\textwidth]{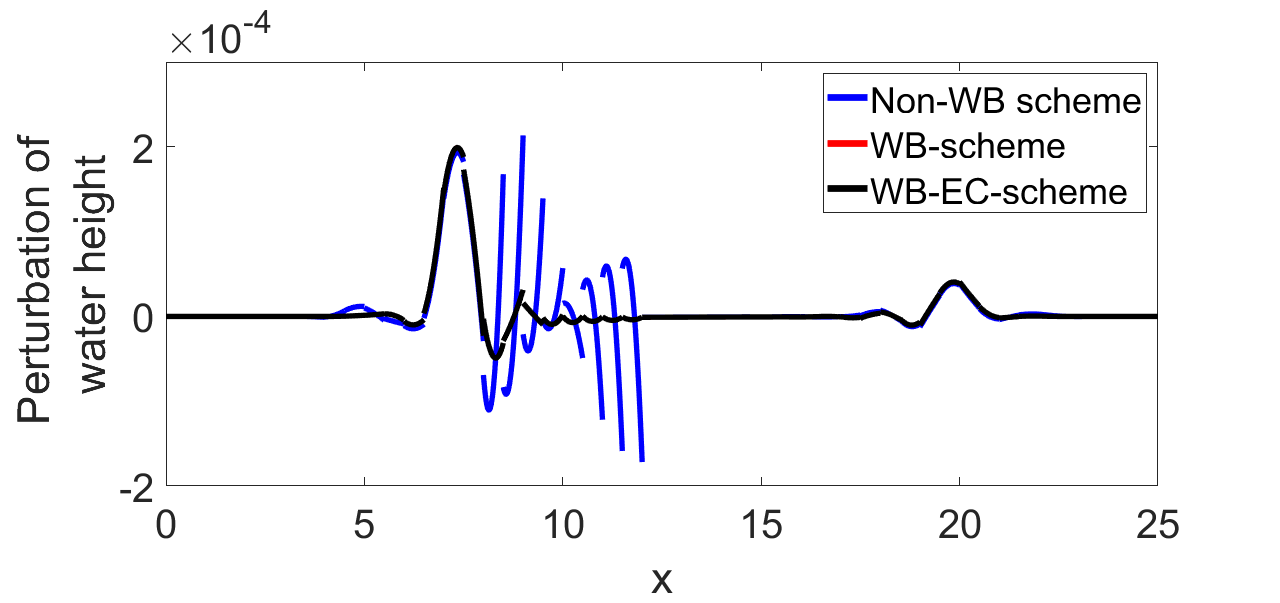}}	
%	\subfigure[$\eta=10^{-5}$]{\includegraphics[width=0.33\textwidth]{oneD/Bathymetry/diffheight_subcrit_gf_pert5_ECinconsistency}}
%=======
%	\hspace{-0.6cm}\subfloat[$\xi=10^{-1}$]{\includegraphics[width=0.5\textwidth]{oneD/Bathymetry/diffheight_subcrit_gf_pert1}}\hspace{-0.4cm}	
%	\subfloat[$\xi=10^{-3}$]{\includegraphics[width=0.5\textwidth]{oneD/Bathymetry/diffheight_subcrit_gf_pert3}}	
%>>>>>>> 8633838fbb44d758f3df865c7cef7dc43c507d06
		\caption{Perturbation of frictionless 1d sub-critical equilibrium. Perturbation plot  for the  NWB DGSEM (blue), and for the DGSEM with global flux quadrature   with (black) and without (red) cell entropy correction for $p=2$. }
	\label{fig:subcritical_pert}
\end{figure}

\begin{figure}[H]
%<<<<<<< HEAD
\centering\subfigure[$\xi=2*10^{-2}$]{\includegraphics[width=0.5\textwidth]{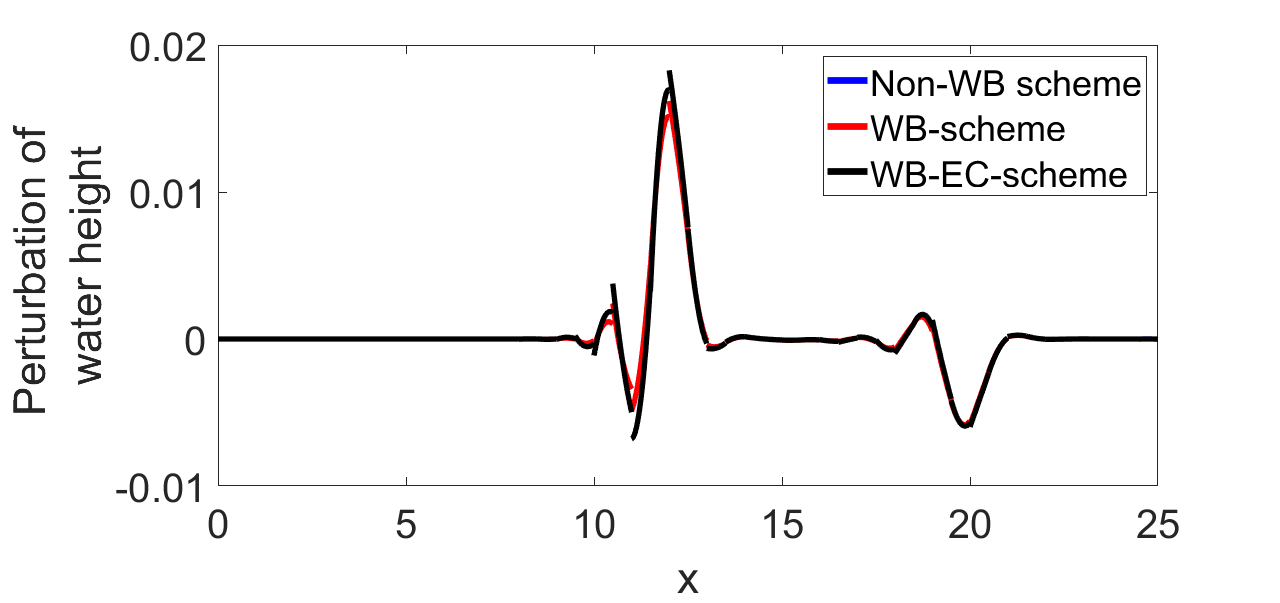}}\hspace{-0.4cm}	
	\subfigure[$\xi=10^{-5}$]{\includegraphics[width=0.5\textwidth]{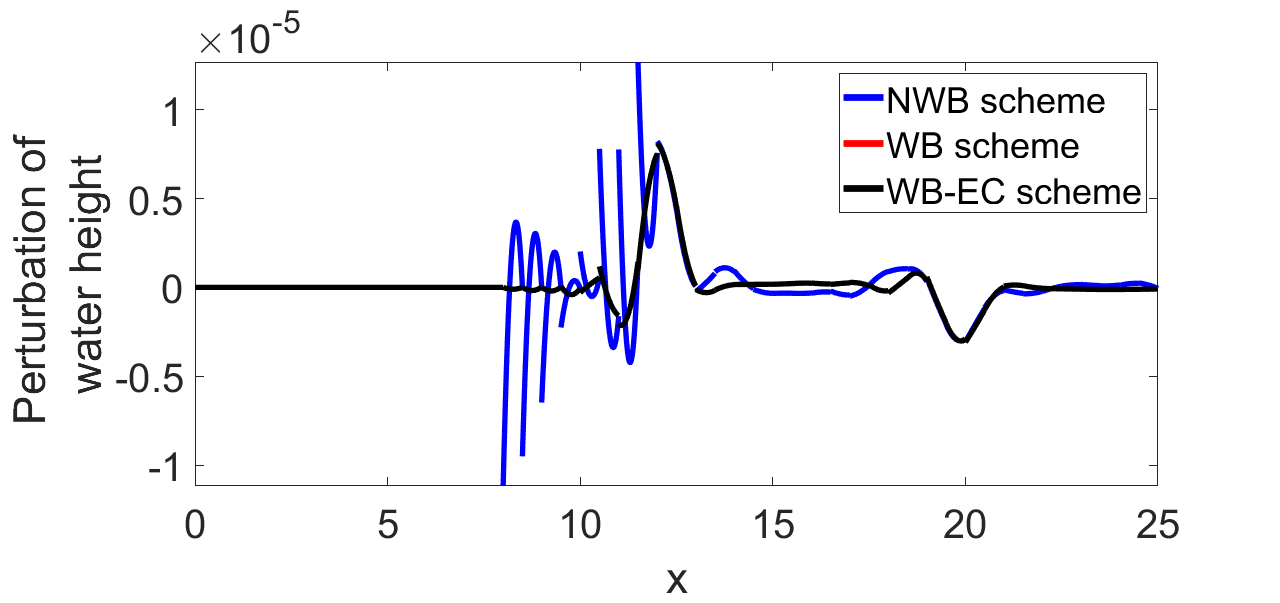}}	
%	\subfigure[$\eta=10^{-6}$]{\includegraphics[width=0.33\textwidth]{oneD/Bathymetry/diffheight_supercrit_gf_pert6_ECinconsistency}}
		\caption{Perturbation of frictionless 1d super-critical equilibrium. Perturbation plot  for the  NWB DGSEM (blue), and for the DGSEM with global flux quadrature   with (black) and without (red) cell entropy correction for $p=2$. }
%=======
%	\hspace{-0.6cm}\subfloat[$\xi=2*10^{-2}$]{\includegraphics[width=0.5\textwidth]{oneD/Bathymetry/diffheight_supercrit_gf_pert0_02}}\hspace{-0.4cm}	
%	\subfloat[$\xi=10^{-5}$]{\includegraphics[width=0.5\textwidth]{oneD/Bathymetry/diffheight_supercrit_gf_pert5}}	
%	\caption{$h-h^*$ with NWB, WB, WB-EC discontinuous Galerkin schemes for perturbation to super-critical equilibrium}
%>>>>>>> 8633838fbb44d758f3df865c7cef7dc43c507d06
	\label{fig:supercritical_pert}
\end{figure}

For the  sub-critical case we  have  considered perturbations of order $\xi =10^{-1}m$, and $\xi=10^{-3}$m, initially set at $x_0=10$. 
The results at   $T=1.5$s  are plotted on figure \ref{fig:subcritical_pert} in terms of $h(x,t) - h^*(x)$. We perform a similar exercise for the   supercritical flow, 
only considering $\xi=2\times 10^{-2}$m, and $\xi=10^{-5}$m. The initial perturbation is this time set at $x_0=6.25$. The results are reported  in figure \ref{fig:supercritical_pert}. 

\begin{figure}[H]
	%<<<<<<< HEAD
	\centering\subfigure[$\xi=10^{-1}$]{\includegraphics[width=0.5\textwidth]{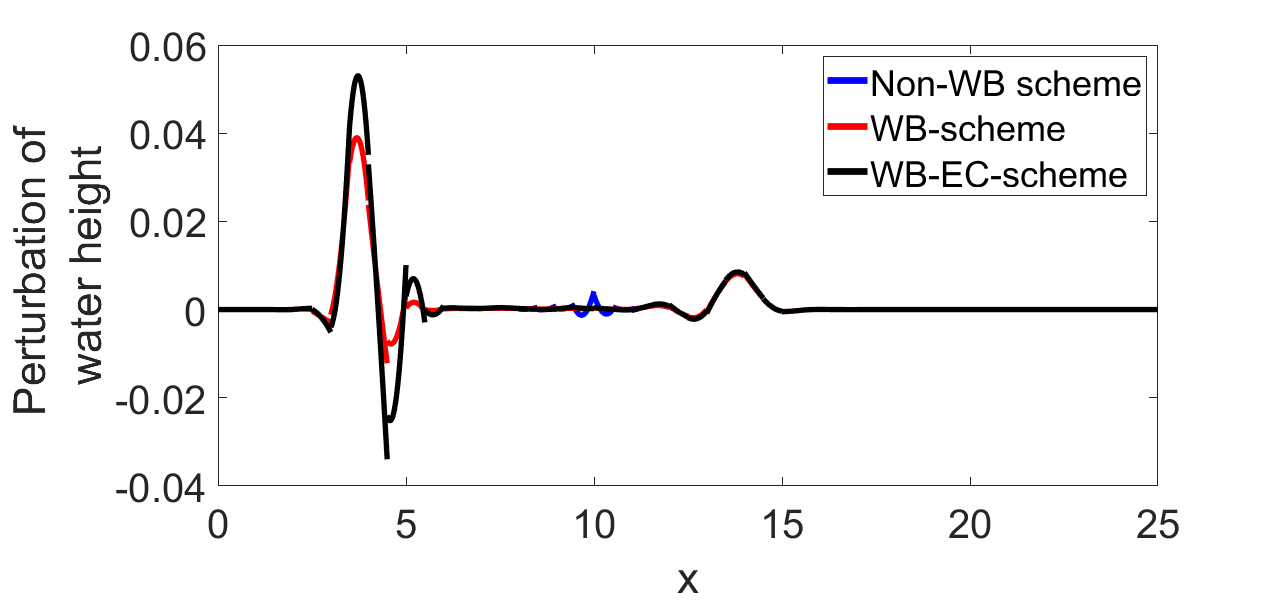}}\hspace{-0.4cm}	
	\subfigure[$\xi=10^{-3}$]{\includegraphics[width=0.5\textwidth]{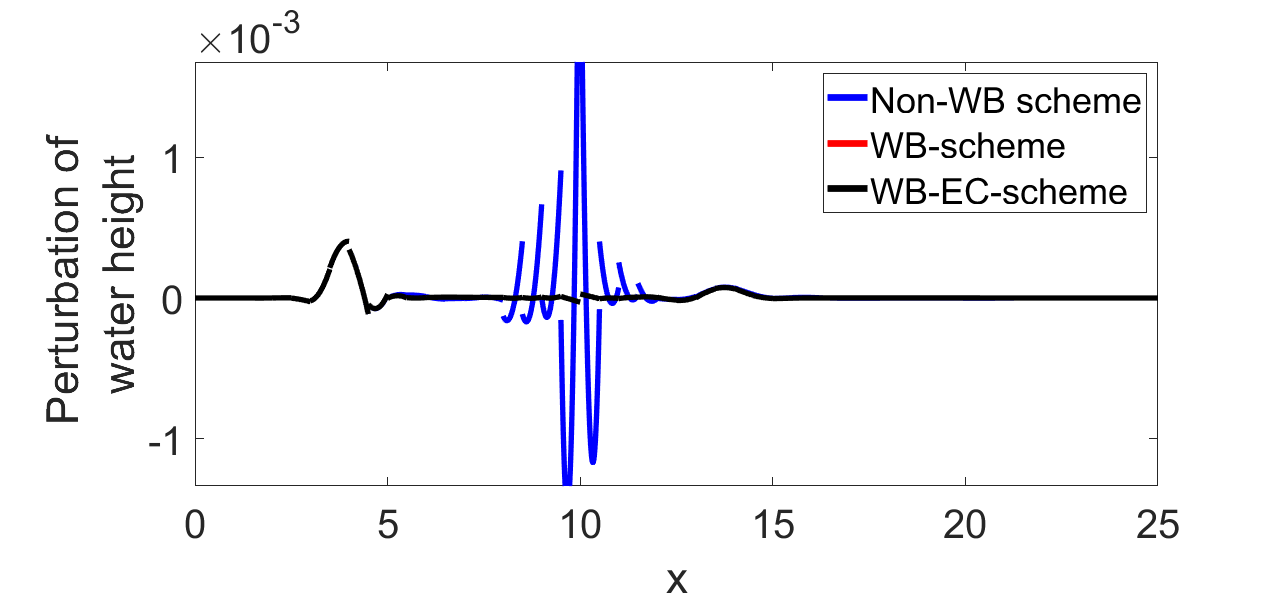}}	
	\caption{Perturbation of frictionless 1d trans-critical equilibrium. Perturbation plot  for the  NWB DGSEM (blue), and for the DGSEM with global flux quadrature   with (black) and without (red) cell entropy correction for $p=2$. }
	%=======
	\label{fig:transcritical_pert}
\end{figure}

\revt{In the same way, we also perform a test for trans-critical flow, with $(q_0,E_0)=(1.53m^2/s,11.0863m^2/s^2)$ and perturbations of of order $\xi =10^{-1}m$, and $\xi=10^{-3}$m, initially set at $x_0=6.25$. It is to be noted that in this case the flow for $x<10$ is sub-critical and is super-critical for $x>10$ with a critical point at $x=10$. The boundary conditions for $(h,hu)$ are calculated accordingly with respect to the given steady state. The results for the perturbation $h(x,t) - h^*(x)$ at $T=1.5$s  are shown in figure \ref{fig:transcritical_pert}.}

\revt{For all three cases, sub-, trans- and super-critical flow} the behaviour observed is similar to the one  seen for the lake at rest state. 
For $p=2$ relatively large perturbations are not affected by the well-balanced nature of the scheme, however as the 
amplitude of the perturbations decreases, the non-balancing error leads to spurious artefacts of amplitude larger than the physical waves.\\

\begin{figure}[H]
	\centering
	{\includegraphics[width=0.6\textwidth]{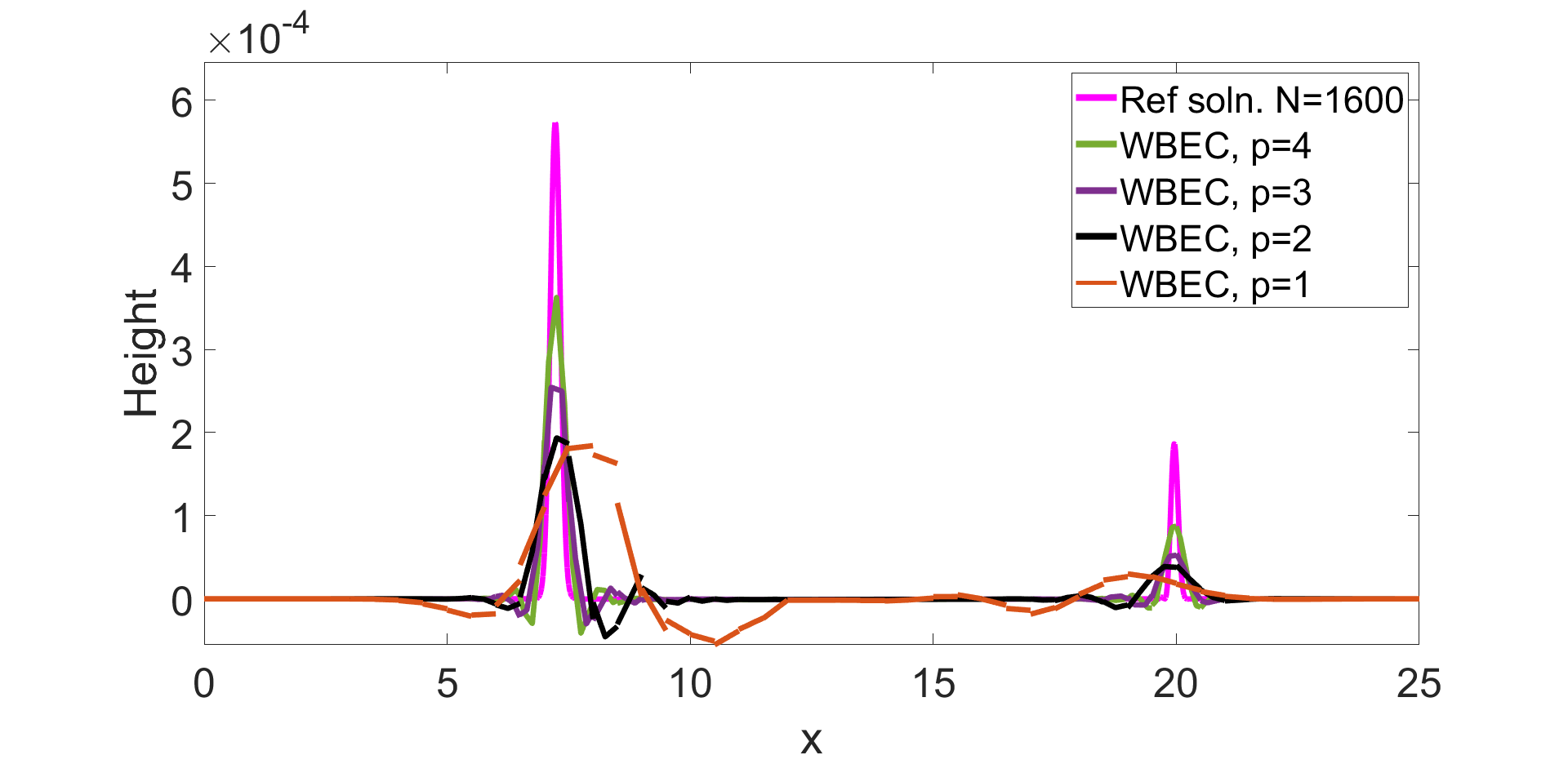}}
	\caption{Perturbation of frictionless 1d sub-critical equilibrium. $p$-convergence study of the perturbation at $T=1.5$s.}
	\label{fig:DG_order_subcritical}
\end{figure}

For completeness we also report a $p-$convergence study for the  sub-critical \revp{with $\xi=10^{-3}$m} case using the global flux quadrature scheme with entropy correction.
The result is reported  in figure \ref{fig:DG_order_subcritical}  showing a classical behaviour: passing from $p=1$ to $p=2$ almost removes the phase error;
passing to $p=3$ and $p=4$ allows to  remove both phase, and amplitude errors.

\subsubsection{Pseudo-one dimensional equilibrium at rest  with transverse Coriolis effects} 
\begin{figure}[H]
%%<<<<<<< HEAD
%	%	%\centering	
%	\includegraphics[width=0.55\textwidth]{oneD/Coriolis/height_cor_pert_4}
%	\caption{Left: $L_1$ error- mesh size with NWB, WB, WB-EC discontinuous Galerkin scheme for flow at  equilibrium given by \eqref{eq:SWE_Coriolis_equilibrium_1}. Right: 
%		$h-h^*$ for  $\eta=10^{-4}$ perturbation o to equilibrium}
%	\label{fig:Coriolis__moving_pert}
%\end{figure} 
%
%%=======
	%	\centering	
	\centering\subfigure[$\xi=0.5$]{\includegraphics[width=0.5\textwidth]{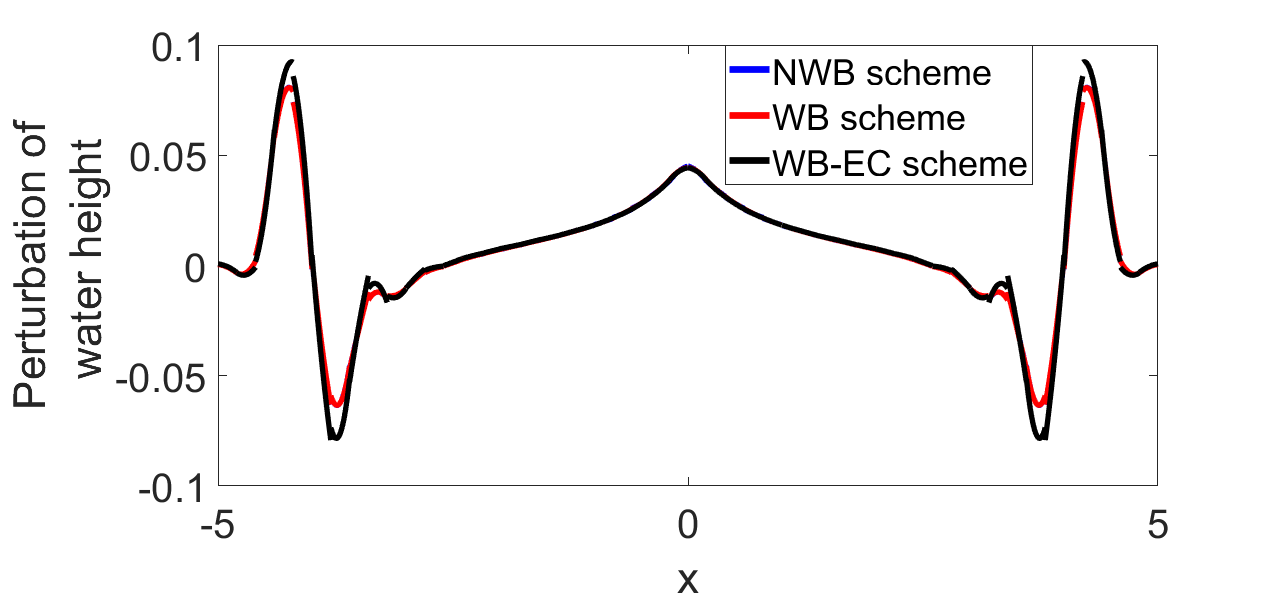}}\subfigure[$\xi=10^{-3}$]{\includegraphics[width=0.5\textwidth]{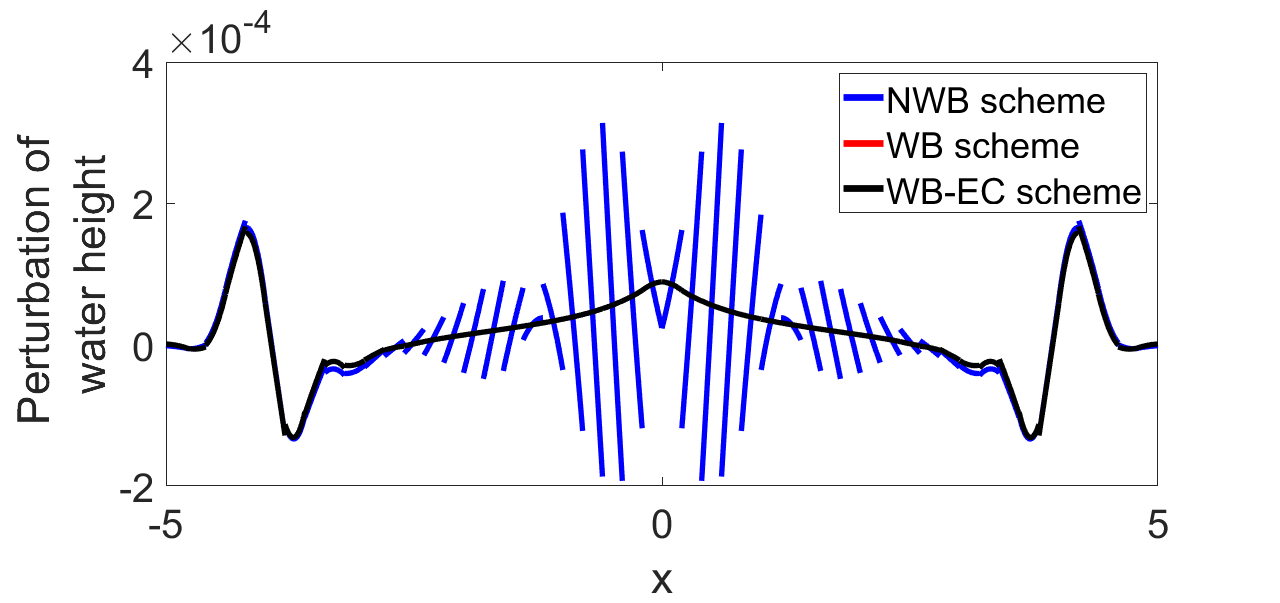}}
	\caption{Pseudo-1d equilibrium at rest with transverse Coriolis effects. Perturbation computed with the $p=2$ DGSEM using the non well-balanced (blue), and the global flux quadrature forms with (black) and without (red) entropy correction.}
	\label{fig:Coriolis_stationary_pert}
\end{figure}

We consider here the first configuration of the two steady states discussed in section \S2.1,
given analytically by the relations \eqref{steady4}. In particular, on the domain  $[-5,5]$ we consider the transverse velocity field
\begin{align}
	v(x)=\frac{g x}{2}e^{-x^2}.
	\label{eq:SWE_Coriolis_equilibrium_1}
\end{align}
and set $g=1\text{m}/\text{s}^2$, and $\omega=2\text{s}^{-1}$. The values for $h^*(x)$ are calculated analytically from \eqref{steady4} \revp{with $\zeta_0=2\text{m}$. The steady state is then  used to determine the boundary conditions for $U$ at $x=-5,5$.}
We consider perturbations  with $\xi=0.5$m, and $\xi=10^{-3}$m, initially centered at $x_0=0$. The   perturbation  $h-h^*(x)$ at $T=2$s
is studied. We compare the global flux quadrature DGSEM with   $p=2$ with the corresponding non well-balanced formulation.
Results are plotted  in figure \ref{fig:Coriolis_stationary_pert}. As for the previous cases
the largest perturbation is evolved in a similar manner by all the schemes. However, only the well-balanced formulations 
correctly capture small amplitude variations as shown by the right plot in the figure.

%>>>>>>> 8633838fbb44d758f3df865c7cef7dc43c507d06

\subsubsection{Frictionless pseudo-one dimensional equilibria with Coriolis effects}
\begin{figure}[H]
%<<<<<<< HEAD
%	\hspace{-0.6cm}\subfigure[$\eta=10^{-1}$]{\includegraphics[width=0.33\textwidth]oneD/Coriolis/height_cor_pert_0_5}}\hspace{-0.4cm}	
%	\subfigure[$\eta=10^{-3}$]{\includegraphics[width=0.33\textwidth]{ooneD/Coriolis/height_cor_pert_4}}	
%	\subfigure[$\eta=10^{-6}$]{\includegraphics[width=0.33\textwidth]{oneD/Bathymetry/height_lakerest_pert6zoom}}
%	\caption{$h-h^*$ with NWB, WB, WB-EC discontinuous Galerkin schemes for perturbation to lake at rest equilibrium}
%	\label{fig:friction_subcrit_pert}
%\end{figure}
%=======
	%	\centering	
	\subfigure[$\xi=0.5$]{\includegraphics[width=0.5\textwidth]{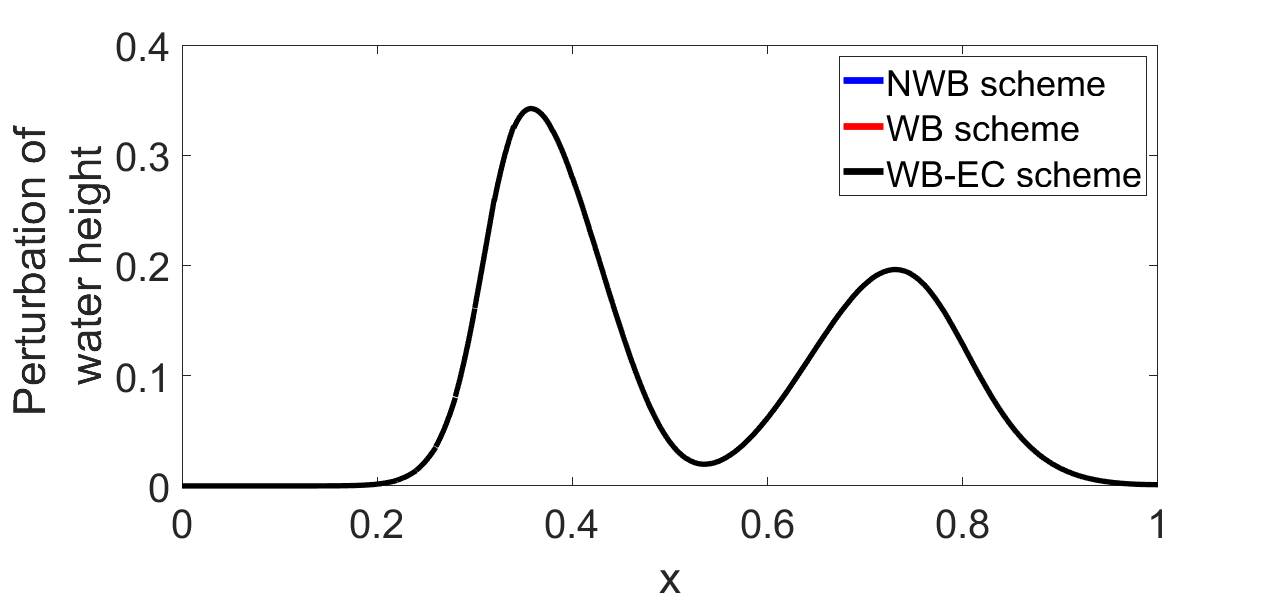}}	
	\subfigure[$\xi=10^{-4}$]{\includegraphics[width=0.5\textwidth]{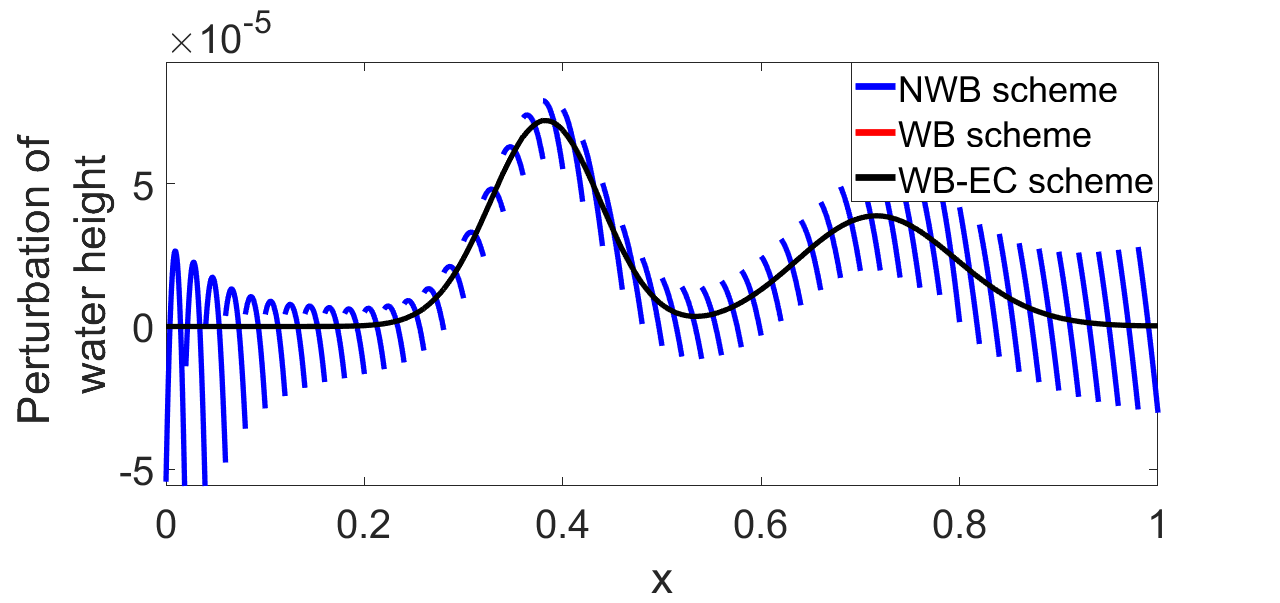}}
	\caption{Pseudo-1d equilibrium  with Coriolis effects. Perturbation computed with the $p=2$ DGSEM using the non well-balanced (blue), and the global flux quadrature forms with (black) and without (red) entropy correction.}
	\label{fig:Coriolis_moving_pert}
\end{figure}

We  repeat the above test for the steady solution given by \eqref{eq:SWE_Coriolis_equilibrium_mov_an}. We add to the analytical equilibrium
studied in section \S7.1 perturbations  of amplitudes   $\xi=0.5$m, and $\xi=10^{-4}$m at $x_0=0.5$.
We compute the evolution with $p=2$ until $T=0.1$s with all the different DGSEM formulations with  $p=2$.
In particular,   figure \ref{fig:Coriolis_moving_pert} compares the perturbation $h-h^*(x)$ obtained with the global quadrature formulations
and with the non-well-balanced scheme.   
We can see just as the other examples that the non-well-balanced scheme fails to capture a small perturbation, which is not the case with the 
global flux quadrature approach.

%>>>>>>> 8633838fbb44d758f3df865c7cef7dc43c507d06

\subsubsection{Moving equilibria with friction} 
\begin{figure}[H]
%<<<<<<< HEAD
%	\hspace{-0.6cm}\subfigure[$\eta=10^{-1}$]{\includegraphics[width=0.33\textwidth]{oneD/friction/diffheight_fric_pert1}}\hspace{-0.4cm}	
%	\subfigure[$\eta=10^{-5}$]{\includegraphics[width=0.33\textwidth]{oneD/friction/diffheight_fric_supercrit_pert5}}	
%	%	\subfigure[$\eta=10^{-6}$]{\includegraphics[width=0.33\textwidth]{oneD/Bathymetry/height_lakerest_pert6zoom}}
%	\caption{$h-h^*$ with NWB, WB, WB-EC discontinuous Galerkin schemes for perturbation to lake at rest equilibrium}
%	\label{fig:friction_supercrit_pert}
%=======
	%	\centering
	\subfigure[Subcritical, $\xi=10^{-1}$]{\includegraphics[width=0.5\textwidth]{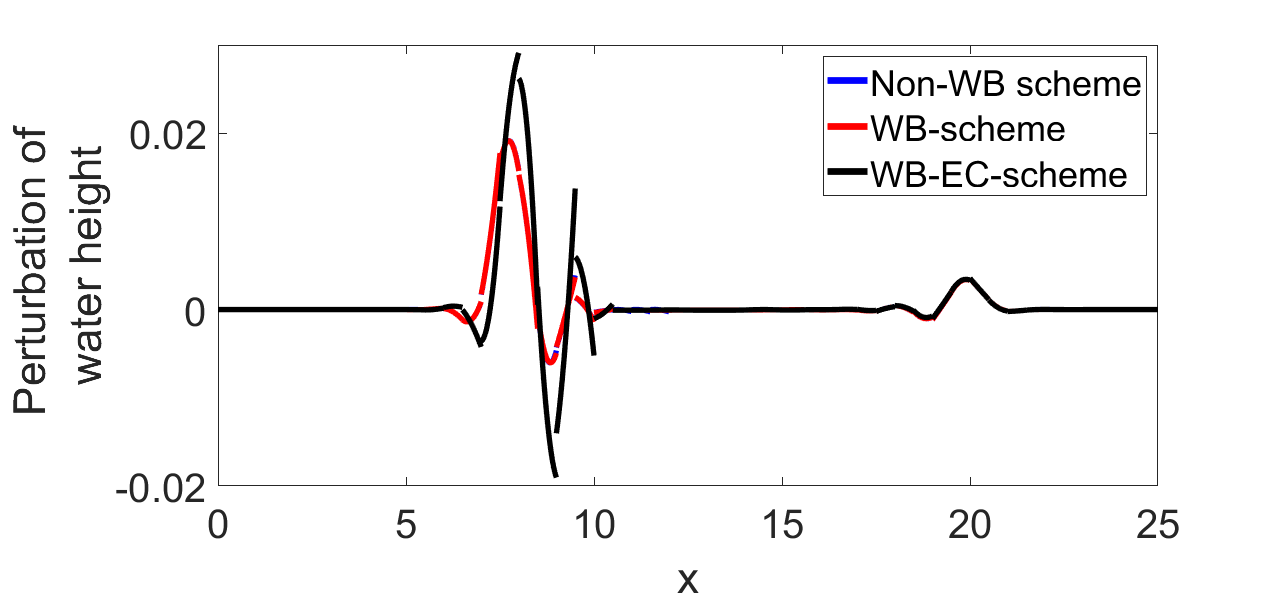}}
	\subfigure[Subcritical, $\xi=10^{-3}$]{\includegraphics[width=0.5\textwidth]{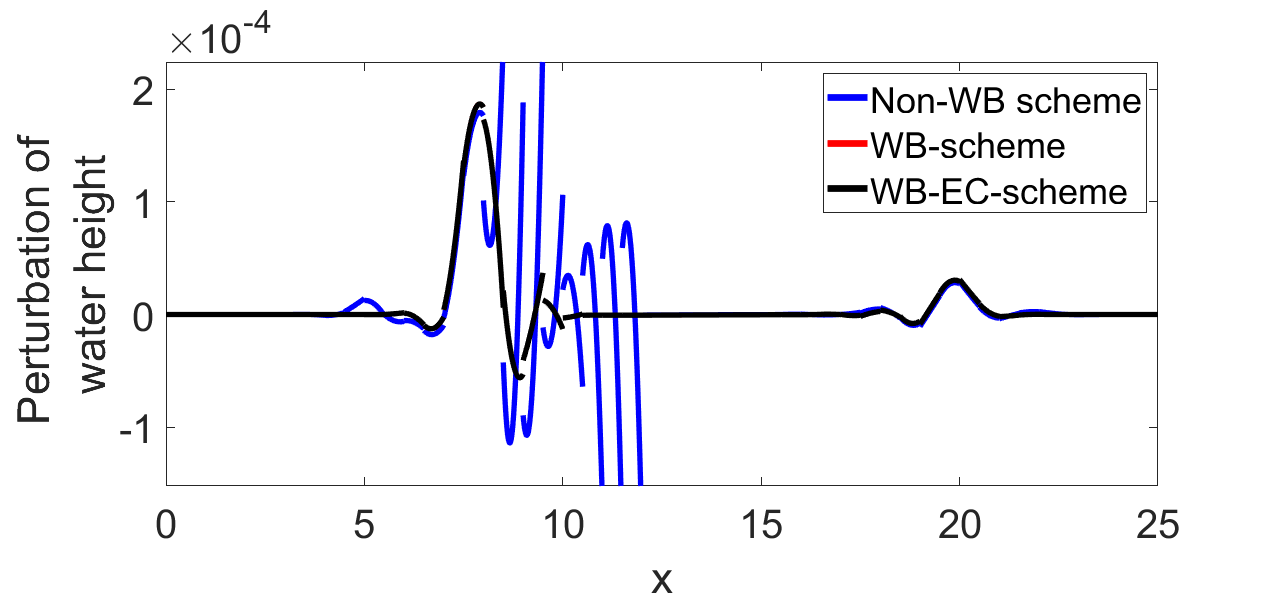}}
	
	\subfigure[Supercritical, $\xi=2\times 10^{-2}$]{\includegraphics[width=0.5\textwidth]{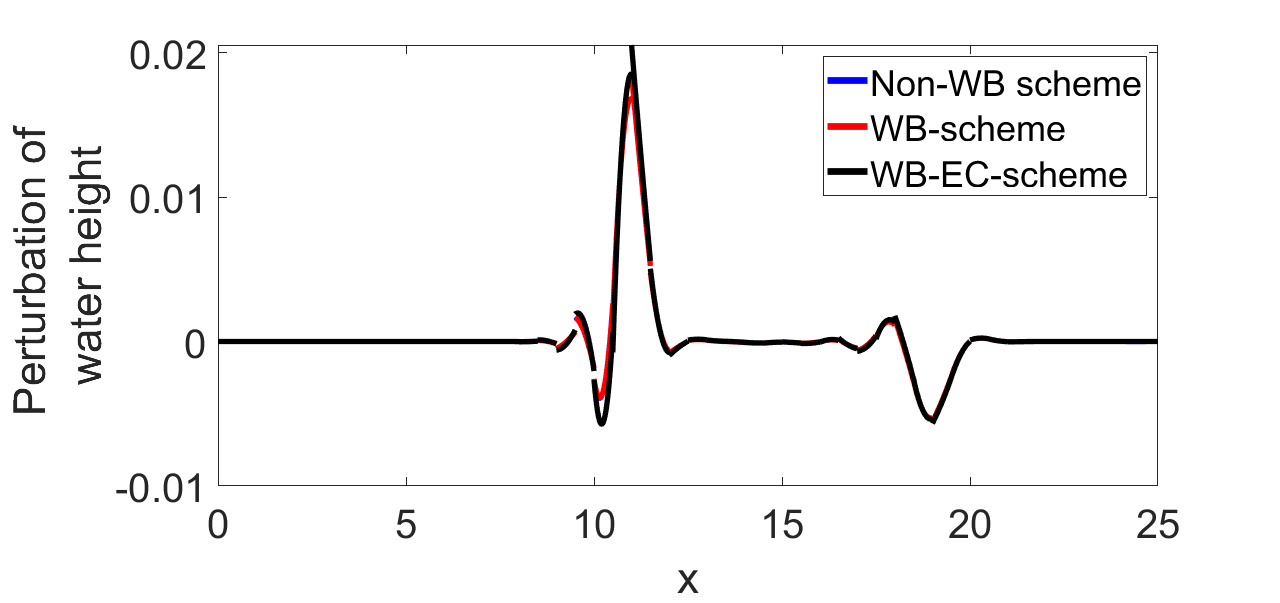}}
	\subfigure[Supercritical, $\xi=10^{-5}$]{\includegraphics[width=0.5\textwidth]{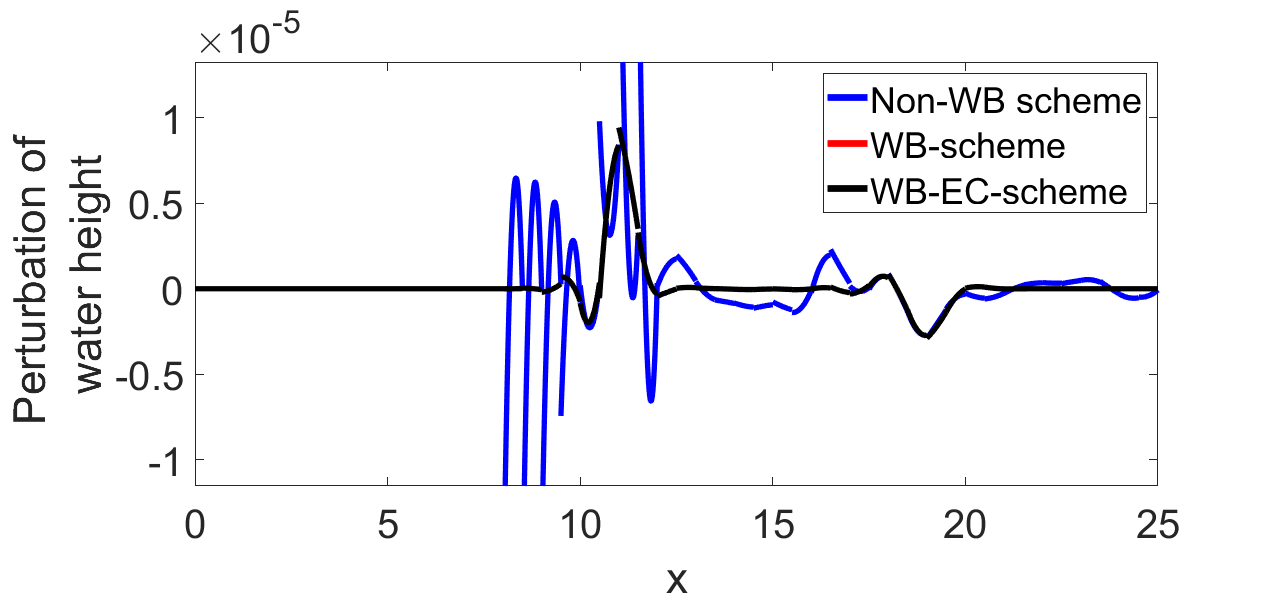}}
	
	\caption{Perturbation of moving equilibria with friction. Top: sub-critical case. Bottom: super-critical case. Comparison of the DGSEM with $p=2$  global quadrature formulations  
	  with (black) and without (red) cell entropy correction, with the non well-balanced one (blue). }
	\label{fig:friction_perturbation}
%>>>>>>> 8633838fbb44d758f3df865c7cef7dc43c507d06
\end{figure}

As a last case we consider a variation of the sub- and super-critical flows given by \eqref{steady0} and already investigated,
only with the addition of friction effects. In this case, the analytical solution is characterized by \eqref{steady6}, and a reference solution can be computed e.g. evaluating the
integrals on a refined mesh. We consider the same sub-critical and super-critical states as before on the domain $[0,25]$m, and bathymetry \eqref{bathymetry_1}. \revp{The boundary conditions are calculated from the steady state solution of $h,hu$ at $x=0$ and $x=25$.}
For simplicity we use constant values of the friction coefficient, set to  $c_f=0.03$ in the sub-critical case, and $c_f=0.05$ in the super-critical case.
As before we add   perturbations of different amplitudes.  To best visualize the difference between all the schemes, we have studied this time
$\xi=10^{-1}$m and $\xi=10^{-3}$m in the sub-critical case. For  the super-critical equilibrium we have set
 $\xi=2\times 10^{-2}$m, and $\xi=10^{-5}$m. We compute the  solution with $p=2$ until $T=1.5$s. We compare
 the solutions obtained with the different schemes   in figure \ref{fig:friction_perturbation}.
 The results  lead to similar conclusions as for all the other equilibria, and source term types.

\subsection{Entropy control vs well-balancing}
We have so far  considered  entropy correction terms  always  consistent with the initialization used. 
In this section we    check the impact of this consistency, as well as  the effect of the adding the cell correction on the fully discrete evolution of the entropy.\\

\begin{figure}[H]
	%	\centering
	\subfigure[$\xi=10^{-4}$]{\includegraphics[width=0.5\textwidth]{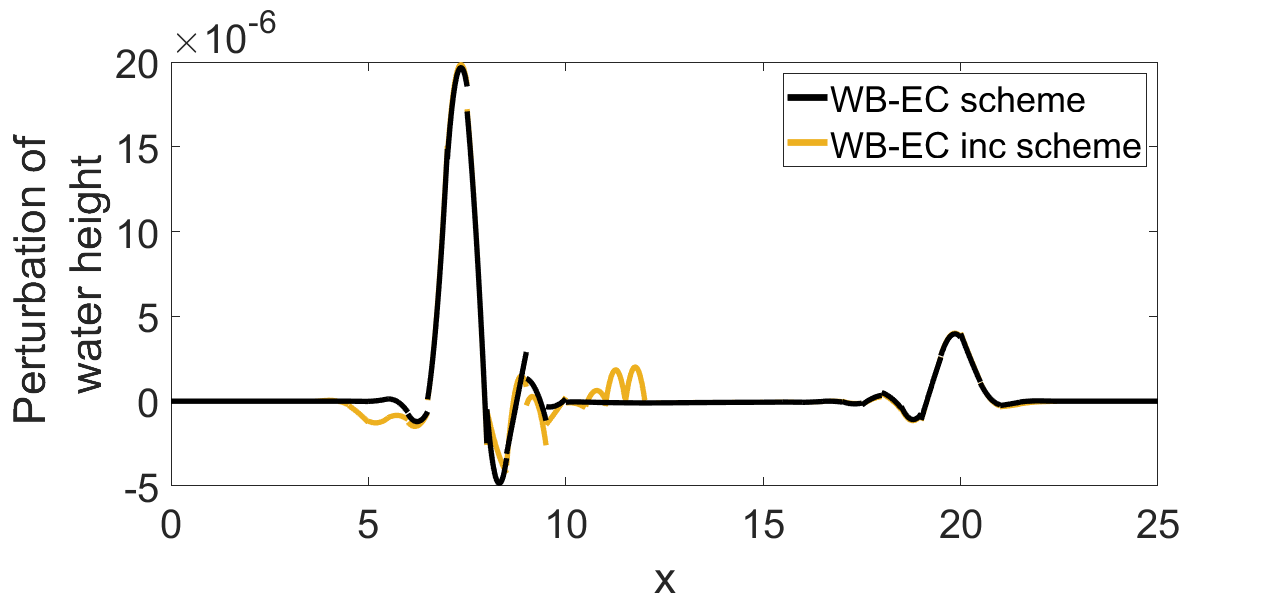}}	
	\subfigure[$\xi=10^{-5}$]{\includegraphics[width=0.5\textwidth]{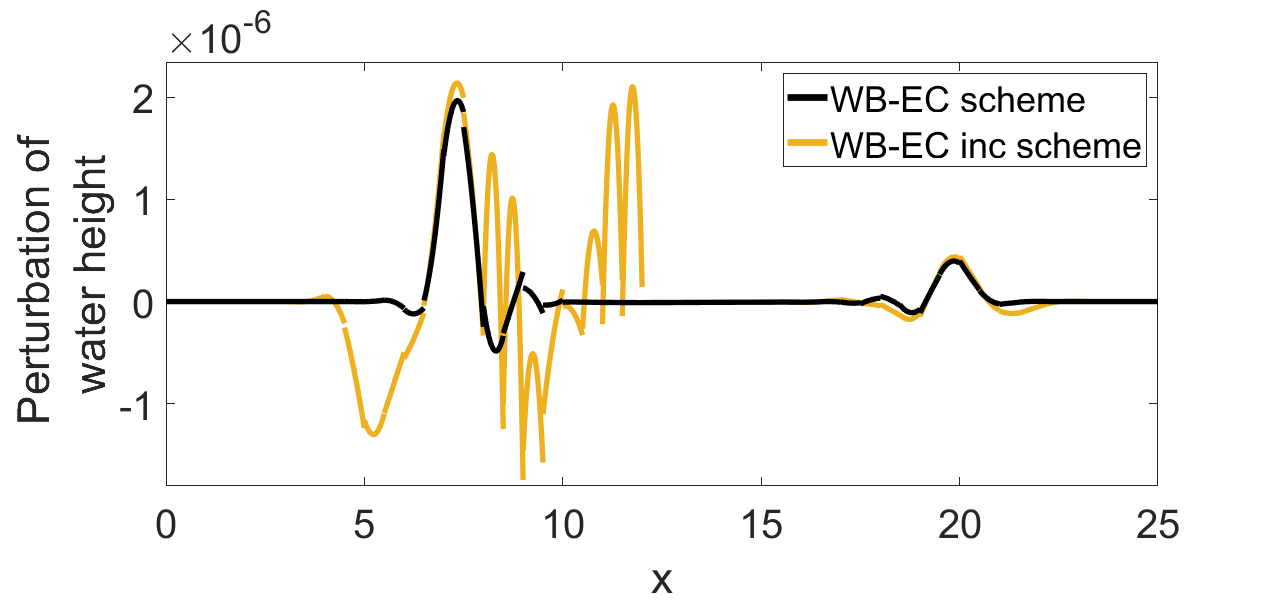}}
	
	\subfigure[$\xi=10^{-5}$]{\includegraphics[width=0.5\textwidth]{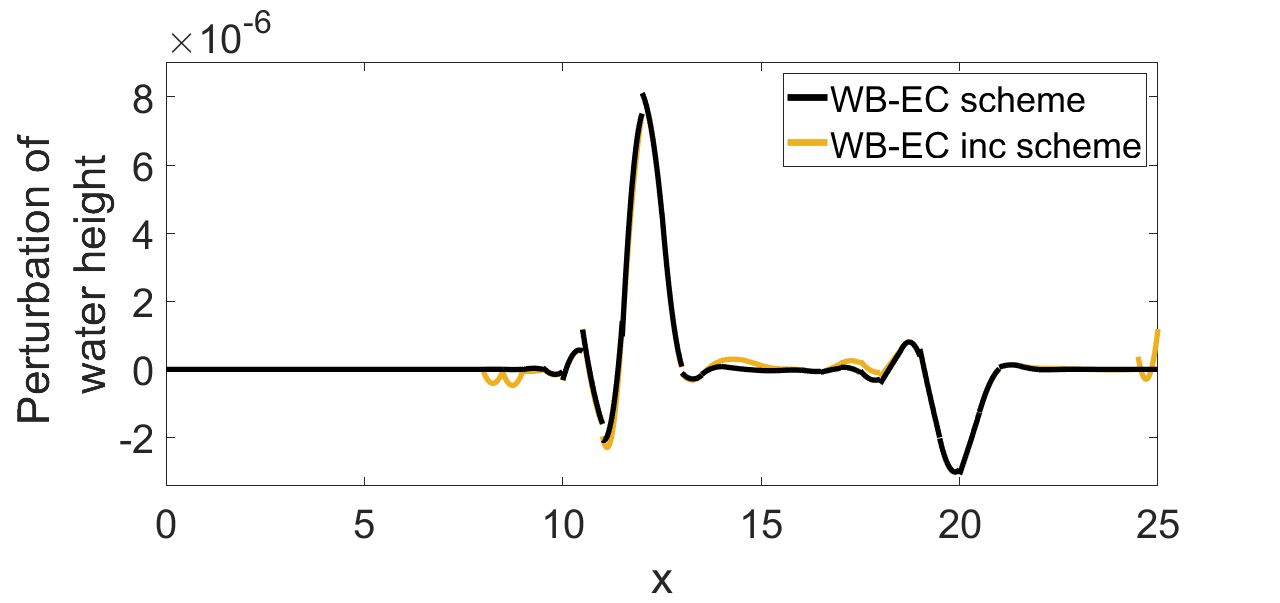}}	
	\subfigure[$\xi=10^{-6}$]{\includegraphics[width=0.5\textwidth]{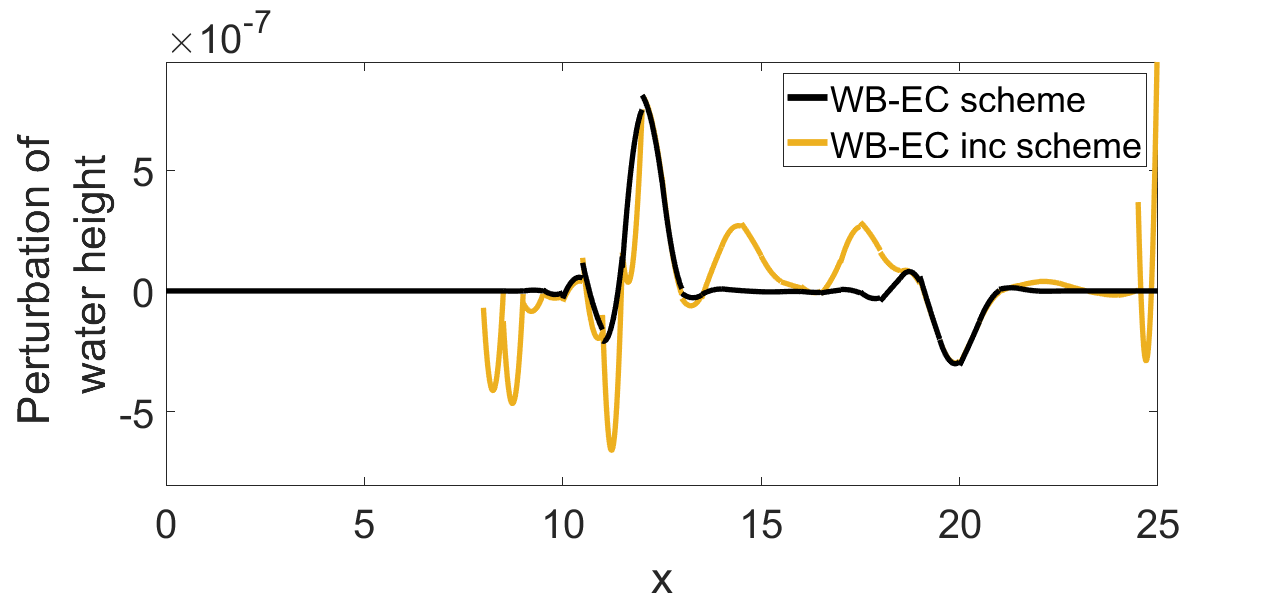}}
	\caption{Frictionless 1D moving equilibria. Perturbations of the discrete global flux solution in the sub-critical (top) and super-critical (bottom) case. 
	   Comparison between using  consistent (black) and inconsistent (yellow) entropy corrections.}
	\label{fig:energy_inconsistency}
\end{figure}

To begin with, we consider a slight variation of the study of the previous section, and
we perturb the initial data given by the global flux solutions for the moving  one dimensional sub-critical and supercritical flows without friction.
We consider small perturbations with amplitudes $10^{-4}$m and $ 10^{-5}$m for the sub-critical case, and $10^{-5}$m and $10^{-6}$m for the super-critical one. The results,
in figure \ref{fig:energy_inconsistency}, show that for small amplitude perturbations properly accounting for the initialization strategy is important.
Despite of the fact that the superconvergence of the  global flux solution to the analytical one, 
the use of analytical entropy fluxes in this case spoils the well-balanced character of the scheme, leading to spurious waves comparable to 
the physical ones. Note that this effect is a much weaker effect than the one observed with the non well-balanced method.
However, in  some cases   it may be necessary to have  machine accuracy preservation  of a physically relevant initial state, and initializing with
the global flux solution may be relevant. In these cases, the use of \eqref{eta-corr15} in the entropy correction is necessary.

\begin{figure}[H]
%<<<<<<< HEAD
%	\includegraphics[width=0.55\textwidth]{oneD/Bathymetry/ent_subcrit_pert1}\hspace{-0.4cm}
%	% 
%	%  \includegraphics[width=0.5\textwidth]{oneD/Bathymetry/ent_subcrit_pert3}\\
%	%		\subfigure[$\eta=10^{-1}$]{
%	\includegraphics[width=0.5\textwidth]{oneD/Bathymetry/ent_supercrit_pert1}%}	
%	%	\subfigure[$\eta=10^{-3}$]{\includegraphics[width=0.5\textwidth]{oneD/Bathymetry/ent_supercrit_pert5}}
%	
%	\caption{Energy with NWB, WB, WB-EC discontinuous Galerkin scheme for $eta=10^{-1}$ perturbation to sub-critical  (left) and super-critical equilibrium (right)}
%	\label{fig:eq_flow_subsonic_supersonic_energy}
%=======
	\includegraphics[width=0.5\textwidth]{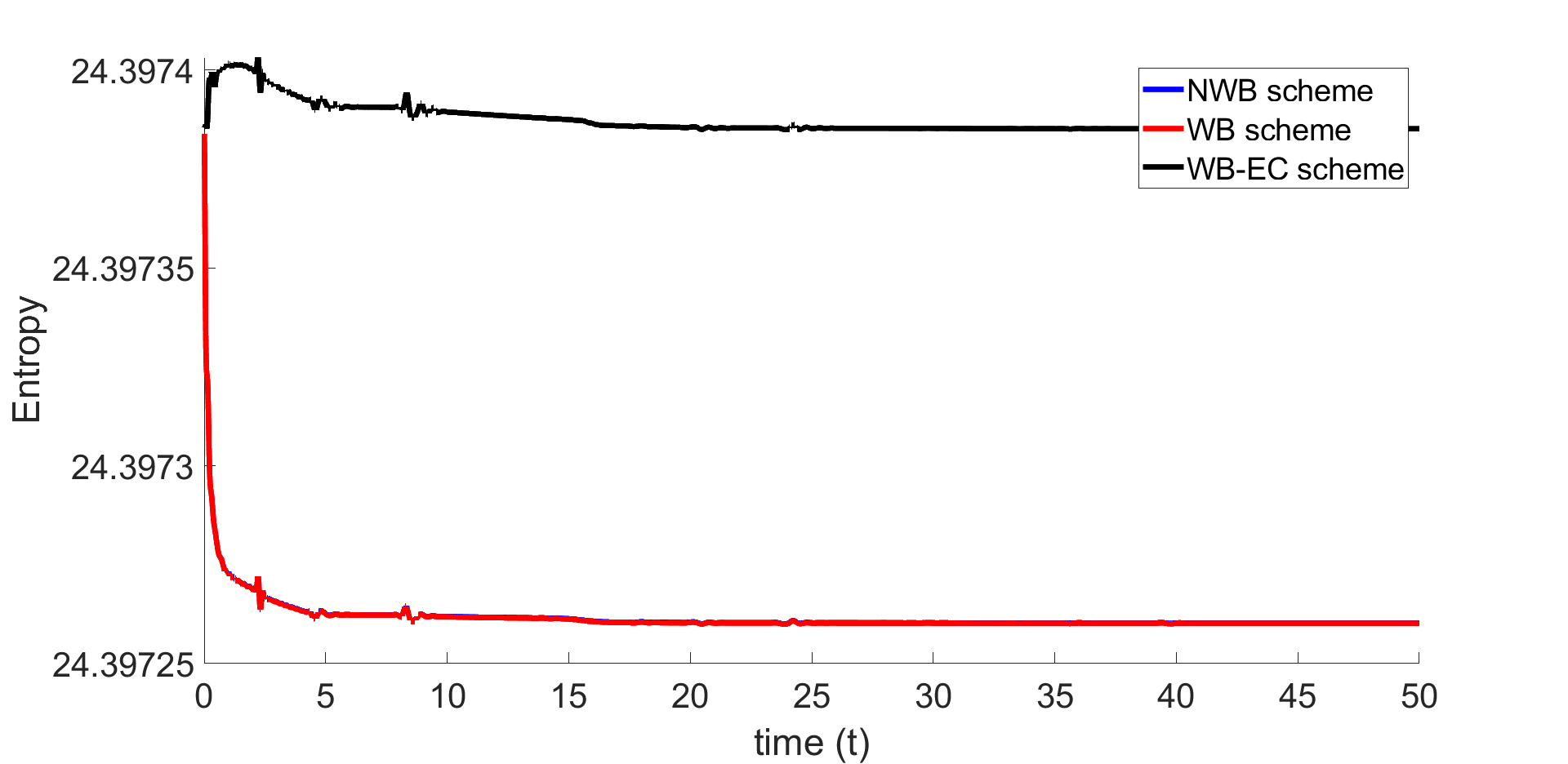}\hspace{-0.2cm}\includegraphics[width=0.5\textwidth]{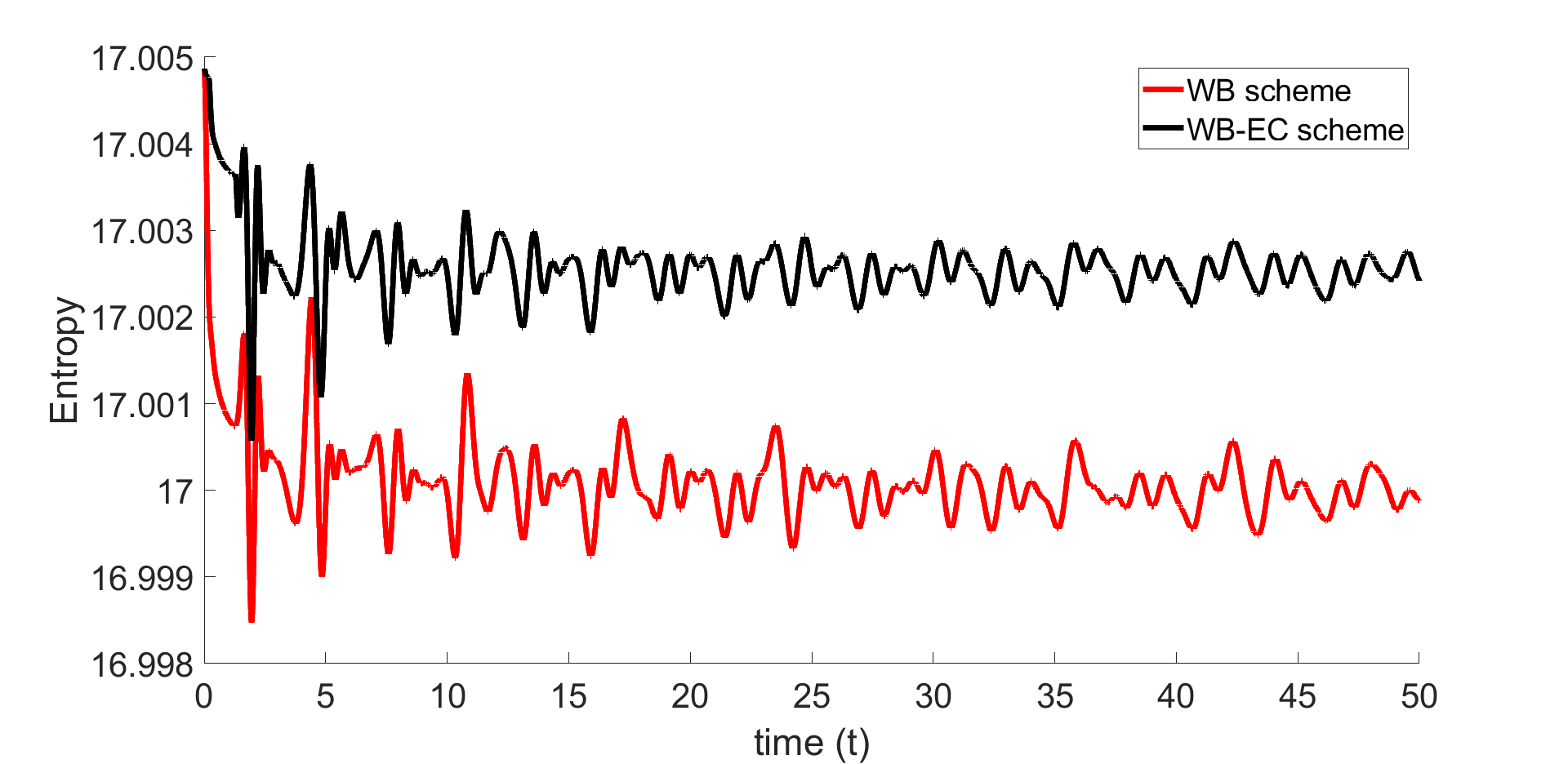}
	\caption{Time discrete evolution of the total energy in the domain for a perturbation of amplitude $\xi=10^{-1}$m to  sub-critical  (left) and super-critical
	  frictionless moving 1d equilibrium (right). Global quadrature DGSEM with entropy correction (black) compared to the basic global quadrature DGSEM (red) 
	    and to the non well-balanced DGSEM (blue).}
	\label{fig:energy_longtime}
%>>>>>>> 8633838fbb44d758f3df865c7cef7dc43c507d06
\end{figure}

We now look at the impact of the entropy correction term on the fully discrete time evolution of the total entropy/energy 
$\sum_K \int_K \eta_{\hh}$. Note that nothing special is done concerning the time integration, so the conservation property holds exactly only in the time continuous case. 
Fully discrete variants can be obtained by other means, but are not considered here (see e.g.\cite{ranocha2020relaxationRKEuler,abgrall2022relaxation,gaburro2023high}). 
Time integration is in particular performed  with  standard RK schemes of the same formal order of the underlying spatial discretization.
We will comment on two examples. The first is the 
%We start  by looking at   a 
perturbation of amplitude  $\xi=10^{-1}$m of the sub-critical and supercritical moving equilibria without friction.
The evolution of entropy until time $T=50$s   in the $p=2$  case for the NWB, and for   the global flux quadrature with and without cell entropy correction 
are  shown in figure \ref{fig:energy_longtime}.  % \ref{fig:moving_eq_energy}.
We can see that  in both cases the entropy correction allows to reduce the variation of entropy in time, improving its conservation. 
%
%\begin{figure}[H]
%	\includegraphics[width=0.5\textwidth]{oneD/Bathymetry/ent_subcrit_pert1}\hspace{-0.2cm}
%	\includegraphics[width=0.5\textwidth]{oneD/Bathymetry/ent_supercrit_pert1}	
%	\caption{Time discrete evolution of the total energy in the domain for a perturbation of amplitude $\xi=10^{-1}$m to  sub-critical  (left) and super-critical
%	  frictionless moving 1d equilibrium (right). Global quadrature DGSEM with entropy correction (black) compared to the basic global quadrature DGSEM (red) 
%	    and to the non well-balanced DGSEM (blue).}
%	\label{fig:moving_eq_energy}
%\end{figure}
%
%We see from figure \ref{fig:moving_eq_energy}, that the entropy correction term is able to reduce the numerical diffusion in the scheme and preserve the initial energy f the flow for both sub-critical and supercritical flow.
%We also check the energy with the three schemes over long time  $T=50$ in figure \ref{fig:energy_longtime}.  

%%Next we check the evolution of entropy for perturbation $\xi=0.5$ for both stationary and moving water equilibria with transverse flow with Coriolis term with initial conditions as given in  \eqref{eq:SWE_Coriolis_equilibrium_mov_gf} and \eqref{eq:SWE_Coriolis_equilibrium_mov_gf} respectively. The results are as given in figure \ref{fig:coriolis_pert_entropy}.
%%
%%\begin{figure}[H]
%%%<<<<<<< HEAD
%%%	%	\centering
%%%	\subfigure[Supercritical, $\eta=10^{-1}$]{\includegraphics[width=0.55\textwidth]{oneD/Coriolis/energy_pert_cor_0_5}}
%%%=======
%%		\centering
%%		\subfigure[Stationary equilibria, $\xi=0.5$]{\includegraphics[width=0.5\textwidth]{oneD/Coriolis/energy_test1_pert0_5_rot}}\subfigure[Moving equilibria, $\xi=0.5$]{\includegraphics[width=0.5\textwidth]{oneD/Coriolis/energy_pert_cor_0_5}}
%%%>>>>>>> 8633838fbb44d758f3df865c7cef7dc43c507d06
%%	
%%	\caption{Entropy evolution with NWB, WB, WBEC over time for shallow water equation with Coriolis}
%%	\label{fig:coriolis_pert_entropy}
%%\end{figure}
%%
%%We again see for both the examples, that there is an improvement in energy preservation with the entropy correction term while there is more numerical diffusion with non-well-balanced scheme and well-balanced scheme without entropy correction.

\begin{figure}[H]
	%	\centering
%<<<<<<< HEAD
%	 \subfigure[Subcritical, $\eta=10^{-1}$]{\includegraphics[width=0.55\textwidth]{oneD/friction/energy_fric_pert1}}
%\subfigure[Supercritical, $\eta=10^{-1}$]{\includegraphics[width=0.55\textwidth]{oneD/friction/energy_fric_supercrit_pert1}}
%=======
	 \subfigure[Sub-critical, $\xi=10^{-1}$]{\includegraphics[width=0.55\textwidth]{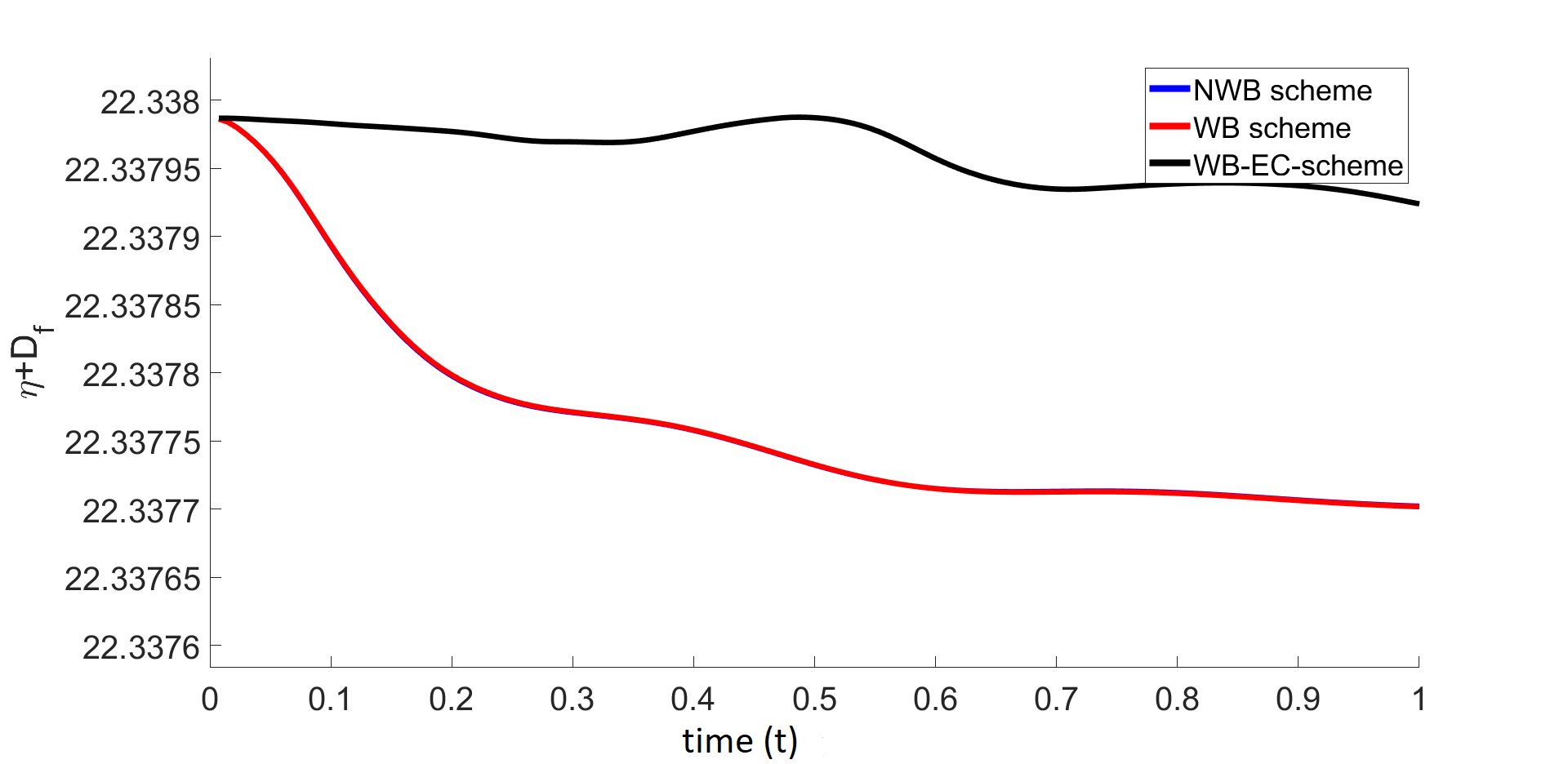}}
\subfigure[Super-critical, $\xi=10^{-1}$]{\includegraphics[width=0.55\textwidth]{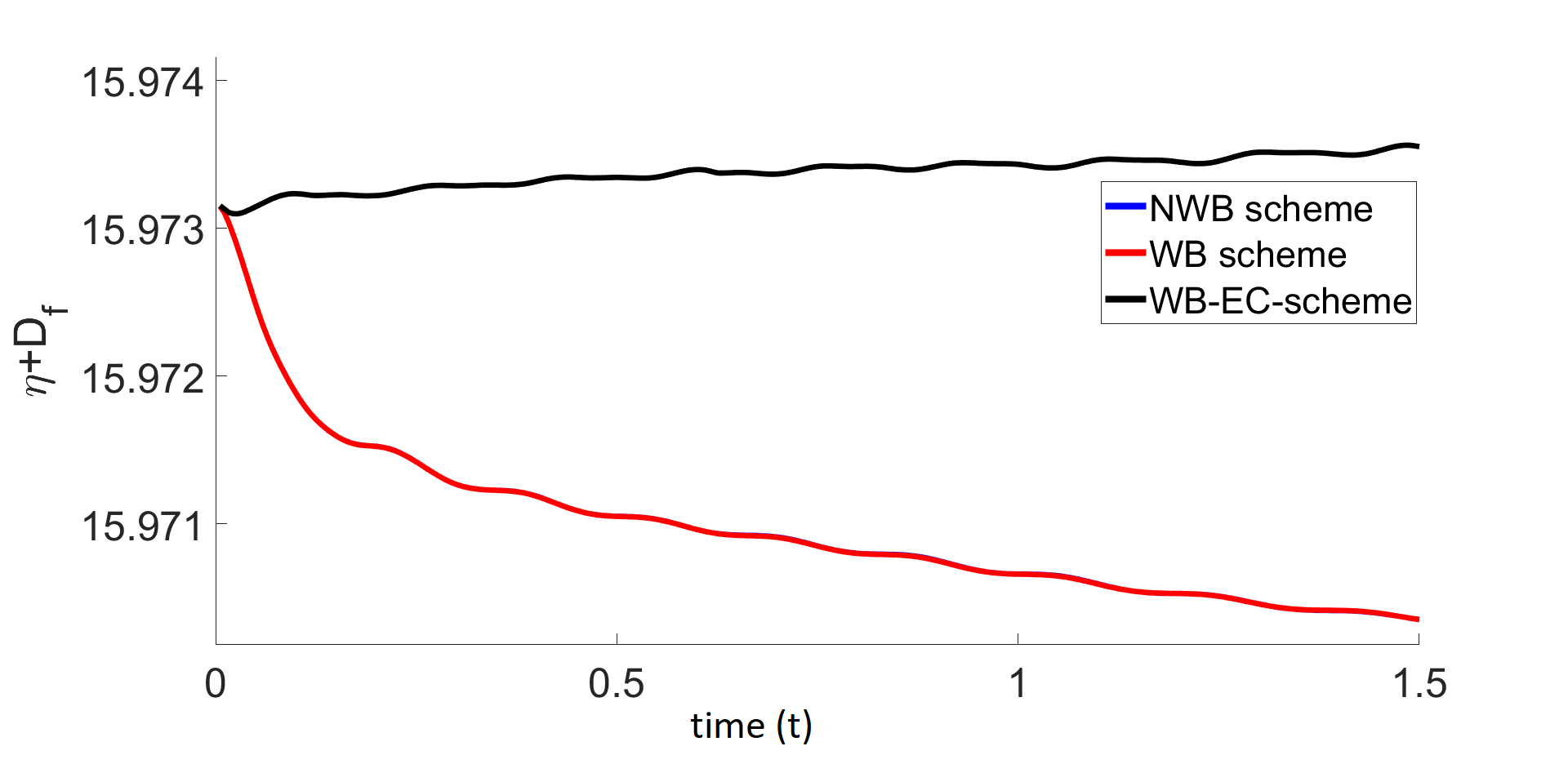}}
%>>>>>>> 8633838fbb44d758f3df865c7cef7dc43c507d06

\caption{Perturbation of  moving equilibria  with friction: evolution of the  entropy balance. Left: sub-critical flow. Right: supercritical flow.}
\label{fig:friction_pert_entropy}
\end{figure}

As a second example, we consider a perturbation $\xi=0.1$ of   moving equilibria with friction in  sub-critical and supercritical conditions.
The entropy evolution is   shown in figure \ref{fig:friction_pert_entropy}. In the figures, following \revp{Remark \ref{remark_3}}, we   study
 the evolution   of the corrected total entropy balance 
  $$N_{\text{tot}}:=\sum_K\left(|K|\bar  \eta_{\varphi} +\int_{t_0}^t\mathcal{D}_f\right)$$ 
  which should constant and equal to $\sum_K|K|\bar\eta_{\varphi}(t=0)$. We see from figure \ref{fig:friction_pert_entropy} that with the entropy correction term allows to obtain
good control of the dissipation of the scheme.

\subsection{2d examples and applications}

We consider now  the application of the scheme proposed to more complex tests. On the one hand we want to investigate
the advantage of using the tensor product extension proposed, which is only  well-balanced along mesh lines,
on genuinely 2D problems. On the other, we want to test the method on problems involving sharp propagating fronts,
and complex structures to assess its robustness.

\subsubsection{2d perturbations of 1d equilibria}

\begin{figure}[H]
\centering
	\subfigure[Steady state and bathymetry]{\includegraphics[width=0.5\textwidth]{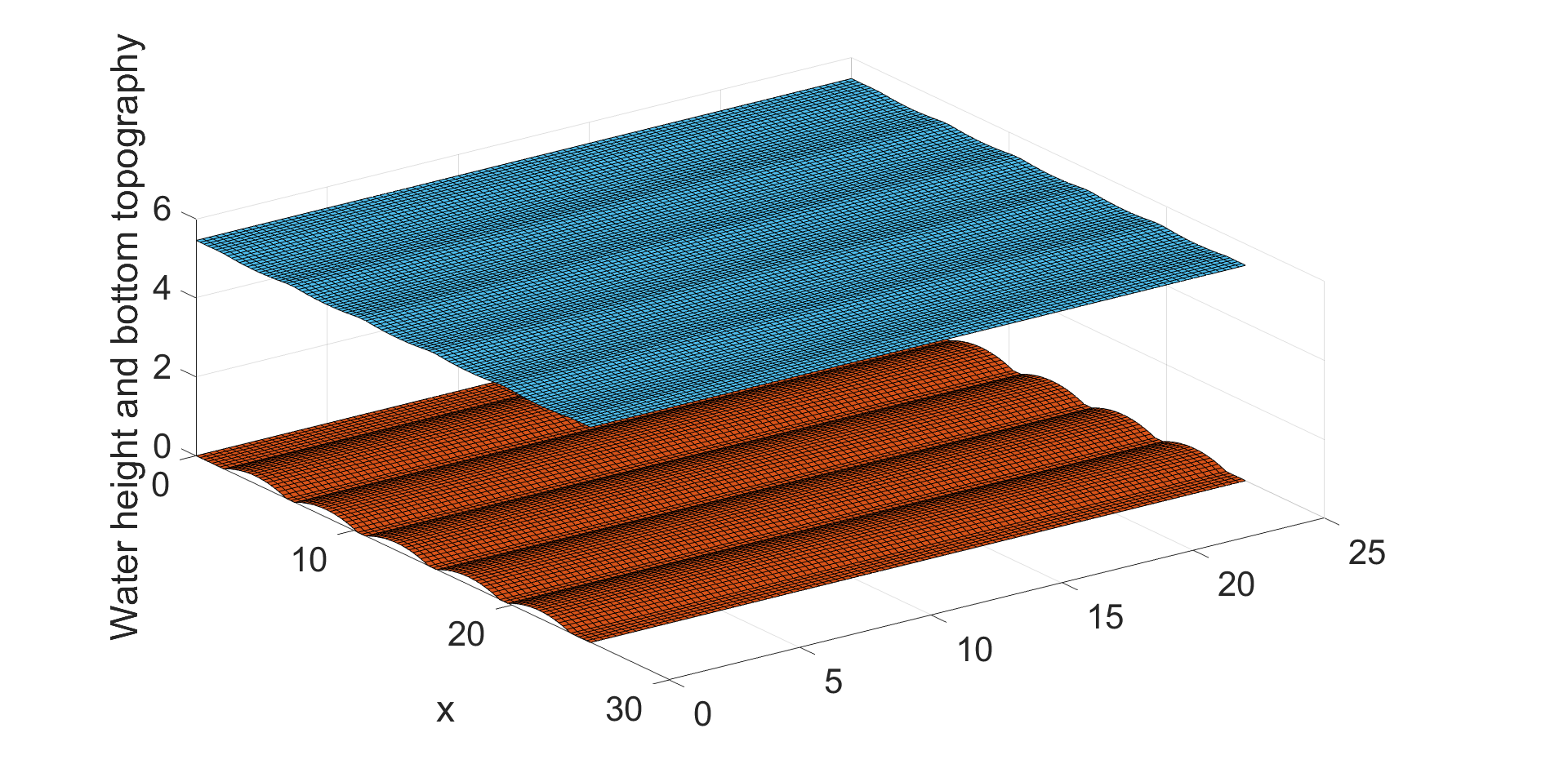}}
	
	\subfigure[NWB]{\includegraphics[width=0.5\textwidth]{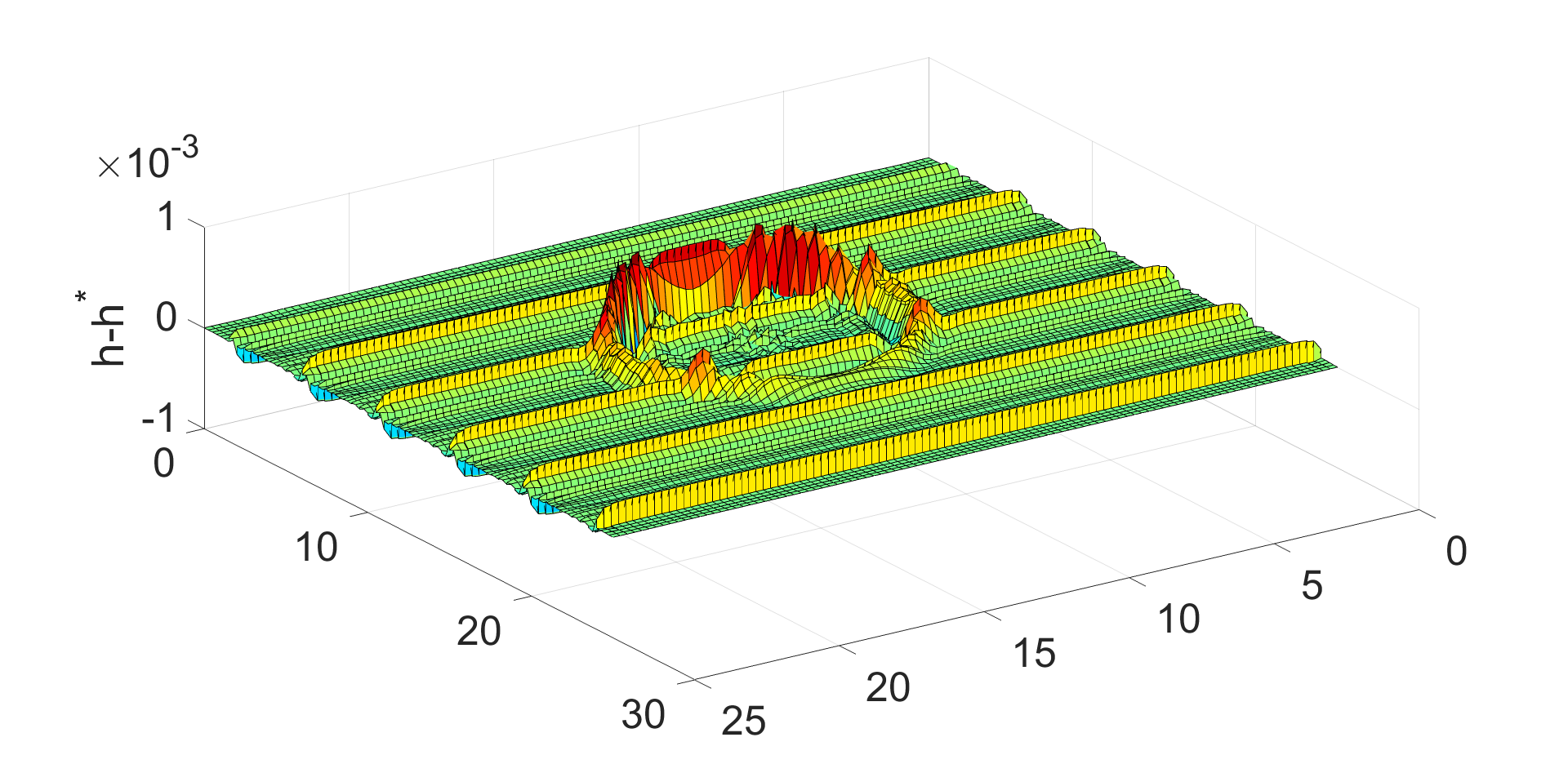}}\subfigure[WB]{\includegraphics[width=0.5\textwidth]{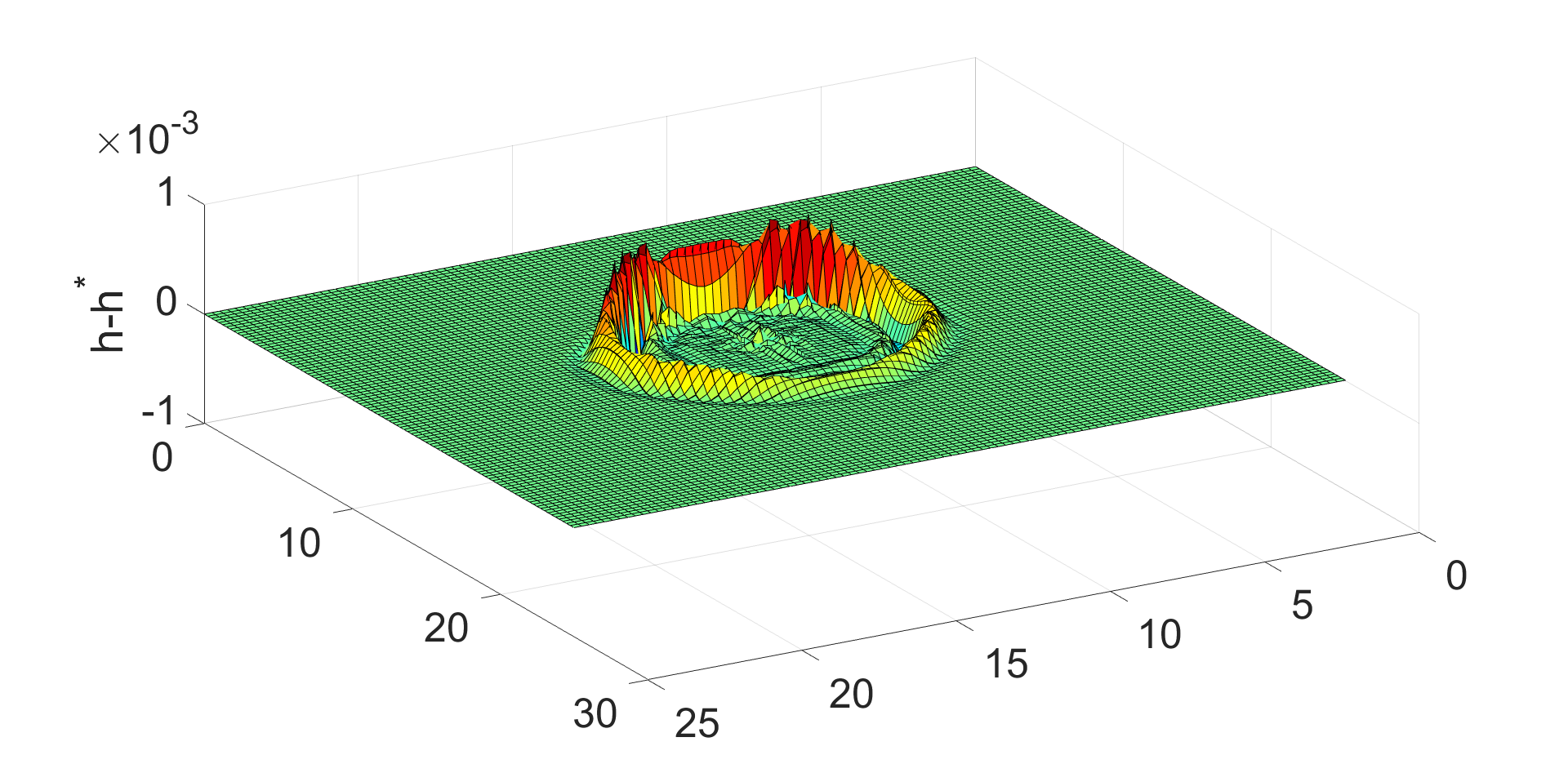}}	
%	\subfigure[$WB-EC$]{\includegraphics[width=0.5\textwidth]{twoD/oneD_equilibrium/diffheight_2dpert1deq_wbec}}
	
%<<<<<<< HEAD
%%	\subfigure[$Entropy$]{\includegraphics[width=0.75\textwidth]{twoD/oneD_equilibrium/entropy_1dequi_2d}}
%	\caption{Numerical solution for $h-h^*$ with NWB, WB, WB-EC discontinuous Galerkin scheme}
%=======
%	\subfloat[$Entropy$]{\includegraphics[width=0.75\textwidth]{twoD/oneD_equilibrium/entropy_1dequi_2d}}
	\caption{Numerical solution for $h-h^*$ with NWB, WB discontinuous Galerkin scheme}
%>>>>>>> 8633838fbb44d758f3df865c7cef7dc43c507d06
	\label{Fig:SWE_1d_eq_pert}
\end{figure}

We start from a simple benchmark: the evolution of 2D perturbations to one dimensional exact  steady solutions. 
We consider first the lake at rest with initial \revp{and boundary} condition  with $\zeta_0=5.47$m, and $hu=hv=0$. The   bathymetry
is given by a series of bumps defined as 
\begin{align*}
	b(x)=\begin{cases}0.2-(x-(4.5 k-0.75))^2/20&4.5 k-3<x<4.5 k+1.5,\: k=1,2,3,4,5\\ 0 & \text{otherwise}\end{cases}
\end{align*} 
\revp{on the domain $x\in [0,25]$ and $y\in [0,25]$.}
We  add a two dimensional perturbation to this equilibrium flow defined as (see also figure \ref{Fig:SWE_1d_eq_pert}(a))
\begin{equation}\label{pert2d}
h=h^*+0.05e^{-100((x-10)^2+(y-12.5)^2)}. 
\end{equation}

\begin{figure}[H]
	\centering
%	\subfigure[$Steady state$]{\includegraphics[width=0.5\textwidth]{twoD/lake_at_rest/lake_rest_2d_ic}}
%	
	\subfigure[NWB]{\includegraphics[width=0.5\textwidth]{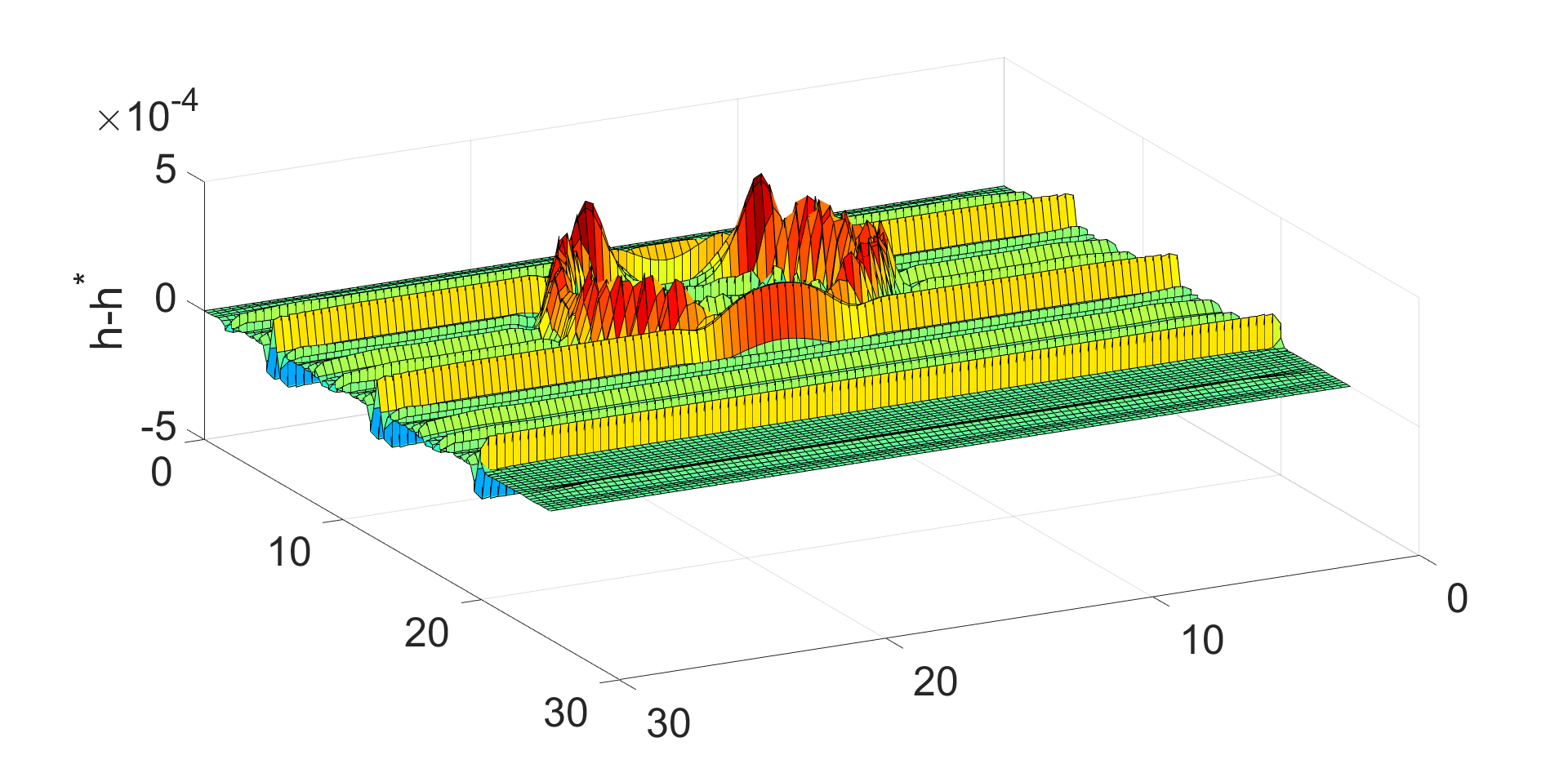}}\subfigure[WB]{\includegraphics[width=0.5\textwidth]{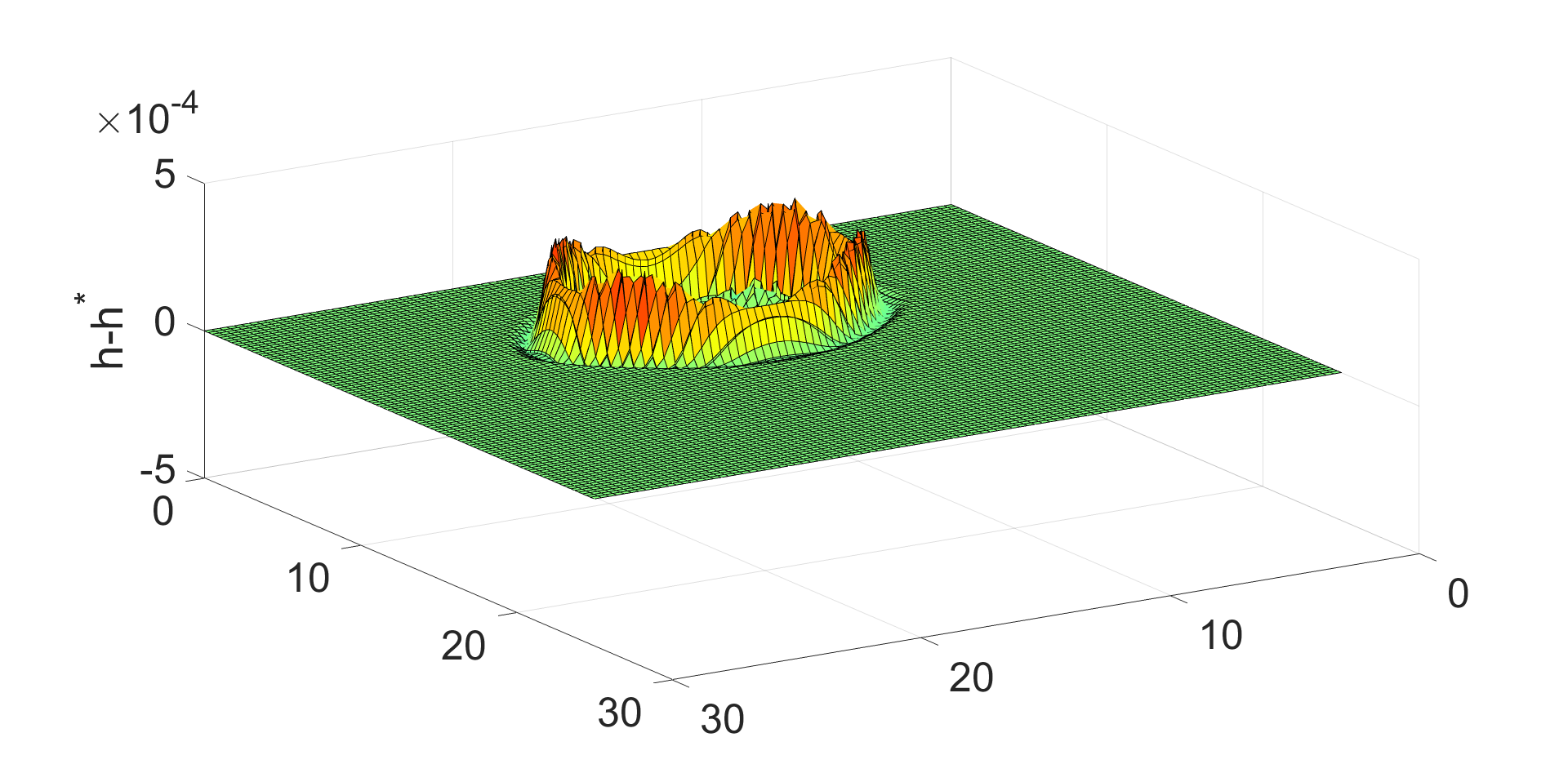}}	
	%	\subfigure[$WB-EC$]{\includegraphics[width=0.5\textwidth]{twoD/oneD_equilibrium/diffheight_2dpert1deq_wbec}}
	
	%	\subfigure[$Entropy$]{\includegraphics[width=0.75\textwidth]{twoD/oneD_equilibrium/entropy_1dequi_2d}}
	\caption{Numerical solution for $h-h^*$ with NWB, WB, WB-EC discontinuous Galerkin scheme}
	\label{Fig:lake_at_rest_2d_pert}
\end{figure}

We evolve the solution until $T=2$s on a $50\times 50$ mesh with both the non-well-balanced and global flux quadrature   schemes with $p=2$. The results are compared    on figure \ref{Fig:lake_at_rest_2d_pert}.  As expected, the non well-balanced formulation introduced spurious numerical oscillations  of size comparable to the amplitude of the perturbation. These perturbations are completely absent in the   well-balanced formulation which
preserves the steady state exactly.

Next we consider a 1D equilibrium where initial condition is calculated using \eqref{steady0} taking  $q_0=5.6865$ and 
 $E_0=54.183738$.   As in the previous example,  we add the perturbation \eqref{pert2d} which is evolved  
 until $T=2$s on a $50\times 50$ mesh with both the non-well-balanced and global flux quadrature   schemes with $p=2$. 
The results are  shown in figure \ref{Fig:SWE_1d_eq_pert}.
The NWB scheme is not able to capture the perturbation accurately as the numerical error is as large as the physical perturbation. 
The global flux quadrature approach allows to obtain a correct preservation of the background flow, 
and a nice evolution for the perturbation. 
 It may be noted that the results of well-balanced scheme with entropy correction are extremely close, so only the first is shown. 

\begin{figure}[H]
	%\centering	
%	\subfigure[Initial condition fr $h+b$] {\includegraphics[width=0.4\textwidth]{twoD/vortex/vortex_ic}}
%	\subfigure[Grid convergence] 
          \begin{center}
	\includegraphics[width=0.35\textwidth]{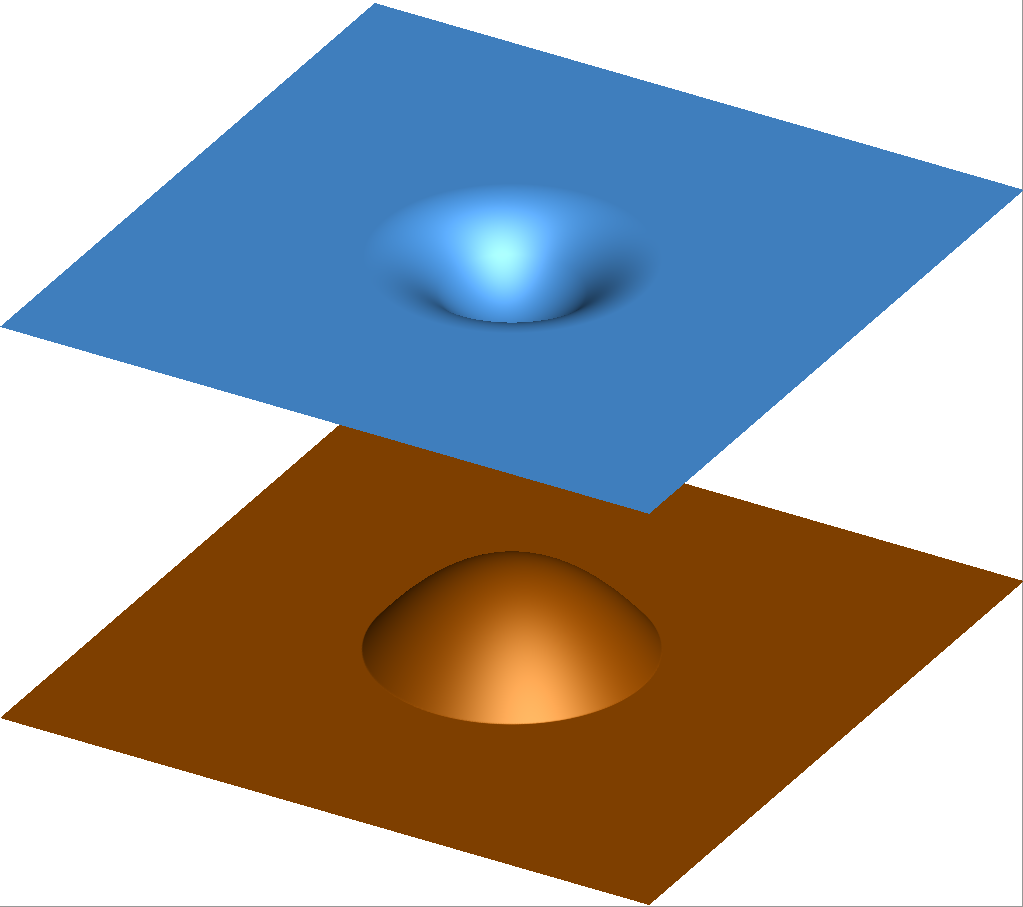}\includegraphics[width=0.6\textwidth]{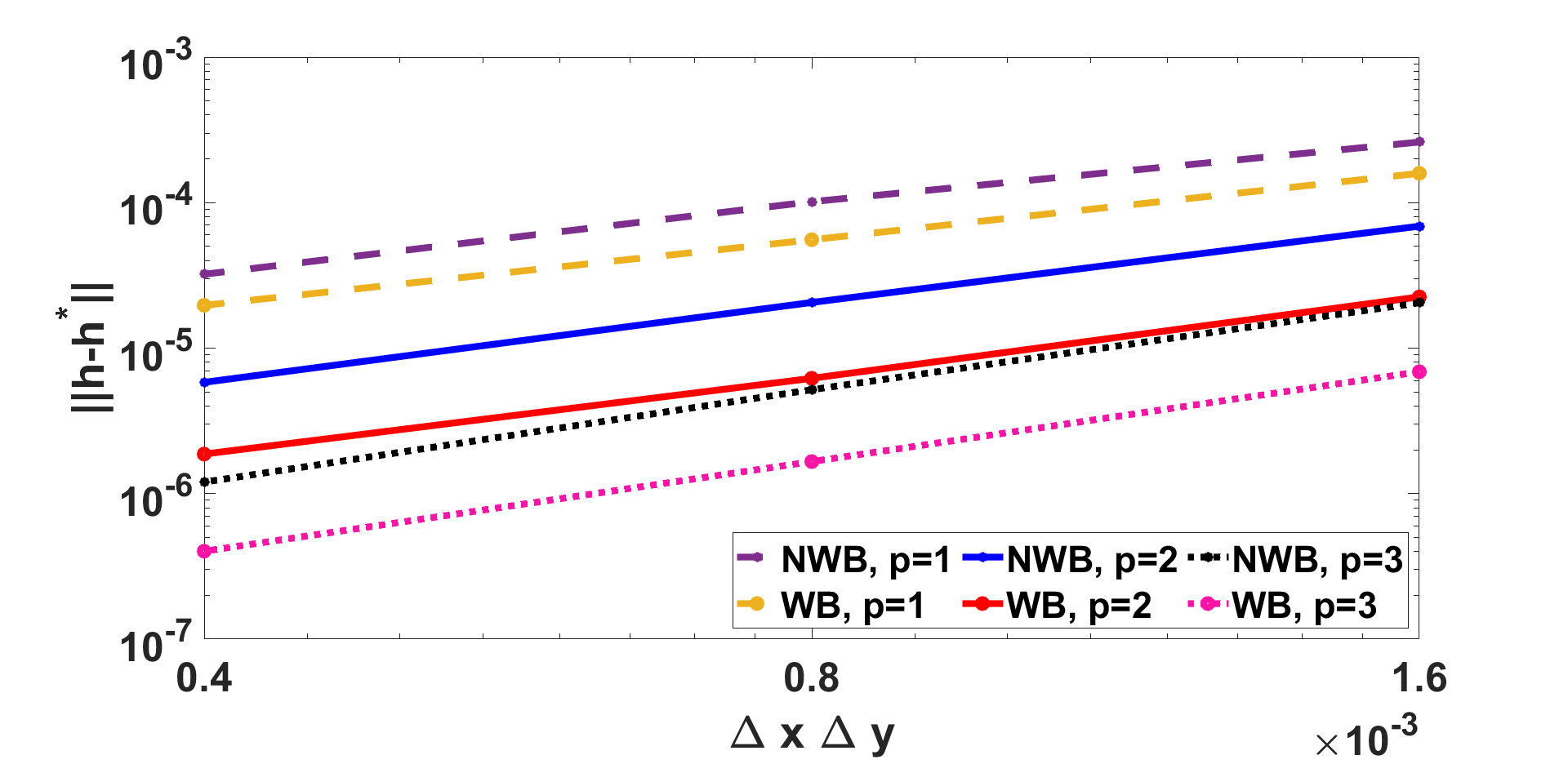}
	\end{center}
%<<<<<<< HEAD
%%	\subfigure[2nd order, $T=1$] {\includegraphics[width=0.33\textwidth]{twoD/vortex/vortex_2ord_t1}}	
%%	\subfigure[3rd order, $T=1$] {\includegraphics[width=0.33\textwidth]{twoD/vortex/vortex_3ord_t1}}
%	\caption{\textcolor{red}{\bf ADD GRID CONVERGENCE with $p=1,2,3$ with WB and NON WB}}
%=======
%	\subfloat[2nd order, $T=1$] {\includegraphics[width=0.33\textwidth]{twoD/vortex/vortex_2ord_t1}}	
%	\subfloat[3rd order, $T=1$] {\includegraphics[width=0.33\textwidth]{twoD/vortex/vortex_3ord_t1}}
	\caption{Steady vortex. Left: scaled 3D view of the  initial solution (free surface and bathymetry). Right:  grid convergence with $p=1,2,3$.}
%>>>>>>> 8633838fbb44d758f3df865c7cef7dc43c507d06
	\label{Fig:SWE_vortex_equi}
\end{figure}

\paragraph{2d stationary vortex with bathymetry} This is a genuinely 2D steady state. It is a variation
of a test often used in literature, see e.g \cite{AUDUSSE20092934} with added bathymetric effects.
The initial condition is given by
\begin{equation}\label{eq:steady_vortex1}
\begin{aligned}
	h^*(r,0)+b(x,y)&=1+\xi^2\begin{cases}
		\frac{5}{2}(1+5\xi^2)r^2 & r\leq \frac{1}{5}\\
		\begin{aligned}
			\frac{1}{10}(1+5\xi^2)&+2r-\frac{3}{10}-\frac{5}{2}r^2\\&+\xi^2[4ln(5r)+\frac{7}{2}-20r+\frac{25}{2}r^2]
		\end{aligned}&\frac{1}{5}<r\leq\frac{2}{5}\\
		\frac{1}{5}(1-10\xi^2+20\xi^2ln(2))& r>\frac{2}{5}
	\end{cases}
	\\
	u^*(x,y,0)&=-\xi y\Upsilon(r),\quad v^*(x,y,0)=\xi x \Upsilon(r), \quad\Upsilon(r)=\begin{cases}
		5 &  r\leq \frac{1}{5}\\\frac{2}{r}-5 & \frac{1}{5}<r\leq\frac{2}{5}\\ 0 &r>\frac{2}{5}
	\end{cases}
	\end{aligned}
\end{equation}
where $r=\sqrt{x^2+y^2}$. For this test we consider $\xi=0.1$ and bathymetry  defined as
$$
b(x,y)=\left\{
\begin{array}{ll}
0.1(1-6.25(x^2+y^2)) \quad & \text{if } (x^2+y^2)<0.16 \\
0 \quad & \text{otherwise }
\end{array}\right.
$$
\revp{on the domain $[-1,1]\times[-1,1]$. The boundary condition is given by $h=1+\frac{1}{5}\xi^2(1-10\xi^2+20\xi^2ln(2))\text{m}, u=v=0\text{m}/\text{s}$ on all the boundaries.}
As scaled visualization of the initial free surface is reported on the left on 
figure \ref{Fig:SWE_vortex_equi}. We report in the same figure the grid convergence
with the global flux quadrature  and with the non well-balanced  method for different degree approximations.
Note that \eqref{eq:steady_vortex1}  only provides $C^0$ continuity.  As shown in \cite{https://doi.org/10.48550/arxiv.2109.10183}
this will limit the convergence attainable with high order schemes to roughly second order, which is the slope obtained here.
The  convergence plot  however also shows that, despite the fact that we do not embed a genuinely two dimensional well-balanced criterion,
the global flux quadrature formulation still provides a considerable decrease in  error especially for  $p=2$ and $p=3$.
In particular, the errors of the WB $p=2$ are quite comparable with those of the  fourth order non well-balanced DGSEM.

%figure \ref{Fig:SWE_vortex} shows the initial conditions and solution for height of water at time $T=1$ and $b(x,y)=0$ and $\xi=0.1$ with $50\times 50 $ grid using 2nd and 3rd order discontinuous Galerkin scheme and the $L_1$ errors are given in Table \ref{Tab:swe_vor}.

%We see from figure \ref{Fig:SWE_vortex_equi} that though we cannot preserve equilibrium flow exactly since the initial conditions do not satisfy the global fluxx equilibrium, the solution with both WB and WBEC scheme is much more accurate than the NWB scheme.

\begin{figure}[H]
	%\centering	
	\subfigure[NWB]{\includegraphics[width=0.5\textwidth]{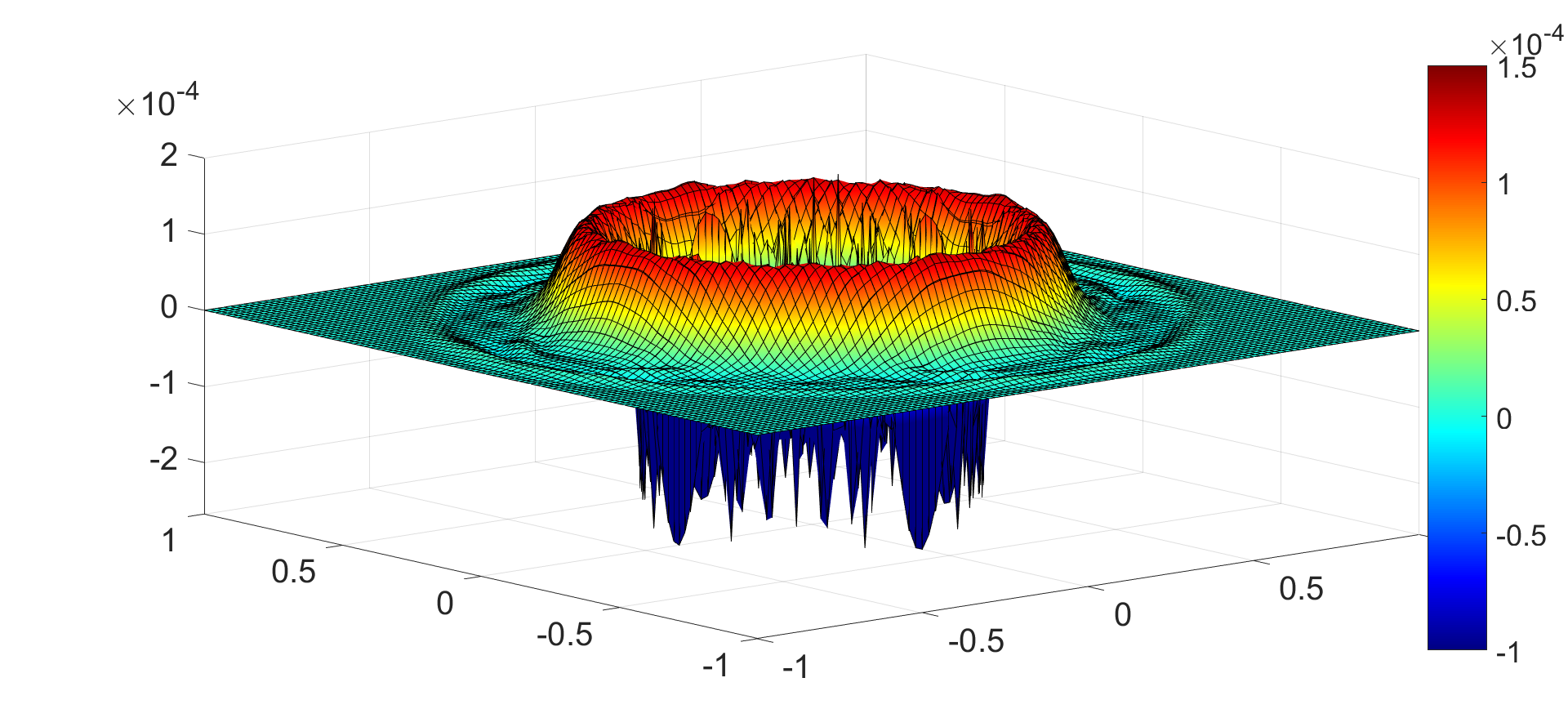}}	
	\subfigure[NWB]{\includegraphics[width=0.325\textwidth]{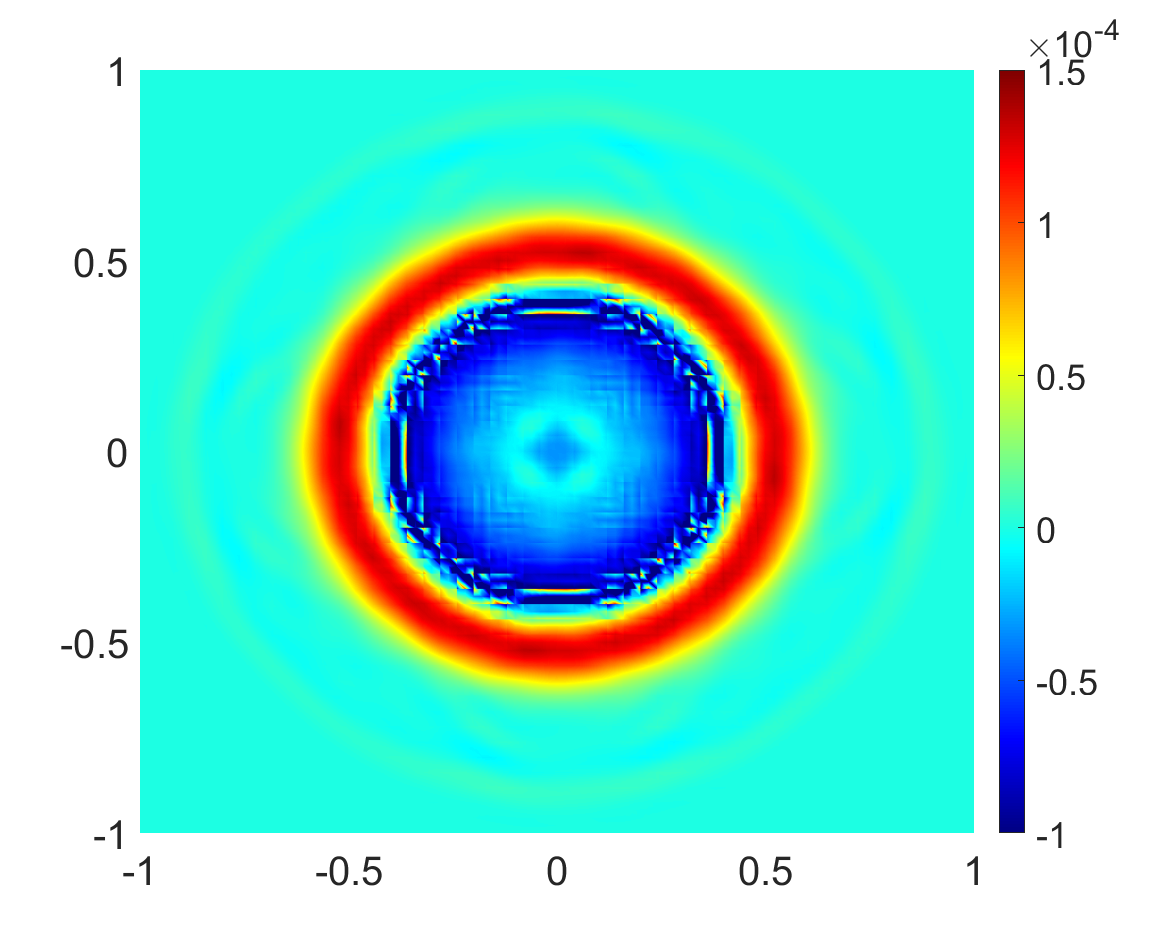}}
	
	\subfigure[WB]{\includegraphics[width=0.5\textwidth]{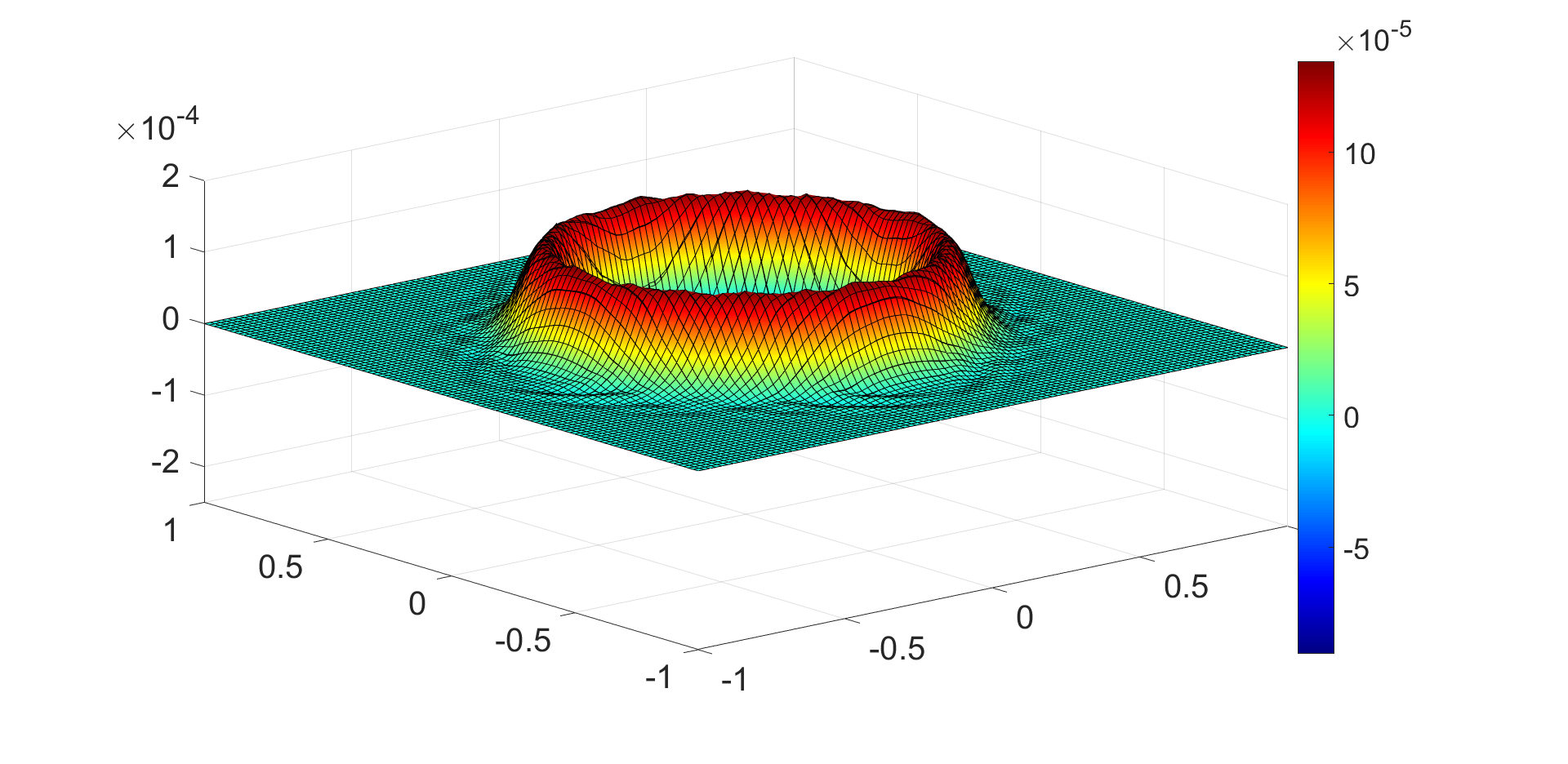}}	
	\subfigure[WB]{\includegraphics[width=0.325\textwidth]{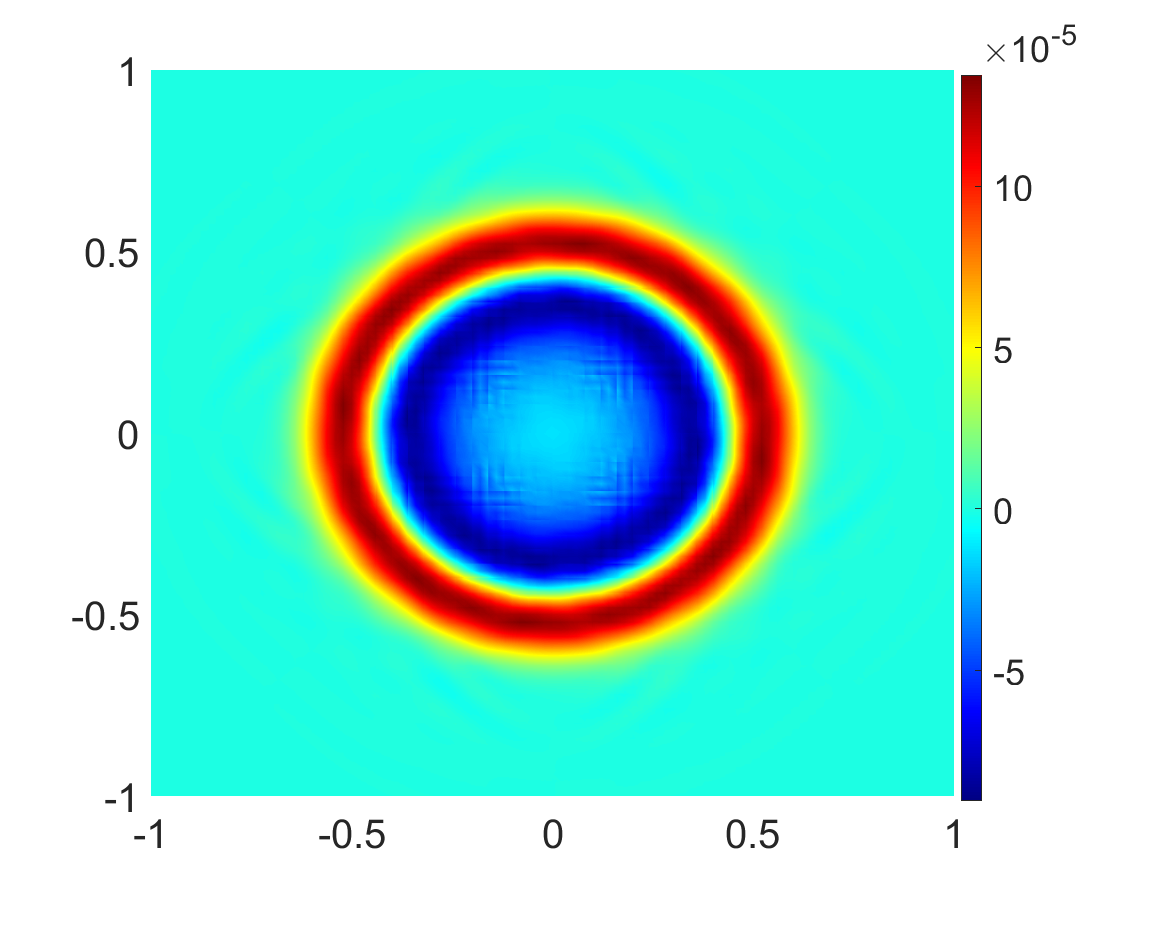}}
	
%<<<<<<< HEAD
%%	\subfigure[$WB-EC$]{\includegraphics[width=0.5\textwidth]{twoD/vortex/vortex_pert_wbec}}	
%%	\subfigure[$WB-EC$]{\includegraphics[width=0.3\textwidth]{twoD/vortex/vortex_pert_wbec_xy}}
%		\caption{Numerical solution for $h-h^*$ with NWB, WB discontinuous Galerkin scheme\textcolor{red}{\bf USE SAME SCALE FOR PERTURBATION? CHANGED}}
%=======
%	\subfloat[$WB-EC$]{\includegraphics[width=0.5\textwidth]{twoD/vortex/vortex_pert_wbec}}	
%	\subfloat[$WB-EC$]{\includegraphics[width=0.3\textwidth]{twoD/vortex/vortex_pert_wbec_xy}}
		\caption{Numerical solution for $h-h^*$ with NWB, WB discontinuous Galerkin scheme}
%>>>>>>> 8633838fbb44d758f3df865c7cef7dc43c507d06
	\label{Fig:SWE_vortex_pert}
\end{figure}

We also test the evolution of a small perturbation added to water height and given by
$$h=h^*+10^{-3}e^{-100(x^2+y^2)} .$$ 
The perturbation is evolved up to $T=0.05$ with the third order  schemes on a  $50\times 50$ grid.
Snapshots of the perturbation at the final time are reported in figure \ref{Fig:SWE_vortex_pert}.
The NWB scheme introduces spurious oscillations in the solution over time, while the well-balanced method
provides a much cleaner resolution of the waves. 
%As before only the results with entropy correction a
% which is not the case with WB and WBEC scheme. The spurious oscillations also lead to an increase in energy of the solution as can be seen in figure \ref{Fig:SWE_vortex_pert-energy}. 

%We can also see that the entropy correction helps improve the energy preservation of the WB scheme over time.

%\begin{figure}[H]
%	%\centering	
%           \includegraphics[width=0.55\textwidth]{twoD/vortex/entropy_vortex}
%	\caption{Energy evolution in time for the    NWB, WB, WB-EC schemes} 
%	\label{Fig:SWE_vortex_pert-energy}
%\end{figure}

 \paragraph{Anticyclonic vortex propagation} This is a moving variant of the previous case.  It  consists of a 
 vortex  propagating westward due to the effect of the variation of the Coriolis coefficient which is modified as
 $$
 \omega = \omega_0 + \beta y
 $$
 to model the  effects curvature effects using  a tangent plane approximations.  The domain a rectangular basin of  $2000 \times 1200$km. The
initial condition is given by a Gaussian distribution of the free surface centered at the origin of the
domain, prescribed together with a velocity field which is in geostrophic balance. We refer to \cite{MR2460785,NavasmoltillaMurillo2018} and references therein 
for details on the test setup, and values of the different parameters. The solution is computed with the third order schemes until a final time corresponding to 8 weeks on a relatively coarse
mesh of $50\times30$km. 

\begin{figure}[H]
	%\centering	Results/twoD/
\hspace{-0.2cm}  \subfigure[NWB]{\includegraphics[width=0.55\textwidth,trim={5cm 0cm 0cm 0cm},clip]{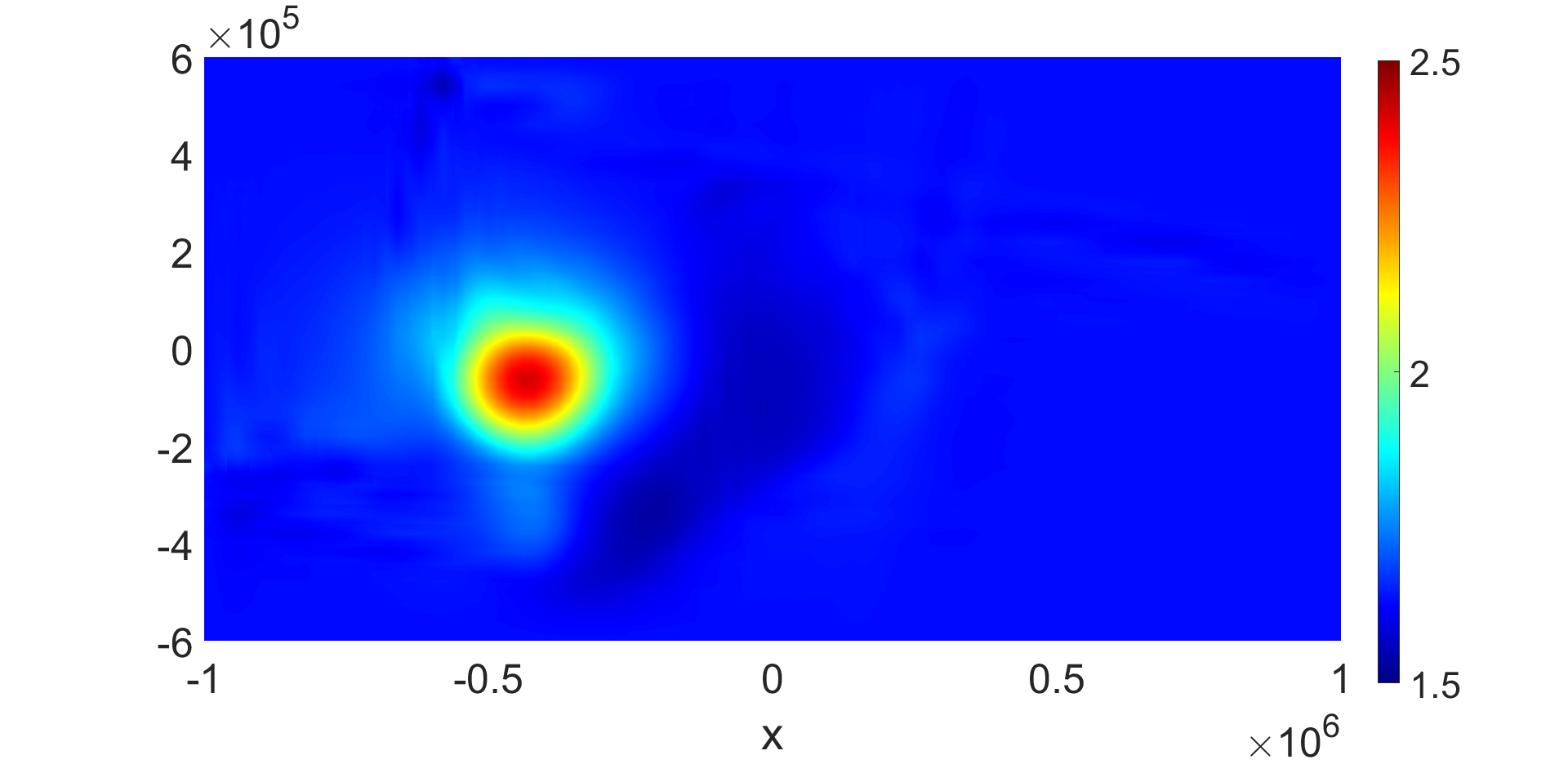}}\hspace{-0.6cm} 
	\subfigure[WB]{\includegraphics[width=0.55\textwidth,trim={5cm 0cm 0cm 0cm},clip]{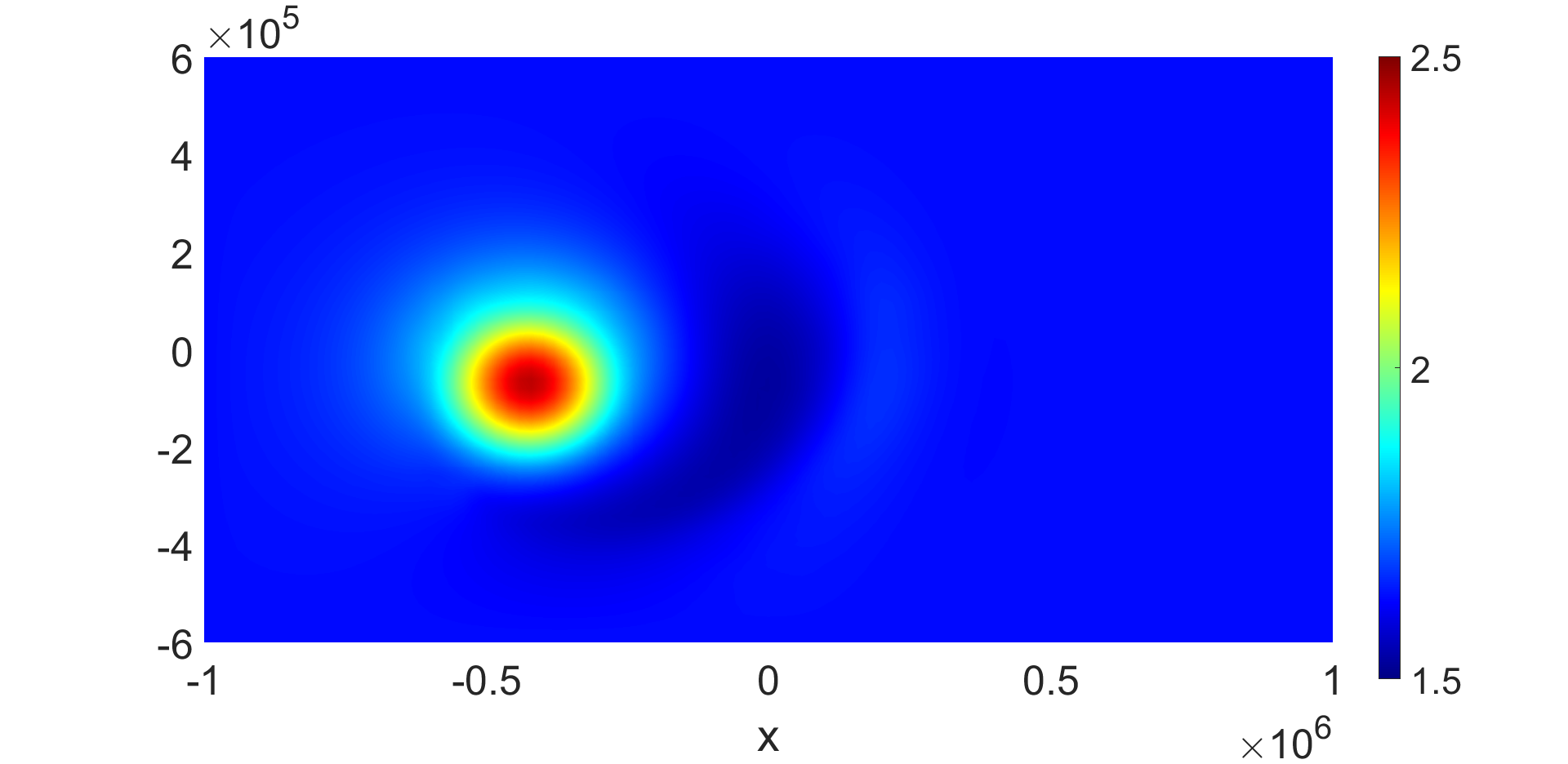}}	
%<<<<<<< HEAD
%%	\subfigure[$WB-EC$]{\includegraphics[width=0.5\textwidth]{twoD/vortex/vortex_pert_wbec}}	
%%	\subfigure[$WB-EC$]{\includegraphics[width=0.3\textwidth]{twoD/vortex/vortex_pert_wbec_xy}}
%		\caption{Numerical solution for $h-h^*$ with NWB, WB discontinuous Galerkin scheme\textcolor{red}{\bf USE SAME SCALE FOR PERTURBATION? CHANGED}}
%=======
%	\subfloat[$WB-EC$]{\includegraphics[width=0.5\textwidth]{twoD/vortex/vortex_pert_wbec}}	
%	\subfloat[$WB-EC$]{\includegraphics[width=0.3\textwidth]{twoD/vortex/vortex_pert_wbec_xy}}
		\caption{Anticyclonic vortex propagation.}
%>>>>>>> 8633838fbb44d758f3df865c7cef7dc43c507d06
	\label{Fig:SWE_vortex_mov}
\end{figure}

We plot the free surface contour levels obtained with the non well-balanced and global flux quadrature methods on figure
\ref{Fig:SWE_vortex_mov}. The results  confirm the observations made for the steady vortex:   the well-balanced scheme provides a much cleaner solution, 
while some spurious waves around the vortex are clearly visible in the non well-balanced result.

\paragraph{Geostrophic adjustment in 1d and 2d}

\begin{figure}[H]
%	\centering
%	\subfigure[$T=0$]{\includegraphics[width=0.75\textwidth]{oneD/Coriolis/sweGE0}}	
%	
%<<<<<<< HEAD
%	\hspace{-0.5cm}\subfigure[$T=1$]{\includegraphics[width=0.65\textwidth]{oneD/Coriolis/sweGE1}}\hspace{-0.75cm}
%	\subfigure[$T=2$]{\includegraphics[width=0.65\textwidth]{oneD/Coriolis/sweGE2}}
%	
%	\hspace{-0.5cm}\subfigure[$T=3$]{\includegraphics[width=0.65\textwidth]{oneD/Coriolis/sweGE3}}\hspace{-0.75cm}
%	\subfigure[$T=4$]{\includegraphics[width=0.65\textwidth]{oneD/Coriolis/sweGE4}}
%	\caption{Solution with 2nd and 3rd order well-balanced DG scheme at $T=\frac{1}{2}\frac{\pi}{f},\frac{\pi}{f},\frac{3}{2}\frac{\pi}{f}, 2\frac{\pi}{f}$}
%=======
\centering\subfigure[$T=\frac{1}{2}\frac{\pi}{\omega}$]{\includegraphics[width=0.5\textwidth]{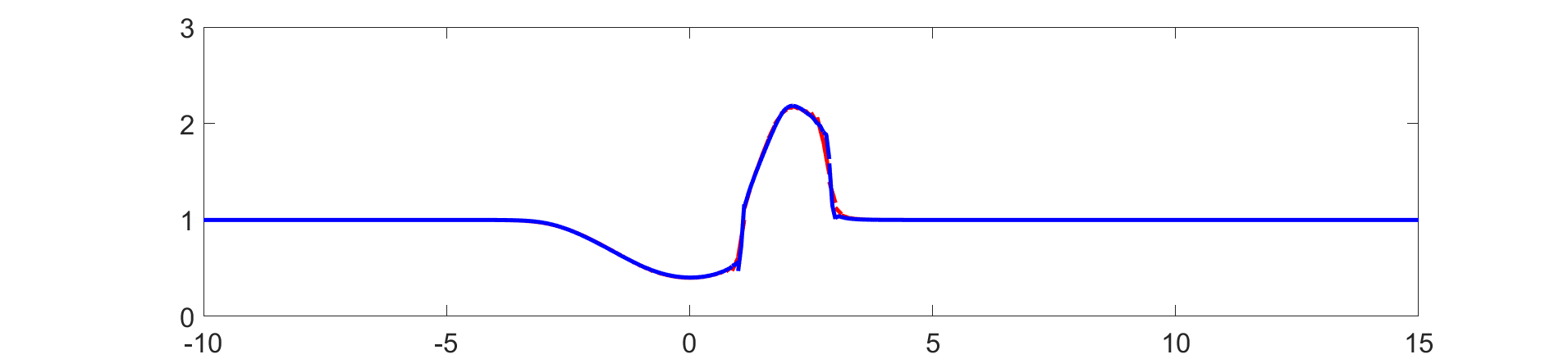}}\hspace{-0.75cm}
	\subfigure[$T=\frac{\pi}{\omega}$]{\includegraphics[width=0.5\textwidth]{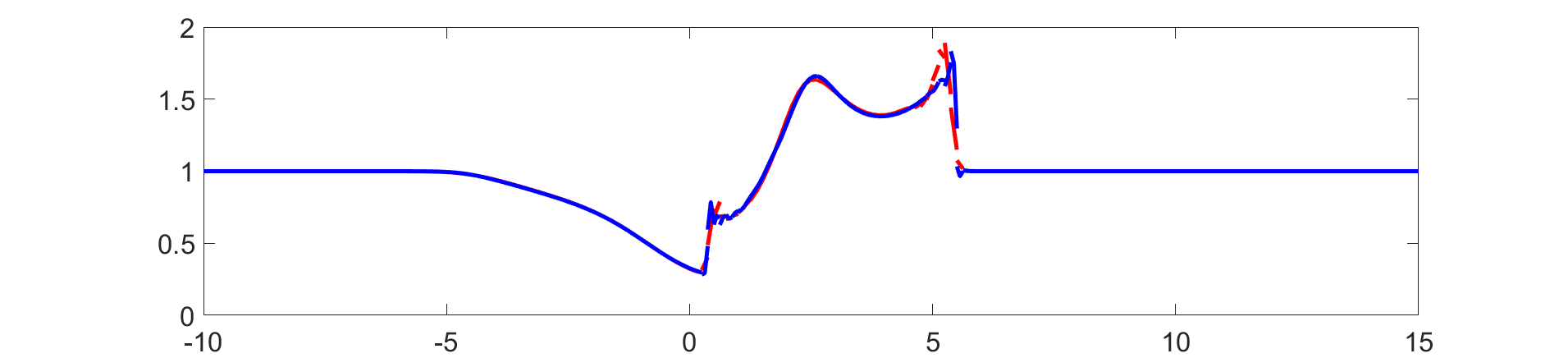}}
	
\centering\subfigure[$T=\frac{3}{2}\frac{\pi}{\omega}$]{\includegraphics[width=0.5\textwidth]{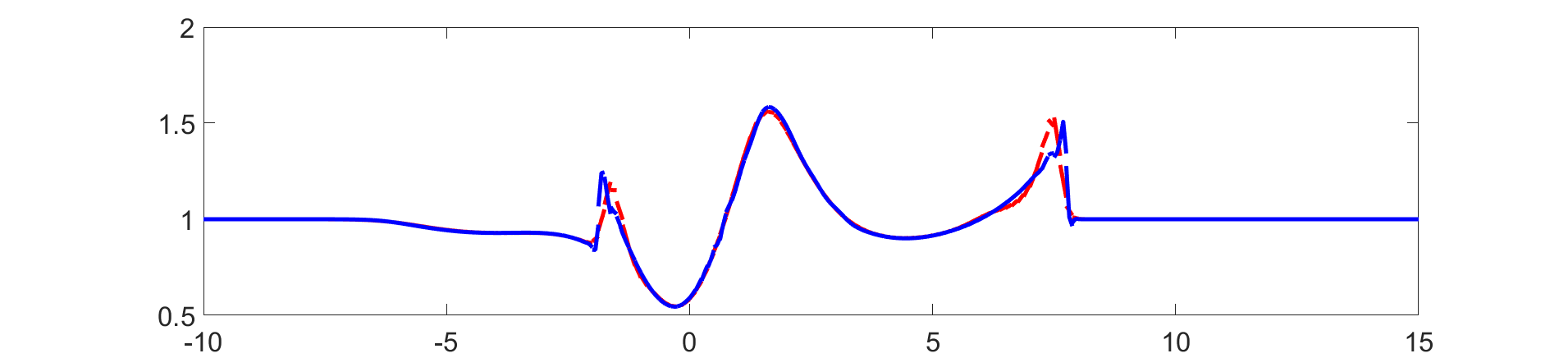}}\hspace{-0.75cm}
	\subfigure[$T=2\frac{\pi}{\omega}$]{\includegraphics[width=0.5\textwidth]{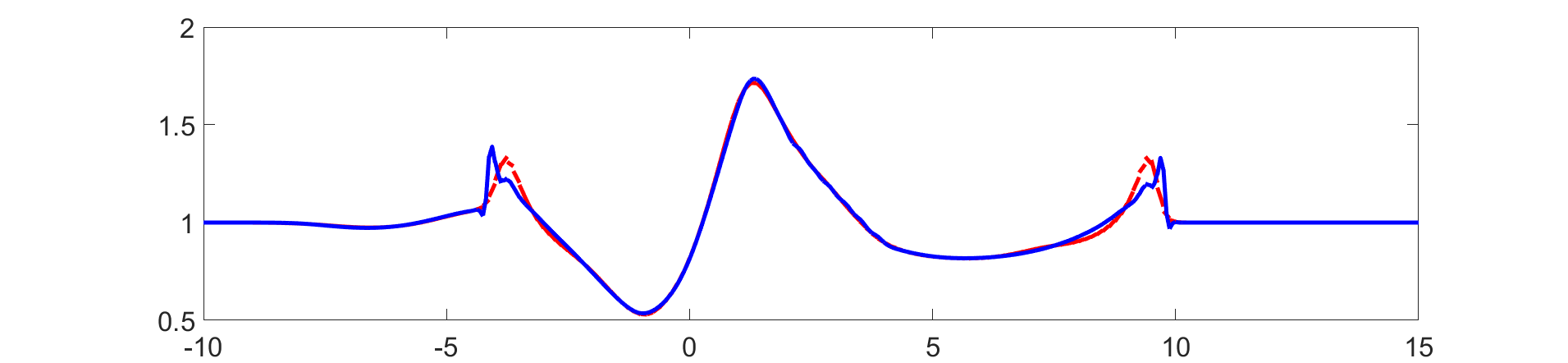}}
	\caption{Solution with 2nd and 3rd order well-balanced DG scheme at $T=\frac{1}{2}\frac{\pi}{\omega},\frac{\pi}{\omega},\frac{3}{2}\frac{\pi}{\omega}, 2\frac{\pi}{\omega}$}
%>>>>>>> 8633838fbb44d758f3df865c7cef7dc43c507d06
	\label{fig:swe_GE}
\end{figure}

We consider now  the simulation of the complex wave dynamics of the geostrophic adjustment in one and two space dimensions.
We start from   the  one dimensional numerical test for geostrophic adjustment used in \cite{bouchut_sommer_zeitlin_2004}. The initial conditions are given by
\begin{align}
	h(x,0)=1,\quad u(x,0)=0\\
	v(x,0)=2\frac{(1+\text{tanh}(4\frac{x}{L}+2))(1-\text{tanh}(4\frac{x}{L}-2))}{(1+\text{tanh}(2))^2}
\end{align}
with $ L=2\text{m}$, and   $\omega=1\text{s}^{-1},\: g=1\text{m}/\text{s}^2$ in the rotating shallow water equations. The initial disturbance in the momentum leads to two fast inertial gravity waves and shock formation. figure \ref{fig:swe_GE} gives the solution for free water surface at time $T=\frac{1}{2}\frac{\pi}{\omega},\frac{\pi}{\omega},\frac{3}{2}\frac{\pi}{\omega}, 2\frac{\pi}{\omega}$ computed using a second and third order well-balanced discontinuous Galerkin scheme with a grid of $N=200$ in the domain $[-10,15]$m. \revp{The boundary condition is given by $(h,u,v)(-10,t)=(h,u,v)(15,t)=(1,0,0)$.}

From the solution we see that both the schemes are able to capture the two fast waves accurately in time, with sharp front almost
free of oscillations despite the absence of any limiter. The third order scheme is less diffusive and is able to capture the shock front more accurately.\\

\begin{figure}[H]
 \centering\includegraphics[width=0.5\textwidth]{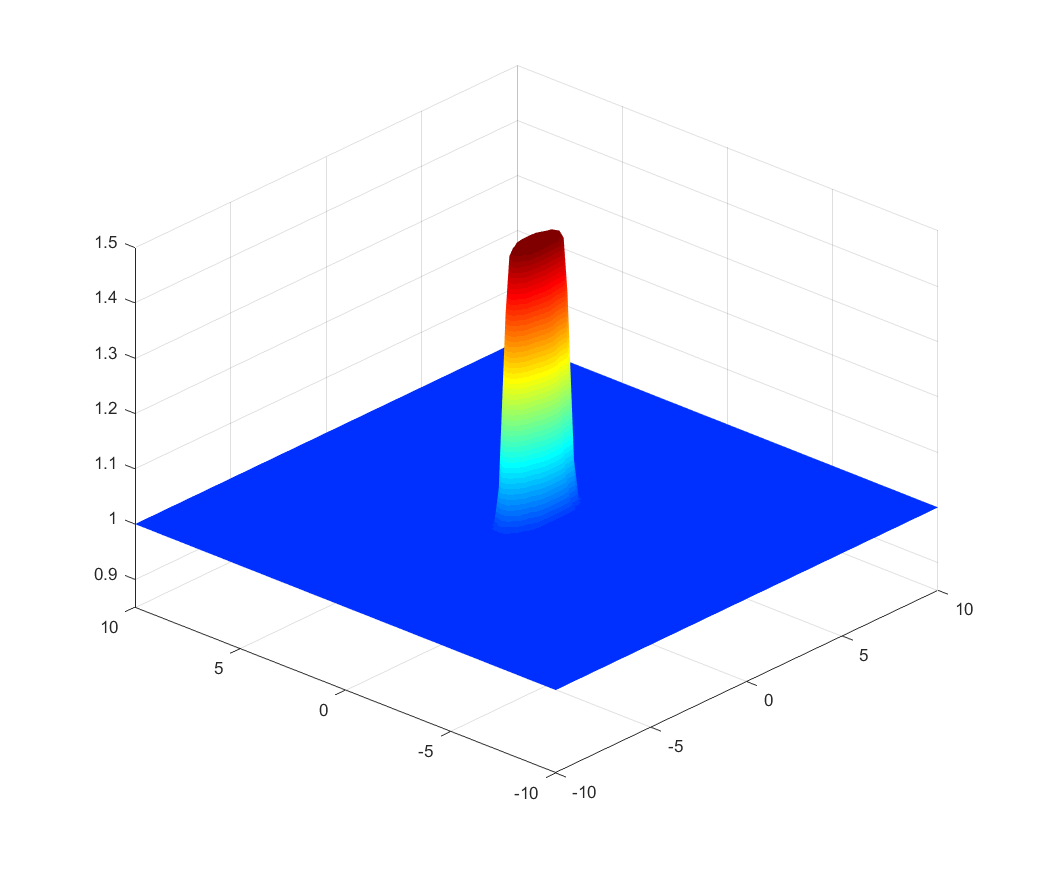} 
	\caption{Initial free surface for the 2d geostrophic adjustment problem by  \cite{MR2460785}}
	\label{Fig:SWE_GE0}
\end{figure}

\begin{figure}[H]
	\centering
%	\subfigure[$T=0$]{\includegraphics[width=0.5\textwidth]{twoD/geostrophic/swee0}}
	\subfigure[$T=4$]{\includegraphics[width=0.5\textwidth,trim={0 1cm 0 10cm},clip]{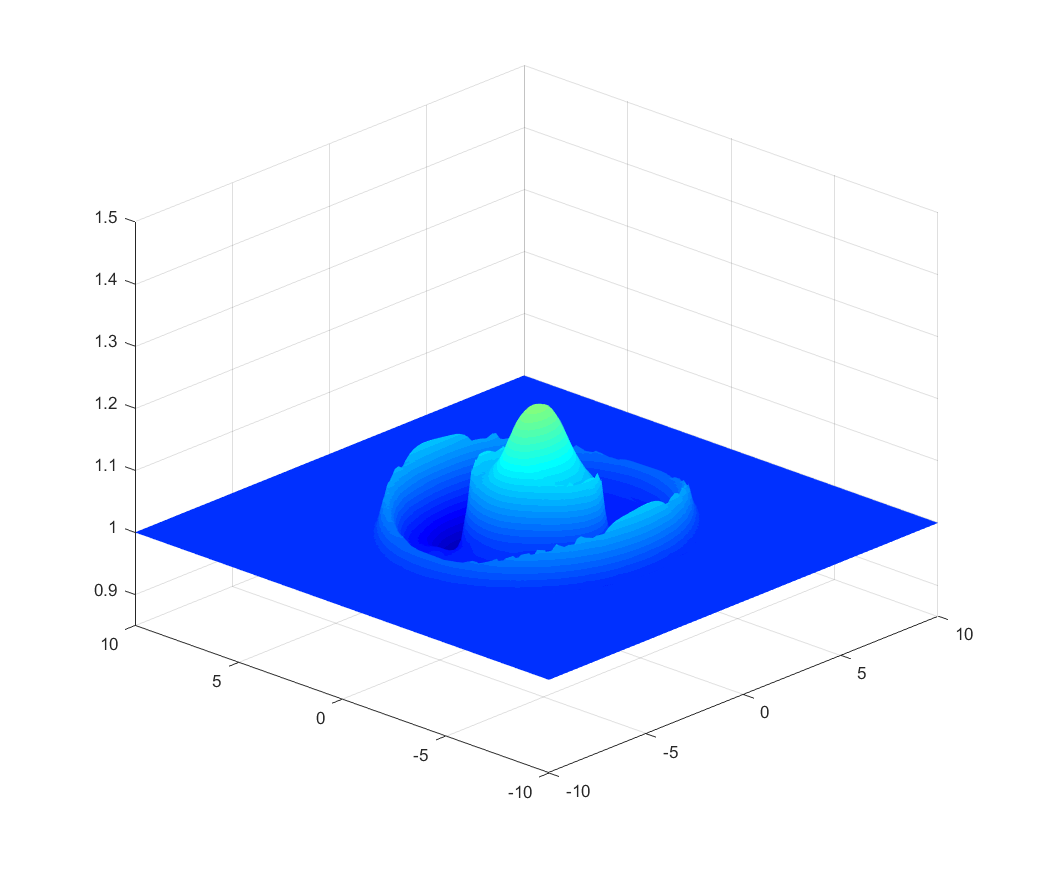}}\subfigure[$T=8$]{\includegraphics[width=0.5\textwidth,trim={0 1cm 0 10cm},clip]{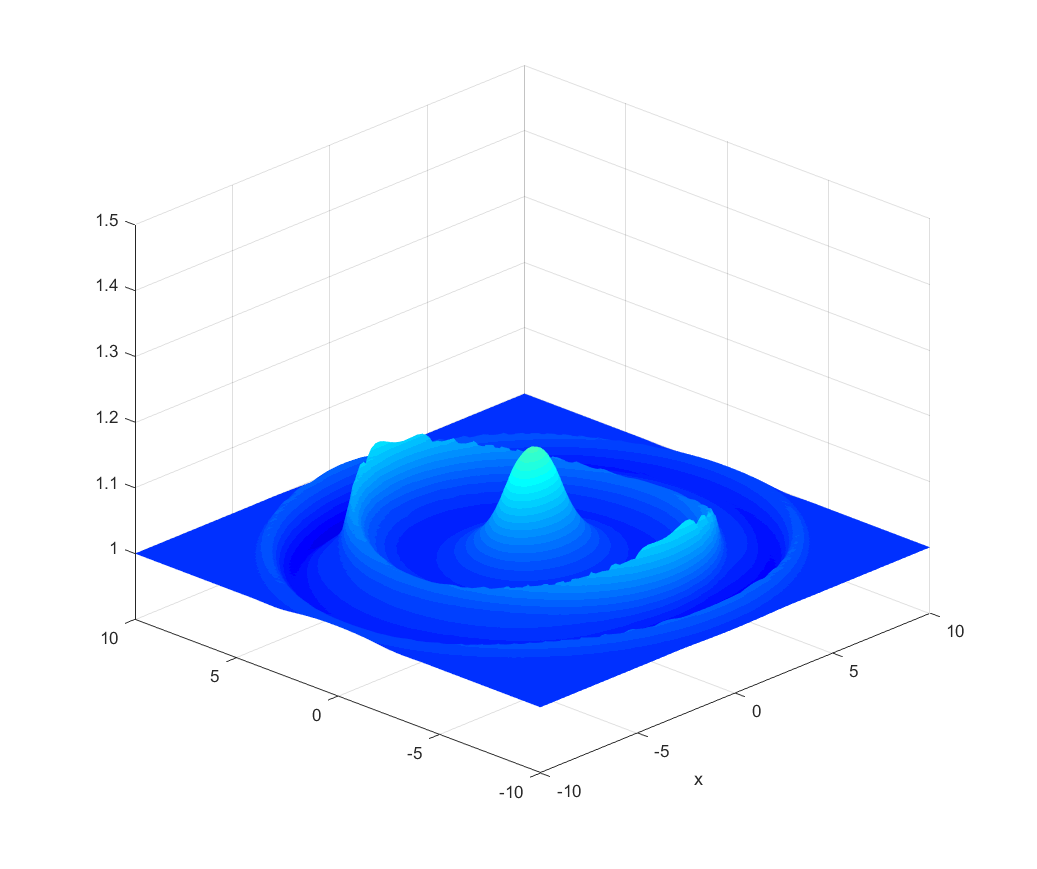}}

	\subfigure[$T=12$]{\includegraphics[width=0.5\textwidth,trim={0 1cm 0 10cm},clip]{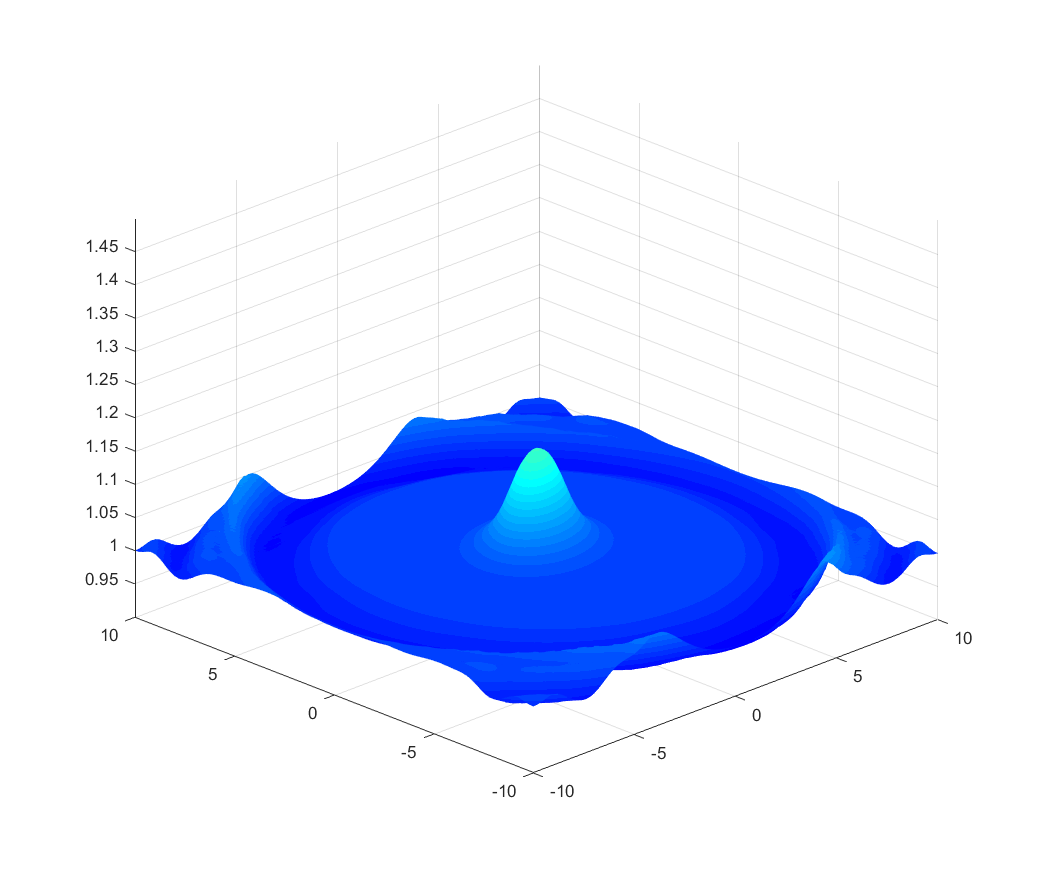}}\subfigure[$T=20$]{\includegraphics[width=0.5\textwidth,trim={0 1cm 0 10cm},clip]{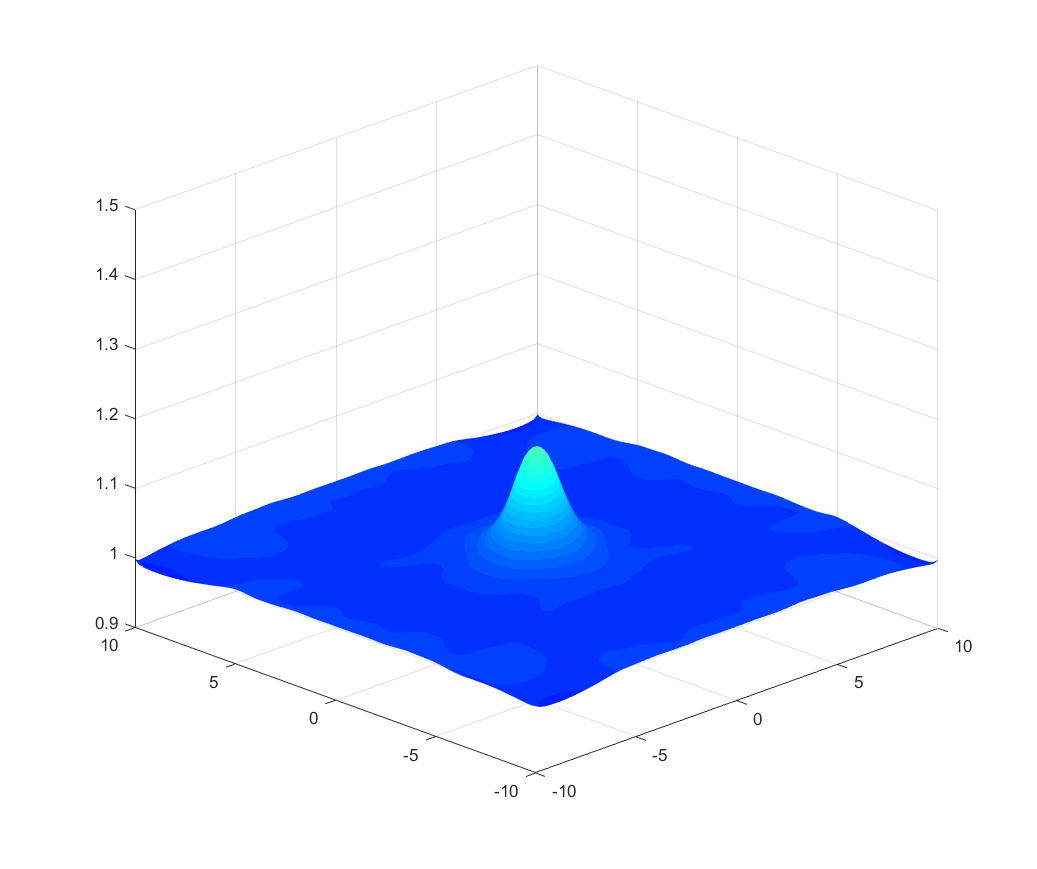}}
	
	\caption{2d geostrophic adjustment problem. Water surface at $T=4,8,12,20$}
	\label{Fig:SWE_GE}
\end{figure}

Next, consider the  two-dimensional extension of this test,  proposed  in  \cite{MR2460785}. In this case    the initial conditions are as follows
\begin{align}
	h(x,y,0)&=1+0.25*(1-\text{tanh}(\frac{\sqrt{2.5 x^2+y^2/(2.5)}-1}{0.1}))\\
	u(x,y,0)&=v(x,y,0)=0
\end{align}
on the spatial domain $[-10,10]\text{m}\times[-10,10]\text{m}$. The initial free surface is visualized in figure \ref{Fig:SWE_GE0},
while figure \ref{Fig:SWE_GE} reports the free surface   at times $T=4,8,12,20$s  obtained with the third order well-balanced method using a  $50\times 50$ mesh.

The solution for this problem is given by two shock waves which move radially away from the center leaving behind some mass which rotates at the initial position. We see that the well-balanced scheme is able to capture these shock waves accurately and the final equilibrium is achieved without any unphysical effects due to the scheme.

\paragraph{Kelvin front generation on the equatorial $\beta$-plane}

\begin{figure}[H]
%	\centering
\centering\subfigure[2nd order]{\includegraphics[width=0.5\textwidth]{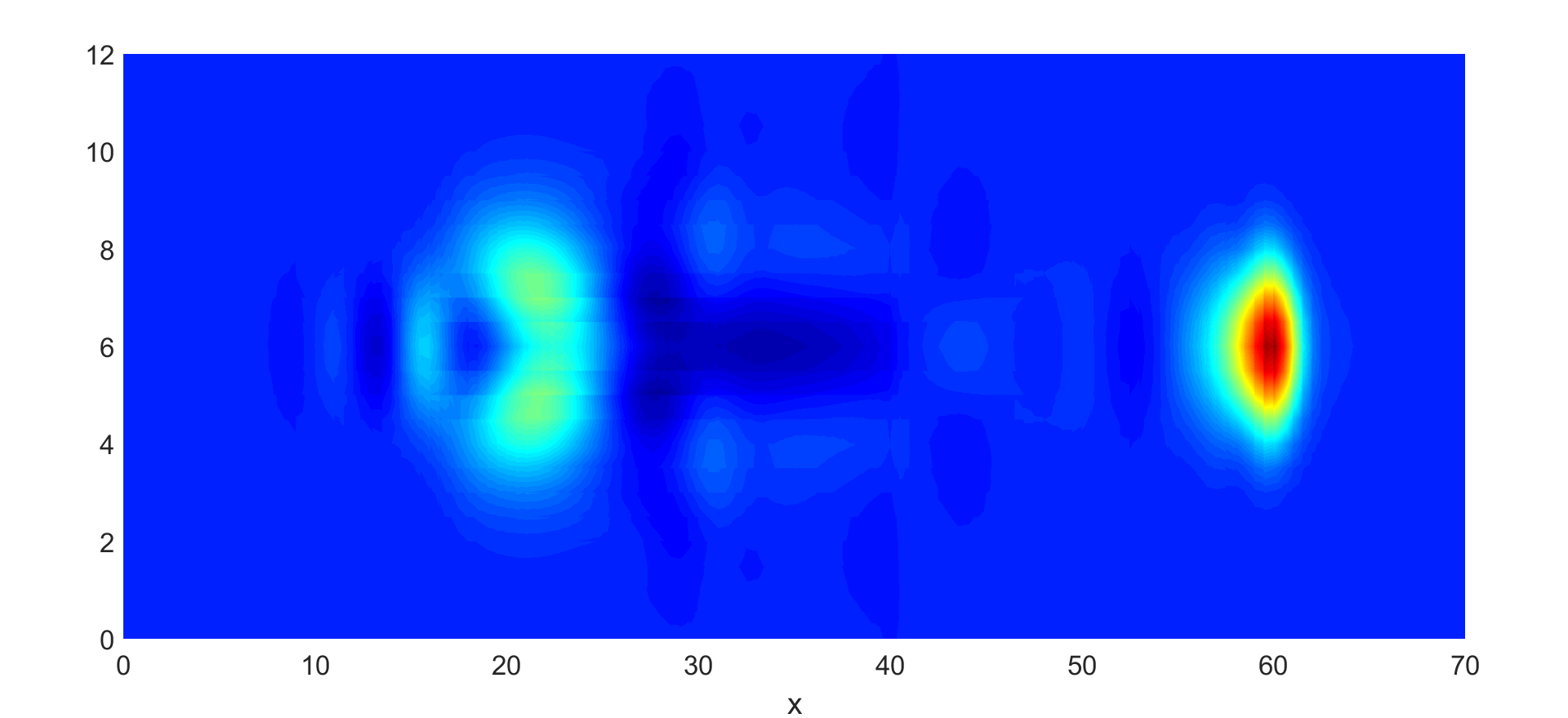}}\subfigure[3rd order]{\includegraphics[width=0.5\textwidth]{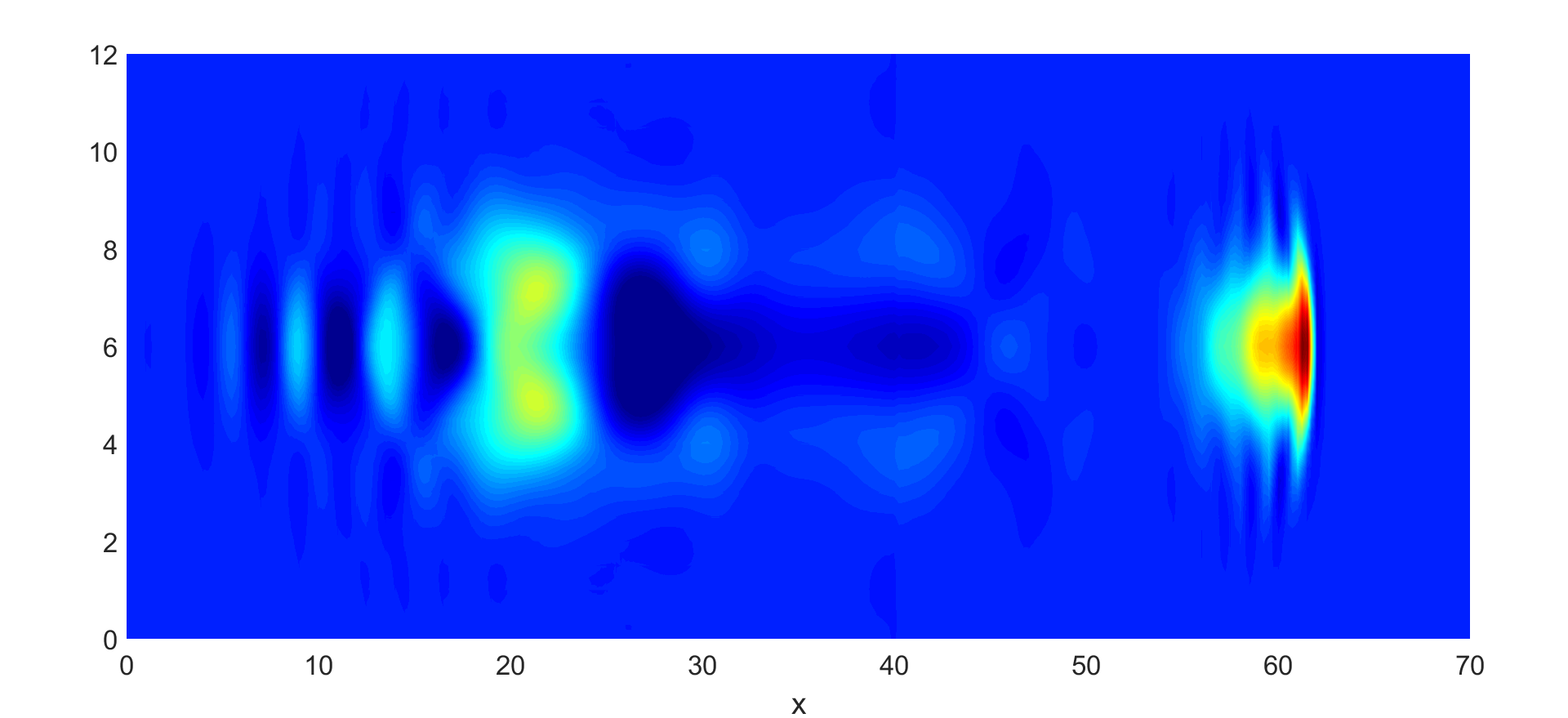}}
	
%	\subfigure[$Entropy$]{\includegraphics[width=0.9\textwidth]{twoD/kelvin/entropy_kelvin}}
	\caption{Numerical solution for $h+b$ with 2nd and 3rd order discontinuous Galerkin scheme}
	\label{Fig:SWE_kel}
\end{figure}

As a last test, we show the capabilities of the scheme proposed to capture the generation of short secondary waves. The test considered involves    the formation and propagation of Kelvin and Rossby waves
with a generation of a  Kelvin front and secondary  Poincar\'e  waves. The set up is the  the same as in   \cite{NavasmoltillaMurillo2018,fedorov00}.
The computational domain $[0,70]\times [0,12]$ is  divided into $140\times 24$ cells. The initial condition and bottom topography are given by,

\begin{align*}
	h(x,y,0)=2-b(x,y)+0.8exp\big(-\frac{(x-30)^2+(y-6)^2}{3}\big)\\
	b(x,y)=\begin{cases}
		0&x\leq 40\\ 0.025x-1& x>40
	\end{cases}
\end{align*}
with $g=1$m/s$^2$. As for the moving vortex case, the Coriolis coefficient is modified to model the effects of curvature. In this case, following  \cite{NavasmoltillaMurillo2018} we set  $\omega(y)=y-6$. 
Figure \ref{Fig:SWE_kel} shows the solution at $T=20$ with a second and third order well-balanced discontinuous Galerkin scheme.

The solution consists of the short wavelength Kelvin waves which carry energy eastwards and the long wavelength Rossby waves carrying energy westward. The nonlinear Kelvin waves steepen and eventually break forming a broken wave front propagating eastward which leads to generation of the secondary Poincaré waves.
From the solution in figure \ref{Fig:SWE_kel}, we see that both the 2nd and 3rd order scheme are able to capture the general physical behavior of both the Kelvin and Rossby waves. However only the third order scheme is able to capture the Poincar\'e waves. Indeed, the solution obtained with the second order method
gives a Kelvin front too diffused to trigger the generation of these secondary waves.

\section{Conclusions and perspectives} \label{sec:Conclusion}

In this work, we have studied a new formulation to construct numerical methods verifying a  discrete fully well-balanced  criterion agnostic
of \revPO{any specific  form of  steady equilibrium}. In one space dimension, the underlying discrete well-balanced criterion is based on an equivalence between 
\revPO{the discontinuous Galerkin spectral element method and a  superconvergent Gauss collocation  integrator 
applied to an  ODE for the flux. In practice, this property is obtained by means of a simple}
%the local differential problem (ODE) defining the steady solution, and the non-local integral equation obtained from the steady ODE.  
%
%%The discrete well-balanced criterion proposed is based on the equivalence of independent discretizations of the  two problems using the same data. 
%The application of this idea to DGSEM schemes using Gauss-Lobatto interpolation shows that 
%with a simple, and yet non-straightforward,
 modification of the integral of the source term.
%  we obtain an equivalence with 
%Gauss collocation methods for the integral equation. 
  This %provides a 
  discrete well-balanced DGSEM approach, \revPO{which we} named  global flux quadrature based DGSEM, \revPO{allows} %for which we have
   a clear characterization of the discrete steady state, with provable superconvergence estimates to general  analytical steady states. 

In addition, a cell correction term is proposed to verify cell entropy balance. The accuracy and the consistency of  this term with the discrete well-balanced
condition proposed is characterized. The numerical benchmarks for the shallow water equations in one space dimension confirm all the theoretical 
properties on a wide variety of equilibria, involving several different source terms. 

%Numerical tests in two space dimensions show  that the generalization
%of the global flux quadrature 
%proposed allows error reductions also on smooth vortex like solutions.
Applications on complex problems involving sharp fronts and complex wave patterns confirm the accuracy and show the robustness of the   method. 

Several perspectives are being explored including the design of non-linear variants allowing to handle genuinely discontinuous solutions,
other discretization techniques (e.g. finite volumes), and more general formulations in particular for   time-dependent and  multidimensional problems.

 \appendix

 \section{Proof of the consistency estimate of Proposition \ref{prop_entropy}}
 \revp{We prove the last part of the proposition namely and in particular  the consistency  estimate $|\mathcal{E}|=\mathcal{O}(h^{p+1})$. To this end,
 following    \cite{AR:17}, we consider a smooth exact solution $U^e\in C^p$, and its  projection
onto the local finite element space $U^e_{\hh}$, as well as a smooth compactly supported 
test function $v\in C^{k}(\Omega)$, with $k\ge 1$, with $v_{\hh}$ denoting its  projection
onto the local finite element space.  Proceeding as in \cite{AR:17}
we then formally replace $U_{\hh}$ by  $U^e_{\hh}$ in the scheme. The reminder is a local measure of the consistency error.
We now multiply the $i$-th nodal  reminder  in each element $K$ by the projected nodal value of the test function $v_i$,
and define the consistency error as in    \eqref{eta-corr6} by summing up over all degrees of freedom $i$,
and over all elements $K$.
Then  as in \cite{AR:17} first we   rewrite the error in weak form and then perform and estimation of each term appearing in 
the expression. 
Using  the definitions of the fluctuations \eqref{eq:dgsem-gf9} appearing in the error   \eqref{eta-corr6},
of the artificial viscosity term $\mathcal{D}_i^K$ (equation \eqref{eta-corr0}),  invoking the  
 SBP property of the DGSEM method, and including a  rest
to account for the inexactness of the underlying  Gauss-Lobatto quadrature,
we can write the error  \eqref{eta-corr6} in the following form
   \begin{equation}\label{eq:proof-e0}
 \begin{split}
 \mathcal{E}= \sum_K  \int_K  v_{\hh}\partial_tU^e_{\hh}  
   -&\sum_K  \int_K   \partial_x v_{\hh} F_{\hh}(U^e_{\hh})     
   - \sum_K  \int_K  \partial_x v_{\hh} R^e_{\hh}    + \sum_f [\![v_{\hh}]\!] \hat G_{\hh}(U^e_{\hh})     
   \\
+&   \sum_K \mathcal{R}_{K}^{\mathsf{GL}}
 + \sum_K \alpha_K(U^e_{\hh})\int\limits_K\partial_x v_{\hh} A_0(U^e_{\hh})\partial_x U^e_{\hh},
\end{split}
 \end{equation}
 where the term involving the jump of the test function is a result of the locality of the finite element projection. }
\revPO{We now consider the PDE applied to the exact solution $\partial_tU^e + \partial_xG^e =0$,
and    in each element we  multiply t by $v_{\hh}$ and integrate by parts to get
 $$
 \int_K  v_{\hh}\partial_tU^e - \int_K    \partial_x v_{\hh} G^e      + 
(v_{\hh} G^e)^R -(v_{\hh} G^e)^L=0 ,
 $$
 where $G^e=F^e-R^e$ with    source flux $R^e$ is obtained by exact integration  of  
 (assuming the initial state on the left hand of the domain to be zero)
 \begin{equation}
 R^e = \int_{x_ 0}^xS(U^e,s).%ds
  \end{equation}
 Summing up over element and subtracting from  \eqref{eq:proof-e0}, and  omitting the boundary conditions
 (due to the compactness of $v$),  we can write   %on each element
% \begin{equation}
% \begin{split}
% \mathcal{E}=\int\limits_{\Omega} v_{\hh}\partial_tU^e_{\hh}  
%   -&\int\limits_{\Omega}  \partial_x v_{\hh} F_{\hh}(U^e_{\hh})     
%   -   \int\limits_{\Omega}  \partial_x v_{\hh} R^e_{\hh}   
%   \\
%+&   \sum_K \mathcal{R}_{K}^{\mathsf{GL}}
% + \sum_K \alpha_K(U^e_{\hh})\int\limits_K\partial_x v_{\hh} A_0(U^e_{\hh})\partial_x U^e_{\hh}
%\end{split}
% \end{equation}}
%
%, and where we 
%remark that  all the  face  terms have cancelled out due to the continuity of the  smooth exact solution $U^e$,
%and to the sum over the elements.   
%Using the continuity of  the exact solution, the smoothness of the test function $v$, 
%and thus the continuity of its interpolated approximation $v_{\hh}$, we can also  write   
 \begin{equation}
 \begin{split}
 \mathcal{E}=\overbrace{\int\limits_{\Omega} v_{\hh}\partial_t(U^e_{\hh} -U^e)  }^{\mathrm{I}}
   -&\overbrace{\int\limits_{\Omega}  \partial_x v_{\hh} (F_{\hh}(U^e_{\hh})-F(U^e)) }^{\mathrm{II}}  -   
   \overbrace{\int\limits_{\Omega}  \partial_x v_{\hh} (R^e_{\hh}  - R^e)  }^{\mathrm{III}}\\
\underbrace{\sum_f [\![v_{\hh}]\!] (\hat G_{\hh}(U^e_{\hh})  -G^e) }_{\mathrm{IV}}     +& 
   \sum_K \mathcal{R}_{K}^{\mathsf{GL}} + \sum_K \alpha_K(U^e_{\hh})\int\limits_K\partial_x v_{\hh} A_0(U^e_{\hh})\partial_x U^e_{\hh}
\end{split}
 \end{equation}}
  We now proceed to a term by term estimate. Following \cite{AR:17}, we start from noting that  for a test function $v\in C^{k}(\Omega)$ with   $k \ge 1$ we can readily  claim that 
    \begin{equation}
  \|v_{\hh}\|\le C_1 < \infty\,,\;\;
    \|\partial_x v_{\hh}\|\le C_2 < \infty\,,\;\;\
    \revPO{ [\![v_{\hh}]\!] \le C_3 \max(\hh^2,\hh^{p+1}) }
      \end{equation}
 having omitted the $L^2_{\Omega}$ subscript from the norms to simplify the notation.  
This allows to bound terms $\mathrm{I}$ and $\mathrm{II}$  by the  approximation error as
    \begin{equation}
  |\mathrm{I}| \le \tilde C_{\mathrm{I}}(\Omega,\|v\|) \,\hh^{p+1} \,,\;\;
 | \mathrm{II}| \le \tilde C_{\mathrm{II}}(\Omega,\|\partial_xv\|) \,\hh^{p+1}  
      \end{equation}
      \revPO{The term $\mathrm{IV}$ can also be easily controlled by the approximation error by using the form
      of the numerical flux \eqref{eq:dgsem-gf10} which allows to write (using $^+$ and $^-$ for quantities on the 
      two sides of an element face) 
\begin{equation*}
   \begin{split}
  | \mathrm{IV}|  =&|\sum_f\left\{ \alpha  [\![v_{\hh}]\!] ( ( G_{\hh}^e)^+  -G^e)+ (1-\alpha)  [\![v_{\hh}]\!] ( ( G_{\hh}^e)^-  -G^e) 
      +\mathcal{D} [\![v_{\hh}]\!] [\![U^e_{\hh}]\!]\right\}| \\ \le&
      \alpha \sum_f | [\![v_{\hh}]\!] ( ( G_{\hh}^e)^+  -G^e)| + (1-\alpha) \sum_f | [\![v_{\hh}]\!] ( ( G_{\hh}^e)^-  -G^e)| 
      +\sum_f | |\mathcal{D} [\![v_{\hh}]\!] [\![U^e_{\hh}]\!]|  .
      \end{split}
      \end{equation*} 
      Now by standard approximation arguments \cite{cr72}, we can claim that  on each side of a given face the projection
      of the exact quantiles  provides en error of $\mathcal{O}(\hh^{p+1})$ so that the jumps are also of the same order.
      Since the number of faces is (in 1d) of order $\hh^{-1}$, this  readily allows to show that 
       \begin{equation}
  |\mathrm{IV}| \le \tilde C_{\mathrm{IV}} (\Omega,\|v\|) \,\hh^{p+1} \hh^{\min(1,p)} \,.\;\;
       \end{equation}
      }
To bound  term $\mathrm{III}$, we first consider the estimate of the quadrature rest.  Under the the smoothness hypotheses
made, for a finite element approximation of degree $p$ in one dimension 
we can use the exactness of the Gauss-Lobatto formulas for polynomials of degree $2p -1$. To estimate the 
integration remainder, one can e.g. introduce a  local  truncated Taylor series of the  exact integrand and consider the first term
not integrated exactly which can be bounded by an $\mathcal{O}( \hh^{2p})$. This leads 
 on each element   to an integration  rest of order  $\mathcal{R}_{K}^{\mathsf{GL}}=\mathcal{O}(\hh\times \hh^{2p})=\mathcal{O}( \hh^{2p+1})$. 
 Since the number of elements  in one space dimension  are of $\mathcal{O}(\hh^{-1}) $ we deduce that 
  $   \sum_K \mathcal{R}_{K}^{\mathsf{GL}}\le C_{\mathsf{GL}} \hh^{2p}$.  
  
  To estimate $(R^e_{\hh} - R^e)(\bar x)$ for a given point
   $\bar x$, we    introduce the set $\mathcal{K}_{\bar x}$ of elements such that $x < \bar x $ within the element,
   and we  denote by $K_{\bar x}$ the element containing the point. We can thus write
  \begin{equation}
R^e_{\hh} - R^e =\sum_{K\in\mathcal{K}_x} (\int_K S_h -  \int_K S(U^e,x))  
+ \int_{x_{K_{\bar x}}^L}^{\bar x} S_h -  \int_{x_{K_{\bar x}}^L}^{\bar x} S(U^e,x)
  \end{equation}
 with all integrals computed exactly. As before for all the elements in  $\mathcal{K}_{\bar x}$ we can use the properties
 of the  Gauss-Lobatto formulas, and bound the first term by a $\mathcal{O}( \hh^{2p+1})$. For the second term
 we can use the approximation error to estimate $S_{\hh}-S$, and we end with
  \begin{equation}
|R^e_{\hh} - R^e| \le C_1^S \hh^{2p} + C^S_2 \hh^{p+2}  = \bar C^S \min(\hh^{2p},\hh^{p+2}).
  \end{equation}
  We are left with the estimate of the cell correction term, evaluated using sampled values of the exact solution. 
  We start by re-writing the coefficient $\alpha_K$, re-writing the second in  \eqref{eta-corr3},  \eqref{eta-corr5},
  and \eqref{eta-corr7} which under the current hypotheses become (cf. also  equation \eqref{eta-corr14})
  \begin{equation}
\Phi_{\eta}^K(U^e_{\hh}) = \int_K (W^e_{\hh})^T\partial_x G_{\hh}\,,\;\;
\Psi_{\eta}^K(U^e_{\hh}) = \int_K  \partial_xF_{\eta}(U^e_{\hh}) \,.
  \end{equation}
  Noting that 
  \begin{equation}
\int\limits_{K}\partial_t \eta_{\hh} = \sum_{i=0,p}w_i\dfrac{d\eta_i}{dt}=\sum_{i=0,p}w_iW_i^T\dfrac{dU_i}{dt},
  \end{equation}
  we can readily write
    \begin{equation}
\Psi_{\eta}^K(U^e_{\hh}) -\Phi_{\eta}^K(U^e_{\hh})  = \int_K\left[(
\partial_t \eta_{\hh}(U^e_{\hh}) +\partial_xF_{\eta}(U^e_{\hh}) ) 
- (W^e_{\hh})^T (\partial_tU^e_{\hh} + \partial_x G_{\hh}(U^e_{\hh}) 
)\right] + \mathcal{O}(\hh^{2p+1}),
  \end{equation}
where, following the previous reasoning, the last term is the Gauss-Lobatto quadrature error on the  $ (W^e_{\hh})^T (\partial_tU^e_{\hh}$ term.
The remaining terms can be easily estimate  based on standard approximation arguments following \cite{AR:17}:
  \begin{equation}
  \int_K(
\partial_t \eta_{\hh}(U^e_{\hh}) +\partial_xF_{\eta}(U^e_{\hh}) )=
\int_K(
\underbrace{\partial_t \eta_{\hh}(U^e_{\hh})-\partial_t\eta^e}_{\mathcal{O}(\hh^{p+1})} +\underbrace{\partial_x(F_{\eta}(U^e_{\hh}) -F_{\eta}^e)}_{\mathcal{O}(\hh^{p})} ) =\mathcal{O}(\hh^{p+1})
    \end{equation}
    and similarly 
      \begin{equation}
  \int_K
(W^e_{\hh})^T (\partial_tU^e_{\hh} + \partial_x G_{\hh}(U^e_{\hh})  )=
\int_K(W^e_{\hh})^T 
\Big( 
\underbrace{\partial_t  U^e_{\hh}-\partial_tU^e}_{\mathcal{O}(\hh^{p+1})} +
\underbrace{\partial_x(G_{\hh}(U^e_{\hh}) -G^e)}_{\mathcal{O}(\hh^{p})} 
\Big) =\mathcal{O}(\hh^{p+1}).
    \end{equation}
Using the boundedness of $\partial_xv_{\hh}$, of $\partial_xU^e_{\hh}$, and of $\partial_xW^e_{\hh}$, we can  
estimate
  \begin{equation}
 \alpha_K(U^e_{\hh})=  \dfrac{\Psi_{\eta}^K(U^e_{\hh}) -\Phi_{\eta}^K(U^e_{\hh}) }{\|\partial_xW_{\hh}^e\|}  = \mathcal{O}(\hh^{p+1})
  \end{equation}
  and 
    \begin{equation}
 \alpha_K(U^e_{\hh})\int\limits_K\partial_x v_{\hh} A_0(U^e_{\hh})\partial_x U^e_{\hh} = \mathcal{O}(\hh^{p+2}).
  \end{equation}
\revPO{  Note that the de-singularization of the value of $ \alpha_K$ mentioned in remark 7 guarantees that these scaling are not
  violated close to regions in which the solution is constant.}
Using again the fact that in one dimension the number of elements is of $\mathcal{O}(\hh^{-1})$, we conclude that 
    \begin{equation}
   | \sum_K \alpha_K(U^e_{\hh})\int\limits_K\partial_x v_{\hh} A_0(U^e_{\hh})\partial_x U^e_{\hh}|,
   = \mathcal{O}(\hh^{p+1})
  \end{equation}
  which used with the previous estimates leads to the desired result.
%  \int_{\Omega}\ v ( \dfrac{d(U-U_e)}{dt} + S - S_e )  - \int_{\Omega}\nabla v (F-F_e)  + \sum_f \int_f (\hat F - F_e) + 
%\sum\limits_K \int\limits_K \alpha_K(U,U_e) \nabla ( v - \bar v)^T  A_0\nabla W(U)\\
%\Psi_{\eta}^K-\Phi_{\eta}^K = 
%\oint_K\hat F_{\eta} - 
%\sum\limits_{i\in K} W_i^T 
%\left\{ \int_{K } \phi_i \partial_x G            + [\phi_i(\hat G - G_e + G_e - G )]_L
%+ [\phi_i(\hat G - G_e + G_e - G )]_R\right\}\\
%\quad =\int_{K }\partial_t \eta + \oint_K\hat F_{\eta} 
%- \sum\limits_{i\in K} W_i^T \left\{   \int_{K } \phi_i (\dfrac{dU}{dt}+ \partial_x G     )       + [\phi_i(\hat G - G_e + G_e - G )]_L
%+ [\phi_i(\hat G - G_e + G_e - G )]_R\right\}
	
 \section{Perturbation tests initialization}
\begin{figure} 
\hspace{-0.2cm}\subfigure[WB discrete steady state]{\includegraphics[width=0.36\textwidth]{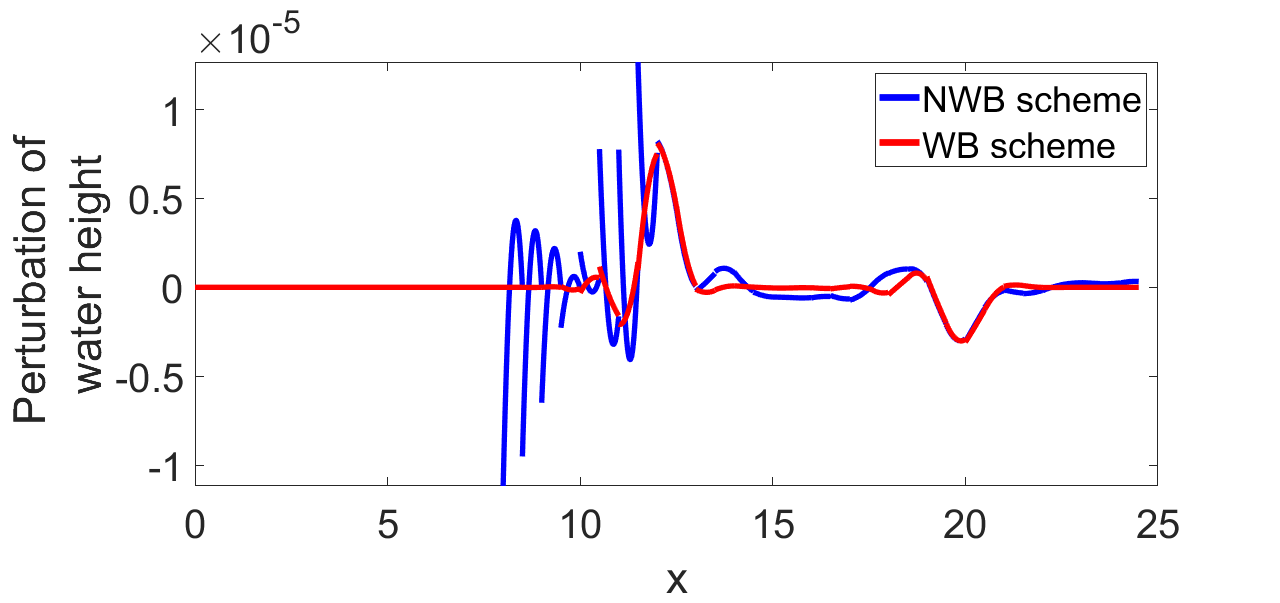}}\hspace{-0.4cm}	
	\subfigure[NWB discrete steady state]{\includegraphics[width=0.36\textwidth]{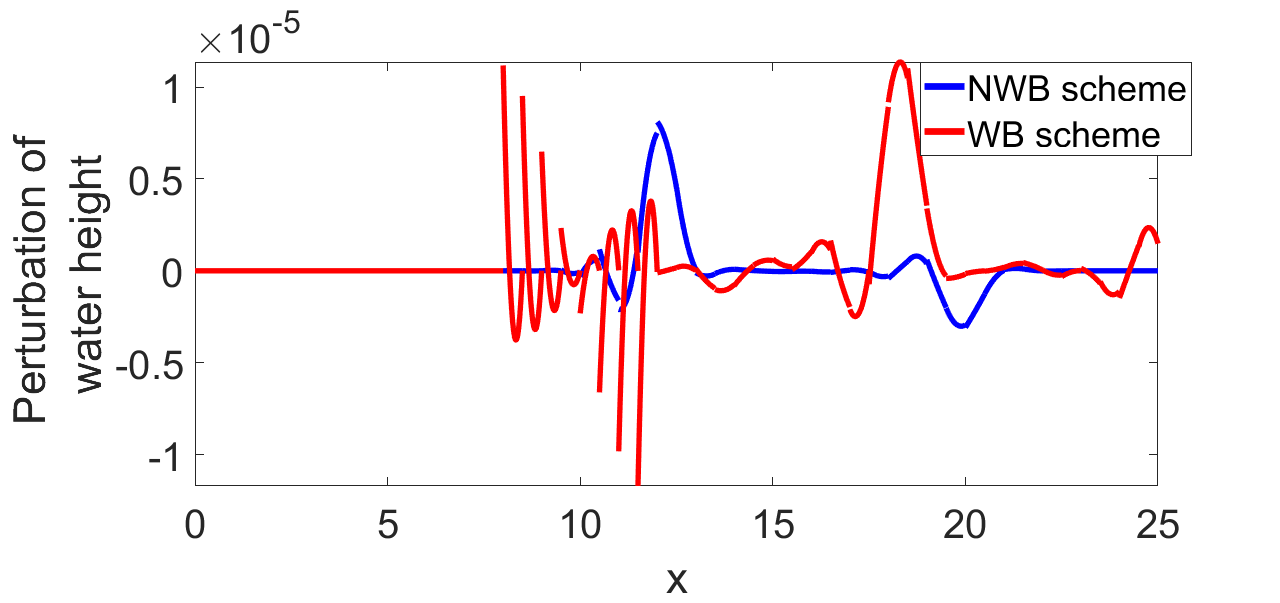}}\hspace{-0.4cm}	
	\subfigure[Analytical]{\includegraphics[width=0.36\textwidth]{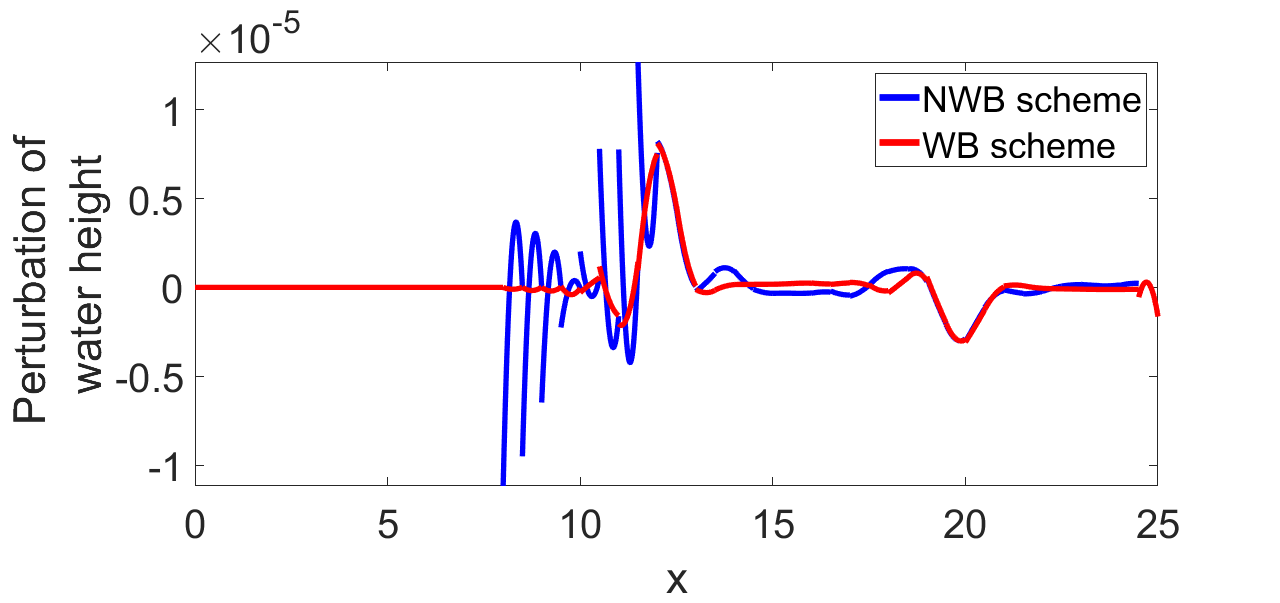}}
	\caption{Frictionless 1d super-critical equilibrium. Perturbation evolution    for the  NWB  (blue)  and WB   (red) scheme   for $p=2$. Initial condition :
	(a)  analytical steady state; (b)  ODE integrator of Proposition 3; (c)  steady state of the NWB scheme}
	%=======
	\label{fig:different_steady_states}
\end{figure}

We consider here an example to show quantitatively why we think the exact equilibrium is to be used in perturbation tests,
and in particular to define $h^*(x) $ and $u^*(x)$ in \eqref{perturbation}.
We consider again the perturbation for the test of section \S7.2.3  with   three different initializations:
the discrete solution of the RK-LobattoIIIA collocation method of proposition 3 (corresponding to the discrete steady state of the well-balanced scheme); 
the
discrete steady state of the non-well-balanced standard DGSEM scheme not using global flux quadrature; the analytical  steady state. 
%We would also like to highlight that the steady state used in the above examples is the analytical steady state given by constant $(q,E)$. The discrete steady state of the well-balanced scheme given by constant global fluxes, could be used and gives similar results. However, one could also calculate a discrete steady state for the non-well-balanced scheme and in this case a perturbation of the non-well-balanced state is better captured by the non-well-balanced scheme. Hence for the examples, we stick to using an analytical steady state, since any of the schemes would be able to better capture a perturbation to a discrete steady state of that scheme. I
In figure \ref{fig:different_steady_states}, we show a comparison of the WB and NWB schemes for the three different initial  states with a very small perturbation $\xi=10^{-5}$. As we can see, initializing with the discrete steady state of one of the scheme favours  the scheme in question. We can hardly say from the first two pictures in figure \ref{fig:different_steady_states} which scheme is better. Only from the perturbation of  the analytical initial state in the rightmost picture in the figure we can realize how closely the well-balanced scheme
reproduces the exact equilibrium.  

%
%
%Concerning the visualization, we report hereafter the solution of the same problem in which
% the plotting 
%\revp{Note that in all the plots in this paper, we have plotted the polynomial approximation of the solution in every cell, instead of an interpolation over the cells to represent the solution. We feel that interpolating over all the cells does not provide a complete picture of the solution for the perturbation as some of the jumps on the cell boundaries are smoothened out. As an example, we have shown in figure \ref{fig:linear_interpolation_subcritical} the perturbation at time $T=1.5$ with initial perturbation $\xi=10^{-3}$ plotted with a piecewise linear curve by averaging over the cell boundaries, as compared to the polynomial approximation in the cell shown in figure \ref{fig:subcritical_pert}.  }
%\begin{figure}[H]
%	%<<<<<<< HEAD
%	\hspace{-0.2cm}\includegraphics[width=0.35\textwidth]{oneD/Bathymetry/super_pert_5_s}\hspace{-0.2cm}
%	\includegraphics[width=0.35\textwidth]{oneD/Bathymetry/diffheight_supercrit_gf_pert5}\hspace{-0.2cm}
%	\includegraphics[width=0.35\textwidth]{oneD/Bathymetry/super_pert_5}
%	\caption{Perturbation of frictionless 1d sub-critical equilibrium with $\xi=10^{-3}$ using piecewise linear plots with average value on cell-boundaries. }
%	\label{fig:linear_interpolation_subcritical}
%\end{figure}
%

\section*{Acknowledgments}

\revs{M. Ricchiuto is a member of the CARDAMOM team,  INRIA and University of Bordeaux research center.
P.Ö. is thankful for the support received from the Gutenberg Research College, JGU Mainz.
Valuable discussions with  G. Russo and S. Boscarino   on the accuracy of Gauss-Lobatto collocation methods are warmly acnkowledged.}

\bibliographystyle{siam}
\bibliography{references}

\end{document}